\documentclass[12pt,leqno]{article}

\usepackage{makeidx}
\makeindex

\def\Span{\mathop{\rm span}}
\def\tr{\mathop{\rm tr}}

\newtheorem{theorem}{Theorem}
\newtheorem{lemma}[theorem]{Lemma}
\newtheorem{proposition}[theorem]{Proposition}
\newtheorem{sublemma}[theorem]{Sublemma}
\newtheorem{definition}[theorem]{Definition}
\newtheorem{corollary}[theorem]{Corollary}
\newtheorem{problem}[theorem]{Problem}
\newtheorem{remark}[theorem]{Remark}
\newtheorem{claim}[theorem]{Claim}
\newtheorem{assumptions}[theorem]{Assumptions}
\newtheorem{sassumptions}[theorem]{Standing Assumptions}
\newtheorem{examples}[theorem]{Examples}
\newtheorem{notation}[theorem]{Notation}

\newcommand{\begintheorem}{\addtocounter{equation}{1}\begin{theorem}}
\newcommand{\beginlemma}{\addtocounter{equation}{1}\begin{lemma}}
\newcommand{\beginproposition}{\addtocounter{equation}{1}\begin{proposition}}
\newcommand{\beginsublemma}{\addtocounter{equation}{1}\begin{sublemma}}
\newcommand{\begindefinition}{\addtocounter{equation}{1}\begin{definition}}
\newcommand{\begincorollary}{\addtocounter{equation}{1}\begin{corollary}}
\newcommand{\beginproblem}{\addtocounter{equation}{1}\begin{problem}}
\newcommand{\beginremark}{\addtocounter{equation}{1}\begin{remark}}
\newcommand{\beginclaim}{\addtocounter{equation}{1}\begin{claim}}
\newcommand{\beginassumptions}{\addtocounter{equation}{1}\begin{assumptions}}
\newcommand{\beginsassumptions}{\addtocounter{equation}{1}\begin{sassumptions}}
\newcommand{\beginexamples}{\addtocounter{equation}{1}\begin{examples}}
\newcommand{\beginnotation}{\addtocounter{equation}{1}\begin{notation}}

\begin{document}

\title{Some topics pertaining to algebras of linear operators}

\author{Stephen Semmes}

\date{}

\maketitle

\begin{abstract}
Here we consider finite-dimensional vector spaces and linear operators
on them.  Some of the basic ingredients will be motivated by
$*$-algebras of linear transformations on inner product spaces, as in
\cite{Arveson, Douglas}.  In particular, let us note that in the
theory of $C^*$ and Von Neumann algebras, the ``double commutant'' of
a $*$-algebra is a kind of basic closure of the algebra.

	A special case concerns algebras of operators generated by
representations of finite groups \cite{Burrow, Curtis2, Serre-rep}.
The use of characters is included in each of these three books, and is
generally avoided in the present article.

	Part of the point of view about the $p$-adic absolute value
and $p$-adic numbers that we follow here is that they can be
interesting in connection with theory of computation and related
settings just because of their nice properties and simplicity, aside
from more traditional number-theoretic uses, even if those are also
very interesting.
\end{abstract}

\tableofcontents

\section{Preliminaries}
\label{Preliminaries}
\setcounter{equation}{0}

	As usual, the real and complex numbers are denoted ${\bf
R}$\index{$R$@${\bf R}$} and ${\bf C}$,\index{$C$@${\bf C}$}
respectively.

\subsection{Finite groups}
\label{Finite groups}

	Normally groups in this article will be finite.

	A finite group\index{groups} is a finite set $G$ with a choice
of identity element $e$\index{$E$@$e$ (identity element of a group)} in
$G$ and a binary operation $G \times G \to G$, $(x,y) \mapsto xy$,
such that $xe = ex = x$ for all $x \in G$, every element $x$ of $G$
has an \emph{inverse} $x^{-1}$, which satisfies
\begin{equation}
	x x^{-1} = x^{-1} x = e,
\end{equation}
and the group operation is associative, so that
\begin{equation}
	x (y z) = (x y) z
\end{equation}
for all $x$, $y$, and $z$ in $G$.

	Two elements $x$, $y$ of $G$ are said to
\emph{commute}\index{commuting elements of a group} if
\begin{equation}
\label{x, y commute}
	xy = yx.
\end{equation}
The group $G$ is said to be \emph{commutative},\index{commutative
groups} or \emph{abelian},\index{abelian groups} if every pair of
elements in $G$ commutes.

	A subset $H$ of $G$ is a \emph{subgroup}\index{subgroups} if
$e \in H$, and if $x^{-1} \in H$ and $xy \in H$ whenever $x, y \in H$.

	If $G$ and $L$ are groups and $f : G \to L$ is a mapping
between them, then $f$ is said to be a
\emph{homomorphism}\index{homomorphisms between groups} if
\begin{equation}
	f(xy) = f(x) \, f(y)
\end{equation}
for all $x$, $y$ in $G$.  Note that $f$ automatically maps
the identity element of $G$ to the identity element of $L$ in this
case, and that 
\begin{equation}
	f(x^{-1}) = f(x)^{-1}
\end{equation}
for all $x$ in $G$.  An \emph{isomorphism}\index{isomorphism between
groups} between two groups is a homomorphism that is one-to-one and
maps the first group onto the second one.  In other words, the
homomorphism should be invertible as a mapping, and one can check that
the inverse mapping is a homomorphism too.

	If $f : G \to L$ is a homomorphism between groups, then
the \emph{kernel of $f$}\index{kernel of a homomorphism between groups}
is the subset $K$ of $G$ consisting of the elements $x$ of $G$ such that
$f(x)$ is the identity element of $L$.  It is easy to see that $K$
is a subgroup of $G$.

	A subgroup $N$ of $G$ is said to be \emph{normal}\index{normal
subgroups} if $x y x^{-1} \in N$ whenever $y \in N$ and $x \in G$.
The kernel of a homomorphism between groups is clearly a normal
subgroup.  Conversely, if $N$ is a normal subgroup of $G$, then $N$ is
the kernel of a homomorhism.  This follows from the standard quotient
subgroup construction.

\subsection{Vector spaces}
\label{Vector spaces}
\index{vector spaces}

	We make the convention that a \emph{vector space} means a
vector space defined over some field $k$ which has positive finite
dimension (unless the contrary is indicated).  If $V$ and $W$ are
vector spaces with the same scalar field $k$, then we let
$\mathcal{L}(V,W)$\index{$L(V,W)$@$\mathcal{L}(V,W)$} denote the set
of linear mappings from $V$ to $W$.  Thus $\mathcal{L}(V,W)$ is also a
vector space over $k$, using ordinary addition of linear operators and
multiplication of them by scalars.  The dimension of
$\mathcal{L}(V,W)$ is equal to the product of the dimensions of $V$
and $W$.

	For $V = W$ we simply write
$\mathcal{L}(V)$\index{$L(V)$@$\mathcal{L}(V)$} for
$\mathcal{L}(V,V)$.  For each $L_1, L_2 \in \mathcal{L}(V)$, the
composition $L_1 \circ L_2$, $(L_1 \circ L_2)(v) = L_1(L_2(v))$ for
all $v \in V$, is also a linear mapping on $V$, and hence an element
of $\mathcal{L}(V)$.  We sometimes write this simply as $L_1 \, L_2$.
This defines a product operation on $\mathcal{L}(V)$, which is not
commutative in general, but which does satisfy the associative law and
distributive laws relative to addition and scalar multiplication.  The
identity mapping $I$ or $I_V$ on $V$, which takes every element of $V$
to itself, is linear and is the identity element of $\mathcal{L}(V)$
with respect to this product.

	Let us write $GL(V)$\index{$GL(V)$} for the set of linear
transformations on $G$ which are invertible.  (If a linear mapping is
invertible simply as a mapping, then the inverse mapping is
automatically linear as well.)  Thus $GL(V)$ is a group which is
infinite.

	Suppose that $v_1, \ldots, v_n$ is a basis for the vector
space $V$ with scalar field $k$.  If $T$ is any linear transformation
on $V$, then there is a unique $n \times n$ matrix $(t_{j,k})$ with
entries in $k$ such that
\begin{equation}
\label{matrix (t_{j,k}) associated to T on V}
	T\Bigl(\sum_{l=1}^n \alpha_l v_l\Bigr) 
		= \sum_{j=1}^n \sum_{k=1}^n t_{j,k} \, \alpha_k \, v_j
\end{equation} 
for all $\alpha_1, \ldots, \alpha_n$ in $k$.  Conversely, if
$(t_{j,k})$ is an $n \times n$ matrix with entries in $k$, then
(\ref{matrix (t_{j,k}) associated to T on V}) defines a linear
transformation $T$ on $V$.

	As usual, if $S$ and $T$ are linear transformations on $V$
corresponding to $n \times n$ matrices $(s_{j,k})$ and $(t_{l,m})$
(with respect to the basis $v_1, \ldots, v_n$), then $S + T$
corresponds to the sum of the matrices, given by $(s_{j,k} +
t_{j,k})$.  If $\gamma$ is any scalar, then the matrix of $\gamma \,
T$ is $(\gamma \, t_{j,k})$.  The composition $S \circ T$ has matrix
$(u_{j,m})$ given by the classical matrix product of $(s_{j,k})$ and
$(t_{l,m})$, namely,
\begin{equation}
	u_{j,m} = \sum_{k = 1}^n s_{j,k} \, t_{k,m}.
\end{equation}
The matrix associated to the identity transformation $I$ is
$(\delta_{j,k})$, where $\delta_{j,k} = 1$ when $j = k$ and
$\delta_{j,k} = 0$ when $j \ne k$.  It is easy to check directly that
the matrix product of $(\delta_{j,k})$ with another matrix (in either
order) is always equal to the other matrix.

\subsection{Representations of finite groups}
\label{Representations of finite groups}

	Let $G$ be a group and let $V$ be a vector space.  A
\emph{representation}\index{representation of a group} of $G$ on $V$
is a mapping $\rho$ from $G$ to invertible linear transformations on
$V$ such that
\begin{equation}
\label{composition rule for representations}
	\rho_{xy} = \rho_x \circ \rho_y
\end{equation}
for all $x$ and $y$ in $G$.  Here we use $\rho_x$ to denote the
invertible linear transformation on $V$ associated to $x$ in $G$, so
that we may write $\rho_x(v)$ for the image of a vector $v \in V$
under $\rho_x$.  As a result of (\ref{composition rule for
representations}), we have that
\begin{equation}
	\rho_e = I,
\end{equation}
where $I$ denotes the identity transformation on $V$, and
\begin{equation}
	\rho_{x^{-1}} = (\rho_x)^{-1}
\end{equation}
for all $x$ in $G$.

	In other words, a representation of $G$ on $V$ is a
homomorphism from $G$ into $GL(V)$.  The dimension of $V$ is called
the \emph{degree}\index{degree of a representation} of the representation.

	Basic examples of representations are the \emph{left regular
representation}\index{left regular representation} and \emph{right
regular representation}\index{right regular representation} over a
field $k$, defined as follows.  We take $V$ to be the vector space of
functions on $G$ with values in $k$.  For the left regular
representation, we define $L_x : V \to V$ for each $x$ in $G$ by
\begin{equation}
\label{def of L_x}
	L_x(f)(z) = f(x^{-1} \, z)
\end{equation}
for each function $f(z)$ in $V$.  For the right regular representation,
we define $R_x : V \to V$ for each $x$ in $G$ by
\begin{equation}
\label{def of R_x}
	R_x(f)(z) = f(z \, x)
\end{equation}
for each function $f(z)$ in $V$.  Thus if $x$ and $y$ are elements of
$G$, then
\begin{eqnarray}
	(L_x \circ L_y)(f)(z) & = & L_x(L_y(f))(z) 
		 =  (L_y(f))(x^{-1} z)				    \\
		& = & f(y^{-1} x^{-1} z)
			= f((xy)^{-1} z) = L_{xy}(f)(z),	\nonumber
\end{eqnarray}
and 
\begin{eqnarray}
	(R_x \circ R_y)(f)(z) = R_x(R_y(f))(z) 
		& = & (R_y(f))(zx)				   \\
		& = & f(z x y) = R_{xy}(f)(z).			\nonumber
\end{eqnarray}

	Another description of these representations which can be convenient
is the following.  For each $w$ in $G$, define the function $\phi_w(z)$ on
$G$ by
\begin{equation}
	\phi_w(z) = 1 \quad\hbox{when } z = w, \qquad
		\phi_w(z) = 0 \quad\hbox{when } z \ne w.
\end{equation}
Thus the functions $\phi_w$ for $w$ in $G$ form a basis for the space of
functions on $G$.  One can check that
\begin{equation}
\label{L_x(phi_w) = phi_{xw}, R_x(phi_w) = phi_{w x^{-1}}}
	L_x(\phi_w) = \phi_{xw},  \qquad R_x(\phi_w) = \phi_{w x^{-1}}
\end{equation}
for all $x$ in $G$.

	Observe that
\begin{equation}
\label{L_x circ R_y = R_y circ L_x}
	L_x \circ R_y = R_y \circ L_x
\end{equation}
for all $x$ and $y$ in $G$.

	More generally, suppose that we have a homomorphism from
the group $G$ to the group of permutations on a nonempty finite set $E$.
That is, suppose that for each $x$ in $G$ we have a permutation $\pi_x$
on $E$, i.e., a one-to-one mapping from $E$ onto $E$, such that
\begin{equation}
	\pi_x \circ \pi_y = \pi_{xy}.
\end{equation}
As usual, this implies that $\pi_e$ is the identity mapping on $E$,
and that $\pi_{x^{-1}}$ is the inverse mapping of $\pi_x$ on $E$.  Let
$V$ be the vector space of $k$-valued functions on $E$.  Then we get a
representation of $G$ on $V$ by associating to each $x$ in $G$ the
linear mapping $\Pi_x : V \to V$ defined by
\begin{equation}
	\Pi_x(f)(a) = f(\pi_{x^{-1}}(a))
\end{equation}
for every function $f(a)$ in $V$.  This is called the
\emph{permutation representation}\index{permutation representation}
corresponding to the homomorphism $x \mapsto \pi_x$ from $G$ to
permutations on $E$.  It is indeed a representation, because for each
$x$ and $y$ in $G$ and each function $f(a)$ in $V$ we have that
\begin{eqnarray}
	(\Pi_x \circ \Pi_y)(f)(a) & = & \Pi_x(\Pi_y(f))(a)
		=  (\Pi_y(f))(\pi_{x^{-1}}(a))			   \\
		& = & f(\pi_{y^{-1}}(\pi_{x^{-1}}(a)))
			= f(\pi_{(xy)^{-1}}(a)).		\nonumber
\end{eqnarray}
Alternatively, for each $b \in E$ one can define $\psi_b(a)$ to be
the function on $E$ defined by
\begin{equation}
	\psi_b(a) = 1 \quad\hbox{when } a = b, \qquad
		\psi_b(a) = 0 \quad\hbox{when } a \ne b.
\end{equation}
Then the collection of functions $\psi_b$ for $b \in E$ is
a basis for $V$, and
\begin{equation}
	\Pi_x(\psi_b) = \psi_{\pi_x(b)}
\end{equation}
for all $x$ in $G$ and $b$ in $E$.

	Suppose that $V_1$ and $V_2$ are vector spaces over the
same field $k$, and that $T$ is a linear isomorphism from $V_1$ onto
$V_2$.  Assume also that $\rho^1$ and $\rho^2$ are representations of
a group $G$ on $V_1$ and $V_2$, respectively.  If
\begin{equation}
	T \circ \rho^1_x = \rho^2_x \circ T
\end{equation}
for all $x$ in $G$, then we say that $T$ determines an
\emph{isomorphism}\index{isomorphism between group representations}
between the representations $\rho^1$ and $\rho^2$.  We may also say
that $\rho^1$ and $\rho^2$ are \emph{isomorphic}\index{isomorphic
group representations} group representations, without referring to $T$
specifically.  Of course isomorphic representations have equal
degrees, but the converse is not true in general.

	For example, let $V_1 = V_2$ be the vector space of $k$-valued
functions on $G$, and define $T$ on $V_1 = V_2$ by $T(f)(a) =
f(a^{-1})$.  This is a one-to-one linear mapping from the space of
$k$-valued functions on $G$ onto itself, and
\begin{equation}
	T \circ R_x = L_x \circ T
\end{equation}
for all $x$ in $G$.  For if $f(a)$ is a function on $G$, then
\begin{eqnarray}
	(T \circ R_x)(f)(a) & = & T(R_x(f))(a) = R_x(f)(a^{-1})
							\\
		& = & f(a^{-1} x) = T(f)(x^{-1} a) 	\nonumber \\
		& = & L_x(T(f))(a) = (L_x \circ T)(f)(a).
							\nonumber
\end{eqnarray}
Therefore the left and right regular representations of $G$ are
isomorphic to each other, in either the real or complex case.

	Now suppose that $G$ is a group and $\rho$ is a representation
of $G$ on a vector space $V$ over the field $k$, and that $v_1,
\ldots, v_n$ is a basis of $V$.  For each $x$ in $G$ we can associate
to $\rho_x$ an $n \times n$ invertible matrix with entries in $k$
using this basis, as in Subsection \ref{Vector spaces}.  We denote this
matrix by $M_x$.  The composition rule (\ref{composition rule for
representations}) can be rewritten as
\begin{equation}
\label{composition rule in terms of matrices}
	M_{xy} = M_x \, M_y,
\end{equation}
where the matrix product is used on the right side of the equation.

	A different choice of basis for $V$ will lead to a different
mapping $x \mapsto N_x$ from $G$ to invertible $n \times n$ matrices.
However, the two mappings $x \mapsto M_x$, $x \mapsto N_x$ will be
\emph{similar},\index{similar (mappings to matrices)} in the sense
that there is an invertible $n \times n$ matrix $S$ with entries in
$k$ such that
\begin{equation}
	N_x = S \, M_x \, S^{-1}
\end{equation}
for all $x$ in $G$.  That is, $S$ is the matrix corresponding to the
change of basis. 

	On the other hand, we can start with a mapping $x \mapsto M_x$
from $G$ into invertible $n \times n$ matrices with entries in $k$
which satisfies (\ref{composition rule in terms of matrices}), and
convert it into a representation on an $n$-dimensional vector space
over $k$ by reversing the process.  That is, one chooses a basis for
the vector space, and then converts $n \times n$ matrices into linear
transformations on the vector space using the basis to get the
representation.  One might as well take the vector space to be $k^n$,
the space of $n$-tuples with coordinates in $k$ (using coordinatewise
addition and scalar multiplication), and the basis $v_1, \ldots, v_n$
to be the standard basis, where $v_i$ has $i$th coordinate $1$ and all
others $0$.  In this way one gets the usual correspondence between $n
\times n$ matrices and linear transformations on $k^n$.

	If one has two representations of $G$ on vector spaces $V_1$,
$V_2$ with the same scalar field $k$, then these two representations
are isomorphic if and only if the associated mappings from $G$ to
invertible matrices as above, using any choices of bases on $V_1$ and
$V_2$, are similar, with the similarity matrix $S$ having entries in
$k$.

\subsection{Reducibility}
\label{Reducibility}

	Let $G$ be a finite group, $V$ a vector space over a field
$k$, and $\rho$ a representation of $G$ on $V$.  Suppose that there is
a vector subspace $W$ of $V$ such that
\begin{equation}
\label{rho_x(W) subseteq W}
	\rho_x(W) \subseteq W
\end{equation}
for all $x$ in $G$.  This is equivalent to saying that
\begin{equation}
	\rho_x(W) = W
\end{equation}
for all $x$ in $G$, since one could apply (\ref{rho_x(W) subseteq W})
also to $\rho_{x^{-1}} = (\rho_x)^{-1}$.  We say that $W$ is
\emph{invariant}\index{invariant subspaces (of a group
representation)} or \emph{stable} under the representation $\rho$.

	Recall that a subspace $Z$ of $V$ is said to be a
\emph{complement}\index{complementary subspaces of a vector space}
of a subspace $W$ if
\begin{equation}
\label{complementary conditions, 1}
	W \cap Z = \{0\}, \quad W + Z = V.
\end{equation}
Here $W + Z$ denotes the span of $W$ and $Z$, which is the subspace of
$V$ consisting of vectors of the form $w + z$, $w \in W$, $z \in Z$.
The conditions (\ref{complementary conditions, 1}) are equivalent to
saying that
\begin{eqnarray}
\label{complementary conditions, 2}
	&& \hbox{every vector $v \in V$ can be written in a}		\\
	&& \hbox{unique way as $w + z$, $w \in W$, $z \in Z$.}
								\nonumber
\end{eqnarray}
Complementary subspaces always exist, because a basis for a vector
subspace of a vector space can be enlarged to a basis of the whole
vector space.

	If $W$, $Z$ are complementary subspace of a vector space $V$,
then we get a linear mapping $P$ on $V$ which is the \emph{projection
of $V$ onto $W$ along $Z$}\index{projection of a vector space onto
a subspace along another subspace} and which is defined by
\begin{equation}
	P(w+z) = w \quad\hbox{for all } w \in W, \ z \in Z.
\end{equation}
Thus $I - P$ is the projection of $V$ onto $Z$ along $W$, where $I$
denotes the identity transformation on $V$.

	Note that $P^2 = P$, when $P$ is a projection.  Conversely,
if $P$ is a linear operator on $V$ such that $P^2 = P$, then $P$ is
the projection of $V$ onto the subspace of $V$ which is the image
of $P$ along the subspace of $V$ which is the kernel of $P$.

\subsection{Reducibility, continued}
\label{Reducibility, continued}

	Again let $G$ be a finite group, $V$ a vector space over a
field $k$, $\rho$ a representation of $G$ on $V$, and $W$ a subspace
of $V$ which is invariant under $\rho$.  In this subsection we assume
that either
\begin{equation}
\label{k has characteristic 0}
	\hbox{$k$ has characteristic $0$}
\end{equation}
or 
\begin{eqnarray}
\label{char(k) > 0, does not divide the number of elements of the group}
   && \hbox{$k$ has positive characteristic and the number of elements}  \\
   && \hbox{of $G$ is not divisible by the characteristic of $k$.}  \nonumber
\end{eqnarray}

	Let us show that there is a subspace $Z$ of $V$ such that $Z$
is a complement of $W$ and $Z$ is also invariant under the
representation $\rho$ of $G$ on $V$.  To do this, we start with any
complement $Z_0$ of $W$ in $V$, and we let $P_0 : V \to V$ be the
projection of $V$ onto $W$ along $Z_0$.  Thus $P_0$ maps $V$ to $W$,
and $P_0(w) = w$ for all $w \in W$.

	Let $m$ denote the number of elements of $G$.  Define a
linear mapping $P : V \to V$ by
\begin{equation}
	P = \frac{1}{m} \sum_{x \in G} \rho_x \circ P_0 \circ (\rho_x)^{-1}.
\end{equation}
Our assumption on $k$ implies that $1/m$ makes sense as an element of
$k$, i.e., as the multiplicative inverse of a sum of $m$ $1$'s in $k$,
where $1$ refers to the multiplicative identity element of $k$.  This
expression defines a linear mapping on $V$, because $P_0$ and the
$\rho_x$'s are.  We actually have that $P$ maps $V$ to $W$, because
$P_0$ maps $V$ to $W$, and because the $\rho_x$'s map $W$ to $W$, by
hypothesis.  If $w \in W$, then $(\rho_x)^{-1}(w) \in W$ for all $x$
in $G$, and then $P_0((\rho_x)^{-1}(w)) = (\rho_x)^{-1}(w)$.  Thus we
conclude that
\begin{equation}
	P(w) = w  \qquad\hbox{for all } w \in W,
\end{equation}
by the definition of $P$.

	The definition of $P$ also implies that
\begin{equation}
\label{rho_y circ P circ (rho_y)^{-1} = P}
	\rho_y \circ P \circ (\rho_y)^{-1} = P
\end{equation}
for all $y$ in $G$.  Indeed,
\begin{eqnarray}
	\rho_y \circ P \circ (\rho_y)^{-1} 
  & = & \frac{1}{m} \sum_{x \in G} \rho_y \circ \rho_x \circ P_0
			  \circ (\rho_x)^{-1} \circ (\rho_y)^{-1}
								\\
 & = & \frac{1}{m} \sum_{x \in G} \rho_{yx} \circ P_0 \circ (\rho_{yx})^{-1}
							\nonumber \\
 & = & \frac{1}{m} \sum_{x \in G} \rho_x \circ P_0 \circ (\rho_x)^{-1}
	= P.						\nonumber
\end{eqnarray}
This would work as well for other linear transformations instead of
$P_0$ as the initial input, but a subtlety is that in general one can
get the zero operator after taking the sum.  The remarks in the
preceding paragraph ensure that this does not happen here, except in
the degenerate situation where $W = \{0\}$.

	Because $P(V) \subseteq W$ and $P(w) = w$ for all $w \in W$,
$P$ is the projection of $V$ onto $W$ along some subspace $Z$ of $V$.
Specifically, one should take $Z$ to be the kernel of $P$.  It
is easy to see that $W \cap Z = \{0\}$, since $P(w) = w$ for all
$w \in W$.  On the other hand, if $v$ is any element of $V$, then
we can write $v$ as $P(v) + (v - P(v))$.  We have that $P(v) \in W$,
and that 
\begin{equation}
	P(v - P(v)) = P(v) - P(P(v)) = P(v) - P(v) = 0,
\end{equation}
where the second equality uses the fact that $P(v) \in W$.  Thus $v -
P(v)$ lies in $Z$, the kernel of $P$.  This shows that $W$ and $Z$
satisfy (\ref{complementary conditions, 1}), so that $Z$ is a
complement of $W$ in $V$.  The invariance of $Z$ under the
representation $\rho$ follows from (\ref{rho_y circ P circ
(rho_y)^{-1} = P}).

	Thus the representation $\rho$ of $G$ on $V$ is isomorphic to
the direct sum of the representations of $G$ on $W$ and $Z$ that are
the restrictions of $\rho$ to $W$ and $Z$.

	There can be smaller invariant subspaces within these
invariant subspaces, so that one can repeat the process.  Before
addressing this, let us introduce some terminology.  We say that
subspaces $W_1, W_2, \ldots, W_h$ of $V$ form an \emph{independent
system}\index{independent system of subspaces} if $W_j \ne \{0\}$ for
each $j$ and if $w_j \in W_j$, $1 \le j \le h$, and
\begin{equation}
	\sum_{j=1}^h w_j = 0
\end{equation}
imply
\begin{equation}
	w_j = 0, j = 1, 2, \ldots, h.
\end{equation}
If, in addition, $\Span (W_1, \ldots, W_h) = V$, then every vector $v$
in $V$ can be written in a unique way as $\sum_{j=1}^h u_j$ with $u_j
\in W_j$ for each $j$.

\begindefinition
\label{def of an irreducible representation}
\index{irreducible representations of finite groups}
Let $G$ be a finite group, $U$ be a vector space, and $\sigma$ be a
representation of $G$ on $U$.  We say that $\sigma$ is
\emph{irreducible} if there are no vector subspaces of $U$ which are
invariant under $\sigma$ except for $\{0\}$ and $U$ itself.
\end{definition}

\beginlemma
\label{decomp into irreducible pieces}
Suppose that $G$ is a finite group, $V$ is a vector space over a field
$k$ which satisfies (\ref{k has characteristic 0}) or (\ref{char(k) >
0, does not divide the number of elements of the group}), and $\rho$
is a representation of $G$ on $V$.  Then there is an independent
system of subspaces $W_1, \ldots, W_h$ of $V$ such that $\Span (W_1,
\ldots, W_h) = V$, each $W_j$ is invariant under $\rho$, and the
restriction of $\rho$ to each $W_j$ is an irreducible representation
of $G$.
\end{lemma}

	In other words, the representation $\rho$ of $G$ on $V$ is
isomorphic to a direct sum of irreducible representations of $G$
(which are restrictions of $\rho$ to subspaces of $V$).  It is not
hard to prove the lemma, by repeatedly finding invariant complements
for invariant subspaces, until one reaches subspaces to which the
restriction of $\rho$ is irreducible.  More precisely, suppose that
one has an independent system of subspaces of $V$ whose span is $V$
and for which each subspace in the system is invariant under $\rho$.
Initially this system could consist of $V$ by itself.  If the
restriction of $\rho$ to each subspace in the system is irreducible,
then we have what we want for the lemma.  Otherwise, for each subspace
where the restriction of $\rho$ is not irreducible, there is a pair of
nontrivial complementary subspaces in that subspace which are
invariant under $\rho$.  We can use these complementary subspaces in
place of the subspace in which they were found, and this leads to a
new independent system of invariant subspaces of $V$ whose span is all
of $V$.  (It is not hard to check that the new collection of subspaces
of $V$ still forms an indepedent system.)  By repeating this process,
we can get irreducibility, as in the lemma.

\beginremark 
{\rm Exercise 16.7 on p136 of \cite{Serre-rep} discusses
examples where representations over fields of positive characteristic
(dividing the order of the group) do not ``come from'' representations
over fields of characteristic $0$.
}
\end{remark}

\subsection{Positive elements and symmetric fields}
\label{Positive elements and symmetric fields}

	Let $k_0$ be a field of characteristic $0$.  We say that a
subset $A$ of $k_0$ is a \emph{set of positive elements}\index{set of
positive elements in a field} if (i) $0$ does not lie in $A$, (ii) $x
+ y$ and $x \, y$ lie in $A$ whenever $x, y \in A$, and (iii) $w^2$
lies in $A$ whenever $w$ is a nonzero element of $k_0$.  Observe that
the multiplicative identity element $1$ in $k_0$ lies in $A$ by (iii),
and that $-1$ does not lie in $A$, since otherwise one could use (ii)
to get that $0$ is in $A$.  Any sum of $1$'s lies in $A$, so that
$k_0$ must have characteristic $0$ in order for a set $A$ of positive
elements to exist.  If $x$ is an element of $A$, then $1/x$ lies in
$A$, because $x \ne 0$ by (i), $1/x^2 \in A$ by (iii), and hence $1/x
= x \, (1/x^2)$ is in $A$ by (ii).

	A field $k$ of characteristic $0$ is said to be a
\emph{symmetric field}\index{symmetric field} if the following
conditions are satisfied.  First, we ask that $k$ be equipped with an
automorphism called \emph{conjugation}\index{conjugation (on a
symmetric field)} which is an involution, i.e., the conjugate of the
conjugate of an element $x$ of $k$ is equal to $x$.  The conjugate of
$x \in k$ is denoted $\overline{x}$, and we write $k_0$ for the
subfield of $k$ consisting of elements $x$ such that $\overline{x} =
x$.  The second condition is that $k_0$ be equipped with a set $A$ of
positive elements.  (Thus $A$ and the conjugation automorphism are
part of the data for a symmetric field.)  The last condition is that
$x \, \overline{x} \in A$ for all nonzero elements $x$ of $k$.

	It may be that the conjugation automorphism is equal to the
identity, so that $k = k_0$.  In this case $k$ is a symmetric field
if it is equipped with a set of positive elements.

	Of course any subfield of the field of real numbers has a
natural set of positive elements, namely the elements that are
positive in the usual sense.  A basic class of symmetric fields are
subfields of the complex numbers which are invariant under complex
conjugation, using complex conjugation as the conjugation automorphism
on the field.  In this case the subfield of elements fized by the
conjugation is the subfield of real numbers in the field, and for the
set of positive elements we again use the elements of the field which
are positive real numbers in the usual sense.

\subsection{Inner product spaces}
\label{Inner product spaces}

	Let $k$ be a \emph{symmetric field}, and suppose that $V$ is a
vector space over $k$.  A function $\langle v, w \rangle$ on $V \times
V$ with values in $k$ is an \emph{inner product}\index{inner product}
on $V$ if it satisfies the following three conditions.  First, for
each $w \in V$, the function
\begin{equation}
	v \mapsto \langle v, w \rangle
\end{equation}
is linear.  Second,
\begin{equation}
	\langle w, v \rangle = \overline{\langle v, w \rangle}
\end{equation}
for all $v, w \in V$, where $\overline{x}$ denotes the conjugate of an
element $x$ of $k$, as in Subsection \ref{Positive elements and symmetric
fields}.  Third, if $v$ is a nonzero vector in $V$, then $\langle v, v
\rangle$ lies in the set of positive elements associated to $k$ as in
Subsection \ref{Positive elements and symmetric fields}, and is nonzero
in particular.  Of course the second condition implies that $\langle
v, v \rangle$ lies in the subfield $k_0$ of $k$ consisting of $x$ such
that $\overline{x} = x$.

	A vector space equipped with an inner product is called an
\emph{inner product space}.  It is easy to see that there are plenty
of inner products on a vector space over a symmetric field, by writing
them down explicitly using a basis.

	Suppose that $\langle \cdot, \cdot \rangle$ is an inner product
on the vector space $V$ over the symmetric field $k$.  Two vectors $u,
w \in V$ are said to be \emph{orthogonal}\index{orthogonal vectors} if
\begin{equation}
	\langle u, w \rangle = 0.
\end{equation}
A collection of vectors $v_1, \ldots, v_m$ in $V$ is said to be
\emph{orthogonal} if $v_j$ and $v_l$ are orthogonal when $j \ne l$.
It is easy to see that any orthogonal collection of nonzero vectors in
$V$ is linearly independent.  A collection $Z_1, \ldots, Z_r$ of
vector subspaces of $V$ is said to be \emph{orthogonal} if any vectors
in $Z_j$ and $Z_p$, $j \ne p$, $1 \le j, p \le r$, are orthogonal.

	Assume that $u_1, \ldots, u_m$ are nonzero orthogonal
vectors in $V$, and let $U$ denote their span.  For each $w$ in $U$,
we have that
\begin{equation}
  w = \sum_{j=1}^m \frac{\langle w, u_j \rangle}{\langle u_j, u_j \rangle}
			\, u_j.
\end{equation}
In other words, $w$ is some linear combination of the $u_j$'s, and
then the inner product and the assumption of orthogonality can be used
to determine the coefficients of the $u_j$'s, as in the preceding
formula.  Define a linear operator $P$ on $V$ by
\begin{equation}
\label{def of P}
   P(v) = \sum_{j=1}^m 
	    \frac{\langle v, u_j \rangle}{\langle u_j, u_j \rangle}
			\, u_j.
\end{equation}
Then $P(v)$ lies in $U$ for all $v$ in $V$, $P(w) = w$ for all $w$ in
$U$, and
\begin{equation}
	\langle P(v), w \rangle = \langle v, w \rangle
		\quad\hbox{for all $v \in V$ and $w \in U$.}
\end{equation}
This last condition is equivalent to saying that $v - P(v)$ is
orthogonal to every element of $U$.

\beginremark
\label{characterizing the orthogonal projection}
{\rm The operator $P : V \to V$ is characterized by the requirements
that $P(v)$ lie in $U$ and $v - P(v)$ be orthogonal to every element
of $U$ for all $v$ in $V$.  Specifically, if $v$ is a vector in $V$
and $z_1$ and $z_2$ are two elements of $U$ such that $v - z_1$ and $v
- z_2$ are orthogonal to all elements of $U$, then $z_1 = z_2$.  This
is because $z_1 - z_2$ lies in $U$, and $z_1 - z_2 = (z_1 - v) + (v -
z_2)$ is orthogonal to all elements of $U$, so that $z_1 - z_2$ is
orthogonal to itself.  As a consequence, $P$ does not depend on the
choice of orthogonal basis for $U$.}
\end{remark}

\beginlemma
\label{lin. indep. vectors to orthogonal vectors}
If $v_1, \ldots, v_m$ are linearly independent vectors in $V$, then
there are nonzero orthogonal vectors $u_1, \ldots, u_m$ in $V$ such
that the span of $v_1, \ldots, v_m$ is equal to the span of $u_1,
\ldots, u_m$.
\end{lemma}

	To prove this, one can argue by induction.  The $m = 1$ case
is trivial, and so we suppose that the statement is true for some $m$
and try to establish it for $m+1$.  Let $v_1, \ldots, v_{m+1}$ be a
set of $m+1$ linearly-independent vectors in $V$.  By the induction
hypothesis, there are nonzero orthogonal vectors $u_1, \ldots, u_m$
such that the span of $v_1, \ldots, v_m$ is equal to the span of $u_1,
\ldots, u_m$.  Let $P$ be the projection onto the span of $u_1,
\ldots, u_m$ defined above, and set $u_{m+1} = v_{m+1} - P(v_{m+1})$.
Then $u_{m+1} \ne 0$, because $v_{m+1}$ does not lie in the span of
$v_1, \ldots, v_m$, by linear independence.  From the properties of
$P$ we know that $u_{m+1}$ is orthogonal to $u_1, \ldots, u_m$, and
hence $u_1, \ldots, u_{m+1}$ is a collection of nonzero orthogonal
vectors.  It is not hard to check that the span of $u_1, \ldots,
u_{m+1}$ is equal to the span of $v_1, \ldots, v_{m+1}$, using the
corresponding statement for $m$ and the definition of $u_{m+1}$.  This
completes the proof of the lemma.

\begincorollary
\label{Every nonzero subspace of V admits an orthogonal basis}
Every nonzero subspace of $V$ admits an orthogonal basis.
\end{corollary}

	(One could say that the subspace $\{0\}$ has the empty
orthogonal basis.)

	Corollary \ref{Every nonzero subspace of V admits an
orthogonal basis} follows from Lemma \ref{lin. indep. vectors to
orthogonal vectors}, by starting with any basis for the subspace.

\begincorollary
\label{subspaces admit orthogonal projections}
If $U$ is a vector subspace of $V$, then there is a unique linear
operator $P = P_U$ on $V$ (the \emph{orthogonal projection} of $V$
onto $U$)\index{orthogonal projections} such that $P(v)$ lies in $U$
for all $v$ in $V$, $v - P(v)$ is orthogonal to all elements of $U$
for any $v$ in $V$, and $P(w) = w$ for all $w$ in $U$.
\end{corollary}

	This follows from the earlier discussion, since we know from
the previous corollary that $U$ has an orthogonal basis.    

	If $U$ is a vector subspace of $V$, then the \emph{orthogonal
complement}\index{orthogonal complement} $U^\perp$\index{$U^\perp$} of
$U$ in $V$ is the vector subspace of $V$ defined by
\begin{equation}
\label{def of W^perp}
	U^\perp = \{v \in V : \langle v, w \rangle = 0 
					\hbox{ for all } w \in U\}.
\end{equation}
Notice that
\begin{equation}
	U \cap U^\perp = \{0\}.
\end{equation}
We can reformulate the characterizing properties of the orthogonal
projection $P_U$ of $V$ onto $U$ as saying that for each vector $v$ in
$V$, $P_U(v)$ lies in $U$ and $v - P_U(v)$ lies in $U^\perp$.  Thus
$V$ is the span of $U$ and $U^\perp$.  Also, $U^\perp$ is exactly the
same as the kernel of $P_U$.

	Let us check that 
\begin{equation}
	(U^\perp)^\perp = U.
\end{equation}
for any subspace $U$ of $V$.  The inclusion $U \subseteq
(U^\perp)^\perp$ is a simple consequence of the definition, and so it
is enough to show that $(U^\perp)^\perp \subseteq U$.  Let $w \in
(U^\perp)^\perp$ be given.  We know that $P_U(w) \in U \subseteq
(U^\perp)^\perp$, and hence $w - P_U(w) \in (U^\perp)^\perp$ as well.
On the other hand, $w - P_U(w)$ lies in $U^\perp$.  Hence $w - P_U(w)$
lies in both $U^\perp$ and $(U^\perp)^\perp$, and is therefore $0$.
This shows that $w = P_U(w)$ is contained in $U$, as desired.

	Notice that
\begin{equation}
	P_{U^\perp} = I - P_U,
\end{equation}
since $I - P_U$ satisfies the properties that characterize
$P_{U^\perp}$.  Also, if $U_1$, $U_2$ are orthogonal subspaces of $V$,
then $U_1 + U_2 = V$ if and only if $U_2 = U_1^\perp$, which is
equivalent to $U_1 = U_2^\perp$.

	Suppose that $T$ is a linear operator on $V$.  The
\emph{adjoint}\index{adjoint of a linear operator}\index{$T^*$
(adjoint of an operator $T$)} of $T$ is the unique linear operator
$T^*$ on $V$ such that
\begin{equation}
	\langle T(v), w \rangle = \langle v, T^*(w) \rangle
\end{equation}
for all $v, w \in V$.  It is not hard to describe the matrix of
$T^*$ with respect to an orthogonal matrix in terms of the matrix
of $T$ with respect to the same basis.

	The adjoint of the identity operator $I$ on $V$ is itself.  If
$S$ and $T$ are linear operators on $V$, then
\begin{eqnarray}
	(T^*)^* & = & T,				\\
	(S + T)^* & = & S^* + T^*,	
\end{eqnarray}
and
\begin{equation}
	(S \circ T)^* = T^* \circ S^*.
\end{equation}
If $a$ is an element of the scalar field $k$, then
\begin{equation}
	(a \, T)^* = \overline{a} \, T^*.
\end{equation}

	A linear operator $S$ on $V$ is said to be
\emph{self-adjoint}\index{self-adjoint linear operators} if $S^* = S$.
Sometimes one says instead that $S$ is \emph{symmetric}.  As a basic
class of examples, if $U$ is a subspace of $V$ and $P_U : V \to V$ is
the orthogonal projection of $V$ onto $U$, then $P_U$ is self-adjoint.
More precisely, if $v_1$ and $v_2$ are arbitrary vectors in $V$, then
\begin{equation}
	\langle P_U(v_1), v_2 \rangle 
		= \langle P_U(v_1), P_U(v_2) \rangle
		= \langle v_1, P_U(v_2) \rangle,
\end{equation}
since $P_U(v) \in U$ and $v - P_U(v) \in U^\perp$ for $v \in V$.

	A linear operator $A$ on $V$ is said to be
\emph{antiself-adjoint},\index{antiself-adjoint linear operators} or
\emph{antisymmetric},\index{antisymmetric linear operators} if $A^* =
- A$.  Any linear operator $T$ on $V$ can be written as the sum
of a self-adjoint operator $S$ and an antiself-adjoint operator $A$
by taking $S = (T + T^*)/2$ and $A = (T - T^*)/2$.

	Suppose that $R$ is a linear operator on $V$ such that
\begin{equation}
	\langle R(v), R(w) \rangle = \langle v, w \rangle
\end{equation}
for all $v, w \in V$.  This is equivalent to saying that $R$ is
invertible and
\begin{equation}
\label{R^{-1} = R^* (preserving the inner product)}
	R^{-1} = R^*.
\end{equation}
One often says that $R$ is an an \emph{orthogonal
transformation}\index{orthogonal transformations} or a \emph{unitary
transformation}.\index{unitary transformations}

	The notion of adjoints can be extended to mappings between two
inner product spaces.  More precisely, let $(V_1, \langle \cdot, \cdot
\rangle_1)$ and $(V_2, \langle \cdot, \cdot \rangle_2)$ be two inner
product spaces, with the same symmetric field $k$ of scalars, and let
$T$ be a linear mapping from $V_1$ to $V_2$.  The
\emph{adjoint}\index{adjoint of a linear operator}\index{$T^*$
(adjoint of an operator $T$)} of $T$ is the unique linear mapping $T^*
: V_2 \to V_1$ such that
\begin{equation}
	\langle T(u), w \rangle_2 = \langle u, T^*(w) \rangle_1
\end{equation}
for all $u \in V_1$ and $w \in V_2$.  The adjoint can be described
easily in terms of orthogonal bases and matrices again, although now
one would use an orthogonal basis in each of $V_1$ and $V_2$.

	As before, $(T^*)^* = T$, $(R + T)^* = R^* + T^*$, and $(a \,
T)^* = \overline{a} \, T^*$ for any two linear operators $R, T : V_1
\to V_2$ and any scalar $a \in k$.  If $(V_3, \langle \cdot, \cdot
\rangle_3)$ is another inner product space with the same scalar field
$k$, and if $S : V_2 \to V_3$ is linear, so that $S \circ T : V_1 \to
V_3$ is defined, then $(S \circ T)^* = T^* \circ S^*$, as linear
mappings from $V_1$ to $V_3$.

	Note that a linear mapping $T : V_1 \to V_2$ preserves the
inner products on $V_1$ and $V_2$, in the sense that
\begin{equation}
	\langle T(u), T(v) \rangle_2 = \langle u, v \rangle_1
\end{equation}
for all $u, v \in V_1$, if and only if $T^* \circ T$ is the identity
operator on $V_1$.  In this case $T$ is one-to-one, but it may not map
$V_1$ onto $V_2$.  If it does, so that $T$ is invertible, then we can
simply say that $T$ preserves the inner products if $T^* = T^{-1}$.

\subsection{Inner products and representations}
\label{Inner products and representations}

	Let $G$ be a finite group, let $V$ be a vector space over a
symmetric field $k$, and let $\rho$ be a representation of $G$ on $V$.
If $\langle \cdot, \cdot \rangle$ is an inner product on $V$, then
$\langle \cdot, \cdot \rangle$ is said to be \emph{invariant under the
representation $\rho$},\index{invariant inner products}\index{inner
products!invariant under a representation} or simply
\emph{$\rho$-invariant}, if every $\rho_x : V \to V$, $x$ in $G$,
preserves the inner product, i.e., if
\begin{equation}
	\langle \rho_x(v), \rho_x(w) \rangle = \langle v, w \rangle
\end{equation}
for all $x$ in $G$ and $v$, $w$ in $V$.

	If $\langle \cdot, \cdot \rangle_0$ is any inner product on
$V$, then we can obtain an invariant inner product $\langle \cdot,
\cdot \rangle$ from it by setting
\begin{equation}
  \langle v, w \rangle = 
	\sum_{y \in G} \langle \rho_y(v), \rho_y(w) \rangle_0.
\end{equation}
It is easy to check that this does define an inner product on $V$
which is invariant under the representation $\rho$.  Notice that the
positivity condition for $\langle \cdot, \cdot \rangle$ is implied by
the positivity condition for $\langle \cdot, \cdot \rangle_0$, which
prevents $\langle \cdot, \cdot \rangle$ from reducing to $0$ in
particular.

	In some situations it is easy to write down an invariant inner
product for a representation directly.  In the case of the left and
right regular representations for a group $G$ over the symmetric field
$k$, as in Subsection \ref{Representations of finite groups}, one can use
the inner product
\begin{equation}
	\langle f_1, f_2 \rangle 
		= \sum_{x \in G} f_1(x) \, \overline{f_2(x)}.
\end{equation}
More generally, for a permutation representation of $G$ relative to a
nonempty finite set $E$, as in Subsection \ref{Representations of finite
groups}, one can use the inner product
\begin{equation}
	\langle f_1, f_2 \rangle 
		= \sum_{a \in E} f_1(a) \, \overline{f_2(a)}.
\end{equation}
This inner product is invariant under the permutation representation,
because the permutations simply rearrange the terms in the sums
without affecting the sum as a whole.

	Let $\langle \cdot, \cdot \rangle$ be any inner product on $V$
which is invariant under the representation $\rho$.  Suppose that
$W$ is a subspace of $V$ which is invariant under $\rho$, so that
\begin{equation}
	\rho_x(W) = W 
\end{equation} 
for all $x$ in $G$.  Let $W^\perp$ be the orthogonal complement of $W$
in $V$ with respect to this inner product $\langle \cdot, \cdot
\rangle$, as in (\ref{def of W^perp}).  Then 
\begin{equation}
	\rho_x(W^\perp) = W^\perp
\end{equation}
for all $x$ in $G$, since the inner product is invariant under $\rho$.
This gives another approach to finding an invariant complement to an
invariant subspace, as in Subsection \ref{Reducibility, continued}.

	One can repeat this to get invariant subspaces on which $\rho$
restricts to be irreducible, as in Subsection \ref{Reducibility,
continued}.  The next lemma is the analogue of Lemma \ref{decomp into
irreducible pieces} in this situation.

\beginlemma
\label{decomp into orthogonal irreducible pieces}
Suppose that $G$ is a finite group, $V$ is a vector space over a
symmetric field $k$, and $\rho$ is a representation of $G$ on $V$.
Assume also that $\langle \cdot, \cdot \rangle$ is an inner product on
$V$ which is invariant under $\rho$.  Then there are orthogonal
nonzero subspaces $W_1, \ldots, W_h$ of $V$ such that $\Span (W_1,
\ldots, W_h) = V$, each $W_j$ is invariant under $\rho$, and the
restriction of $\rho$ to each $W_j$ is an irreducible representation
of $G$.
\end{lemma}

\section{Algebras of linear operators}
\label{Algebras of linear operators}
\setcounter{equation}{0}

\subsection{Basic notions}
\label{Basic notions (Algebras of linear operators)}

	Let $V$ be a vector space over a field $k$, and recall that
$\mathcal{L}(V)$ denotes the collection of linear operators on
$V$.

\begindefinition
\label{def of algebras of linear operators}
\index{algebras of linear operators} 
A subset $\mathcal{A}$ of $\mathcal{L}(V)$ is said to be an
\emph{algebra of linear operators on $V$} if $\mathcal{A}$
is a vector subspace of $\mathcal{L}(V)$ which contains the
identity transformation on $V$ and which contains $S \circ T$
whenever $S, T \in \mathcal{A}$.
\end{definition}

	Thus, for example, $\mathcal{L}(V)$ is an algebra of linear
operators on $V$, as is the set of scalar multiples of the identity
transformation on $V$.

\begindefinition
\label{def of commutant}
\index{commutant (algebras of operators)}
If $\mathcal{A}$ is a subset of $\mathcal{L}(V)$, then the
\emph{commutant of $\mathcal{A}$} is defined to be the set of
operators $T \in \mathcal{L}(V)$ such that $T \circ S = S \circ T$ for
all $S \in \mathcal{A}$.  The commutant of $\mathcal{A}$ is denoted
$\mathcal{A}'$.\index{$A'$@$\mathcal{A}'$ (the commutant of
$\mathcal{A}$)}
\end{definition}

	For example, $\mathcal{L}(V)'$ is the set of scalar multiples
of the identity, and the commutant of the set of scalar multiples of
the identity is $\mathcal{L}(V)$.

	If $\mathcal{A}$ is any subset of $\mathcal{L}(V)$, then the
algebra of operators on $V$ generated by $\mathcal{A}$ is the algebra
consisting of linear combinations of the identity operator and finite
products of elements of $\mathcal{A}$.  It is easy to see that the
commutant of $\mathcal{A}$ is the same as the commutant of the algebra
of operators on $V$ generated by $\mathcal{A}$.  The commutant of any
subset $\mathcal{A}$ of $\mathcal{L}(V)$ is an algebra of operators on
$V$.

	The double commutant\index{double commutant (algebras of
operators)} of a set $\mathcal{A} \subseteq \mathcal{L}(V)$ is simply
the commutant of the commutant of $\mathcal{A}$, $(\mathcal{A}')'$.
For simplicity we also write this as
$\mathcal{A}''$.\index{$A''$@$\mathcal{A}''$ (the double commutant of
$\mathcal{A}$)} By definition, we have that
\begin{equation}
\label{mathcal{A} subseteq mathcal{A}''}
	\mathcal{A} \subseteq \mathcal{A}''.
\end{equation}
We shall be interested in conditions which imply that $\mathcal{A}
= \mathcal{A}''$.  Of course it is necessary that $\mathcal{A}$ be
an algebra of operators on $V$ for this to hold.

	Let us note the following.

\beginlemma 
\label{triple commutant equals commutant}
If $\mathcal{A}$ is any set of operators on $V$, then $\mathcal{A}'''
= \mathcal{A}'$ (where $\mathcal{A}''' = (\mathcal{A}'')'$ denotes the
triple commutant of $\mathcal{A}$).
\end{lemma}

	As in (\ref{mathcal{A} subseteq mathcal{A}''}), $\mathcal{A}'
\subseteq \mathcal{A}'''$.  Conversely, if $T$ lies in
$\mathcal{A}'''$, then $T$ commutes with all elements of
$\mathcal{A}''$, and hence $T$ commutes with all elements of
$\mathcal{A}$, because $\mathcal{A} \subseteq \mathcal{A}''$.  Thus $T
\in \mathcal{A}'$.

	The next general fact about algebras of operators will also
be useful.

\beginproposition
\label{inverse in the algebra when it exists}
Let $V$ be a vector space over the field $k$, and let $\mathcal{A}$ be
an algebra of operators on $V$.  Suppose that $T$ lies in
$\mathcal{A}$, and that $T$ is invertible as a linear operator on $V$.
Then $T^{-1}$ also lies in $\mathcal{A}$.
\end{proposition}

	If $T$ is any linear operator on $V$ and $P(z) = \sum_{j=0}^n
a_j \, z^j$ is a polynomial with coefficients in $k$, then we write
$P(T)$ for the linear operator $\sum_{j=0}^n a_j \, T^j$.  As usual,
$z^0$ is interpreted here as being $1$, and $T^0$ is interpreted as
being the identity operator on $V$.

	Because $\mathcal{L}(V)$ has finite dimension as a vector
space over $V$, there is a nontrivial polynomial $P$ with coefficients
in $k$ such that $P(T) = 0$.  In other words, the operators $T^j$, $j
= 0, \ldots, n$, cannot be linearly independent in $\mathcal{L}(V)$
when $n$ is large enough (e.g., if $n$ is equal to the square of the
dimension of $V$), and a linear relation among them leads to a
polynomial $P$ such that $P(T) = 0$.

	If $T$ is invertible, then we may assume that the constant
term in $P$ is nonzero.  That is, if the constant term were $0$, then
we can write $P(z)$ as $z \, P_1(z)$, where $P_1(z)$ is another
polynomial (with degree $1$ less than the degree of $P$).  In this
case we have that $T \, P_1(T) = 0$, and hence $P_1(T) = 0$, since
$T$ is invertible.  By repeating this process as necessary, we can
reduce to the case where the constant term is not $0$.

	Thus we assume that $P(0) \ne 0$.  Now, $P(z) - P(0)$ is
a polynomial whose constant term vanishes, and so we may write
it as $z \, Q_1(z)$, where $Q_1(z)$ is a polynomial of degree $1$
less than that of $P$.  We obtain that $T \, Q_1(T) = P(0) \, I$,
which is the same as saying that $T^{-1} = P(0)^{-1} \, Q_1(T)$.
Thus $T^{-1}$ can be expressed as a polynomial in $T$, and the
proposition follows.

	For the original polynomial $P(z)$ one can in fact take $P(z)
= \det (T - z \, I)$, because $P(T) = 0$ in this event by the
Cayley--Hamilton theorem.  The constant term $P(0) = \det T$ is already
nonzero when $T$ is invertible.

\subsection{Nice algebras of operators}
\label{Nice algebras of operators}

	Let $V$ be a vector space over a field $k$ again.

\begindefinition
\label{invariant subspaces of an algebra of operators}
\index{invariant subspaces (of an algebra of operators)}
Let $\mathcal{A}$ be an algebra of operators on $V$.  A vector
subspace $W$ of $V$ is \emph{invariant} under $\mathcal{A}$ if $T(W)
\subseteq W$ for all $T$ in $\mathcal{A}$.
\end{definition}

\begindefinition
\label{def of nice algebras of operators}
\index{nice algebras of operators}
An algebra $\mathcal{A}$ of operators on $V$ is said to be
a \emph{nice algebra of operators} if every vector subspace
$W$ of $V$ which is invariant under $\mathcal{A}$ is also
invariant under $\mathcal{A}''$.
\end{definition}

	Note that every vector subspace $W$ of $V$ which is invariant
under $\mathcal{A}''$ is automatically invariant under $\mathcal{A}$,
since $\mathcal{A} \subseteq \mathcal{A}''$.  Thus one can rephrase
Definition \ref{def of nice algebras of operators} as saying that an
algebra of operators $\mathcal{A}$ on $V$ is a nice algebra of
operators if $\mathcal{A}$ and $\mathcal{A}''$ have the same invariant
subspaces.

\beginlemma
\label{invariant subspaces from kernels and images}
Let $S$ and $T$ be linear operators on $V$, and suppose that
$S \circ T = T \circ S$.  If $W$ denotes the image of $S$ and
$Z$ denotes the kernel of $S$, then $T(W) \subseteq W$ and
$T(Z) \subseteq Z$.
\end{lemma}

	This is easy to check.

\beginlemma
\label{invariant subspaces from elements of the commutant}
Let $\mathcal{A}$ be an algebra of linear operators on $V$.
If $S$ lies in the commutant $\mathcal{A}'$ of $\mathcal{A}$,
then the image and kernel of $S$ are invariant subspaces of
the double commutant $\mathcal{A}''$.
\end{lemma}

	This follows from Lemma \ref{invariant subspaces from kernels
and images}.

	Thus, to show that an algebra is nice, one can try to show
that all of its invariant subspaces come from linear operators in
the commutant in this way.  

\beginlemma
\label{complementary subspaces and commuting operators}
Let $\mathcal{A}$ be an algebra of operators on $V$.  Suppose that
$U_1$, $U_1$ are two vector subspaces of $V$ which are complementary,
so that $U_1 \cap U_2 = \{0\}$ and $U_1 + U_2 = V$.  Let $P : V \to V$
denote the projection of $V$ onto $U_1$ along $U_2$, i.e., $P(u_1 +
u_2) = u_1$ for all $u_1 \in U_1$ and $u_2 \in U_2$.  Then $P$ lies in
$\mathcal{A}'$ if and only if $U_1$ and $U_2$ are invariant subspaces
of $\mathcal{A}$.  In this case $U_1$ and $U_2$ are also invariant
subspaces of the double commutant $\mathcal{A}''$.
\end{lemma}

	Assume first that $U_1$ and $U_2$ are invariant subspaces of
$\mathcal{A}$.  Thus if $T$ lies in $\mathcal{A}$, $u_1$ lies in
$U_1$, and $u_2$ lies in $U_2$, then $T(u_1)$ is also an element of
$U_1$ and $T(u_2)$ is an element of $U_2$, and hence
\begin{equation}
	P(T(u_1 + u_2)) = P(T(u_1) + T(u_2)) 
		= T(u_1) = T(P(u_1 + u_2)).
\end{equation}
This shows that $T$ commutes with $P$.  We conclude that $P$ lies
in $\mathcal{A}'$, since this applies to any $T$ in $\mathcal{A}$.

	Conversely, if $P$ lies in $\mathcal{A}'$, then $U_1$ and
$U_2$ are invariant subspaces for $\mathcal{A}''$, and hence for
$\mathcal{A}$, because of Lemma \ref{invariant subspaces from elements
of the commutant}.  More precisely, this uses the fact that $U_1$ is
the image of $P$ and $U_2$ is the kernel of $P$, by construction.

\subsection{$*$-Algebras}
\label{subsection on *-algebras}

	In this subsection we assume that $V$ is a vector space over a
symmetric field $k$ and that $V$ is equipped with an inner product
$\langle \cdot, \cdot \rangle$.

\begindefinition
\label{def of *-algebras of linear operators}
\index{algebras*@$*$-algebras of linear operators}
A subset $\mathcal{A}$ of $\mathcal{L}(V)$ is said to be a
\emph{$*$-algebra of linear operators on $V$} if it is an
algebra and if the adjoint $T^*$ of an operator $T$ lies in
$\mathcal{A}$ whenever $T$ does.
\end{definition}

	The algebra $\mathcal{L}(V)$ and the set of scalar multiples
of the identity are $*$-algebras.

\beginlemma 
If $\mathcal{A}$ is a $*$-algebra of operators on $V$, then the
commutant $\mathcal{A}'$ of $\mathcal{A}$ is a $*$-algebra of
operators on $V$ as well.
\end{lemma}

	This is easy to check.

\beginlemma
\label{W invariant under mathcal{A} implies W^perp is too}
If $\mathcal{A}$ is a $*$-algebra of linear operators on $V$, and if
$W$ is a vector subspace of $V$ which is invariant under
$\mathcal{A}$, then $W^\perp$ is also invariant under $\mathcal{A}$.
\end{lemma}

	The hypothesis that $W$ be invariant under $\mathcal{A}$
can be rewritten as
\begin{equation}
	\langle T(w), u \rangle = 0
\end{equation}
for all $w \in W$, $u \in W^\perp$, and $T \in \mathcal{A}$.  This
implies in turn that
\begin{equation}
	\langle w, T^*(u) \rangle = 0
\end{equation}
for all $w \in W$, $u \in W^\perp$, and $T \in \mathcal{A}$.  In other
words, $T^*(W^\perp) \subseteq W^\perp$ for all $T \in \mathcal{A}$.
Because $\mathcal{A}$ is assumed to be a $*$-algebra, this is equivalent
to saying that $S(W) \subseteq W$ for all $S \in \mathcal{A}$, as
desired.

\begincorollary
\label{*-algebras are nice algebras}
If $\mathcal{A}$ is a $*$-algebra of operators on $V$, then
$\mathcal{A}$ is a nice algebra of operators on $V$.
\end{corollary}

	This follows easily from Lemmas \ref{W invariant under
mathcal{A} implies W^perp is too} and \ref{complementary subspaces and
commuting operators}.

\subsection{Algebras from group representations}
\label{Algebras from group representations}

	Let $G$ be a finite group, $V$ be a vector space over a field
$k$, and $\rho$ a representation of $G$ on $V$.  Consider the set
$\mathcal{A}$ of linear operators on $V$ defined by
\begin{equation}
	\mathcal{A} = \Span \{ \rho_x : x \in G \}.
\end{equation}
Here ``span'' means the ordinary linear span inside the vector space
$\mathcal{L}(V)$.  Thus
\begin{eqnarray}
	&& \hbox{the dimension of $\mathcal{A}$ as a vector space is less}  \\
	&& \hbox{than or equal to the number of elements in $G$}.
								\nonumber
\end{eqnarray}
On the other hand,
\begin{equation}
	\hbox{$\mathcal{A}$ is an algebra of linear operators on $V$.}
\end{equation}
This is easy to check, using the fact that $\rho$ is a representation
of $G$ on $V$ (so that $\rho_x \circ \rho_y = \rho_{xy}$ for all $x$
and $y$ in $G$, and hence $\rho_x \circ \rho_y$ lies in
$\mathcal{A}$ for all $x$ and $y$ in $G$).

\beginlemma
\label{W inv under mathcal{A} iff W inv under rho}
Under the conditions just described, a vector subspace $W$ of $V$ is
invariant under $\mathcal{A}$ if and only if it is invariant under the
representation $\rho$ (as defined in Subsection \ref{Reducibility}).
\end{lemma}

	This is easy to verify.

\beginlemma
\label{mathcal{A} nice if from a group rep and char(k) right}
Under the conditions above, if $k$ has characteristic $0$, or if $k$
has positive characteristic and the number of elements of $G$ is not
divisible by the characteristic of $k$, then $\mathcal{A}$ is a nice
algebra of operators.
\end{lemma}

	This uses the fact that each subspace of $V$ which is
invariant under the representation $\rho$ has an invariant complement,
as in Subsection \ref{Reducibility, continued}, and Lemma
\ref{complementary subspaces and commuting operators}.

	Suppose now that $k$ is a symmetric field, $\langle \cdot,
\cdot \rangle$ is an inner product on $V$, and that $\langle \cdot,
\cdot \rangle$ is invariant under the representation $\rho$, as in
Subsection \ref{Inner products and representations}.  This is the same as
saying that
\begin{equation}
\label{(rho_x)^{-1} = (rho_x)^*}
	(\rho_x)^{-1} = (\rho_x)^*
\end{equation}
for each $x$ in $G$, as in (\ref{R^{-1} = R^* (preserving the inner
product)}) in Subsection \ref{Inner product spaces}.  We can rewrite
(\ref{(rho_x)^{-1} = (rho_x)^*}) as
\begin{equation}
\label{(rho_x)^* = rho_{x^{-1}}}
	(\rho_x)^* = \rho_{x^{-1}}.
\end{equation}
It follows easily that $\mathcal{A}$ is a $*$-algebra in this case, since
the adjoint of each of the $\rho_x$'s lies in $\mathcal{A}$, by
(\ref{(rho_x)^* = rho_{x^{-1}}}).

	In other words,
\begin{eqnarray}
	&& \hbox{the correspondence $T \mapsto T^*$ has the effect of}	\\
	&& \hbox{``linearizing'' the correspondence $x \mapsto x^{-1}$ on $G$.}
								\nonumber
\end{eqnarray}

	Note that because the inner product is assumed to be invariant
under $\rho$, one gets immediately that the orthogonal complement of
an invariant subspace is also invariant, as in Subsection \ref{Inner
products and representations}.  This also follows from the general
result for $*$-algebras in Lemma \ref{W invariant under mathcal{A}
implies W^perp is too}.

\subsection{Regular representations}
\label{Regular representations}

	Let $G$ be a finite group and $k$ be a field.  Define $V$ to
be the vector space of $k$-valued functions on $G$, and consider the
left regular representation $L_x$, $x \in G$, and right regular
representation $R_x$, $x \in G$, of $G$ on $V$, as in Subsection
\ref{Representations of finite groups}.  Let $\mathcal{A}_L$ and
$\mathcal{A}_R$ denote the algebras of operators on $V$ generated by
the linear operators in the left and right regular representations,
respectively.  We would like to show that
\begin{equation}
\label{mathcal{A}_L' = mathcal{A}_R and mathcal{A}_R' = mathcal{A}_L}
	\mathcal{A}_L' = \mathcal{A}_R \quad\hbox{and}\quad
		\mathcal{A}_R' = \mathcal{A}_L.
\end{equation}
Once we do this, it follows that $\mathcal{A}_L'' = \mathcal{A}_L$ and
$\mathcal{A}_R'' = \mathcal{A}_R$.

	The inclusions $\mathcal{A}_R \subseteq \mathcal{A}_L'$ and
$\mathcal{A}_L \subseteq \mathcal{A}_R'$ follow from (\ref{L_x circ
R_y = R_y circ L_x}) in Subsection \ref{Representations of finite
groups}.  Conversely, suppose that $S$ is an operator on $V$ that lies
in $\mathcal{A}_L'$.  This condition is equivalent to
\begin{equation}
\label{S circ L_x = L_x circ S}
	S \circ L_x = L_x \circ S
\end{equation}
for all $x$ in $G$.  We would like to show that $S$ can be written as
a linear combination of the $R_y$'s for $y$ in $G$.  For each $w$ in
$G$, let $\phi_w(z)$ be the function on $G$ defined by $\phi_w(z) = 1$ when
$z = w$ and $\phi_w(z) = 0$ when $z \ne w$, as in Subsection
\ref{Representations of finite groups}.  Thus the $\phi_w$'s form a basis
for $V$, and $L_x(\phi_w) = \phi_{xw}$ for all $x$, $w$ in $G$ (which was
mentioned before in (\ref{L_x(phi_w) = phi_{xw}, R_x(phi_w) = phi_{w
x^{-1}}})).  Using (\ref{S circ L_x = L_x circ S}) we obtain that
\begin{equation}
\label{S(phi_x) = S(L_x(phi_e)) = L_x(S(phi_e))}
	S(\phi_x) = S(L_x(\phi_e)) = L_x(S(\phi_e))
\end{equation}
for all $x$ in $G$.  Of course $S(\phi_e)$ is a function on $G$, which
we can write as
\begin{equation}
	S(\phi_e) = \sum_{w \in G} c_w \, \phi_w,
\end{equation}
where in fact $c_w = S(\phi_e)(w)$ for all $w$ in $G$.  This and
(\ref{S(phi_x) = S(L_x(phi_e)) = L_x(S(phi_e))}) permit us to write $S(\phi_x)$
as
\begin{equation}
	S(\phi_x) = \sum_{w \in G} c_w \, \phi_{xw} 
		= \sum_{w \in G} c_w \, R_{w^{-1}}(\phi_x),
\end{equation}
since $R_u(\phi_x) = \phi_{xu^{-1}}$, as in (\ref{L_x(phi_w) = phi_{xw},
R_x(phi_w) = phi_{w x^{-1}}}).  Thus
\begin{equation}
	S = \sum_{w \in G} c_w \, R_{w^{-1}},
\end{equation}
which is what we wanted.  In the same way one can show that if an
operator commutes with $R_y$ for all $y$ in $G$, then the operator is
a linear combination of $L_u$'s.  Therefore $\mathcal{A}_L' \subseteq
\mathcal{A}_R$ and $\mathcal{A}_R' \subseteq \mathcal{A}_L$.  This
completes the proof of (\ref{mathcal{A}_L' = mathcal{A}_R and
mathcal{A}_R' = mathcal{A}_L}).

\subsection{Expanding an algebra}
\label{Expanding an algebra}

	Let $V$ be a vector space over a field $k$, and let
$\mathcal{A}$ be an algebra of linear operators on $V$.  Fix a
positive integer $n$, and let $V^n$ denote the vector space of
functions $f$ from $\{1, 2, \ldots, n\}$ to $V$.  In effect, $V^n$ is
a direct sum of $n$ copies of $V$, and this is a convenient way to
express it.

	If $T$ is a linear operator on $V$, then we can associate to
it a linear operator $\widehat{T}$ on $V^n$ by
\begin{equation}
	\widehat{T}(f)(j) = T(f(j))
\end{equation}
for each $f \in V^n$ and $j = 1, 2, \ldots, n$.  In other words,
$\widehat{T}$ is the same as $T$ in the $V$ directions in $V^n$, and
does nothing in the other directions.  If one thinks of $V^n$ as a
direct sum of $n$ copies of $V$, then $\widehat{T}$ acts in the same
way as $T$ on each of these copies, with no mixing between the
copies.

	It is not hard to verify that the correspondence $T \mapsto
\widehat{T}$ from $\mathcal{L}(V)$ to $\mathcal{L}(V^n)$ is linear and
cooperates with compositions, so that $\widehat{S \circ T} =
\widehat{S} \circ \widehat{T}$ for all $S, T \in \mathcal{L}(V)$.
Also, if $T$ is the identity operator on $V$, then $\widehat{T}$ is
the identity operator on $V^n$.  Set
\begin{equation}
	\mathcal{A}_n = \{\widehat{T} : T \in \mathcal{A}\}.
\end{equation}
From the preceding remarks it follows that $\mathcal{A}_n$ is an
algebra of operators on $V^n$, since $\mathcal{A}$ is an algebra of
operators on $V$.

\beginproposition
\label{prop about (mathcal{A}'')_n = (mathcal{A}_n)''}
$(\mathcal{A}'')_n = (\mathcal{A}_n)''$.
\end{proposition}

	Here $(\mathcal{A}'')_n = \{\widehat{T} : T \in
\mathcal{A}''\}$, in analogy with $\mathcal{A}_n$, while
$(\mathcal{A}_n)''$ denotes the double commutant of $\mathcal{A}_n$ as
an algebra of operators in $\mathcal{L}(V^n)$.  Thus the proposition
states that one gets the same result whether one first takes the
double commutant of $\mathcal{A}$, and then uses the correspondence $T
\mapsto \widehat{T}$ to pass to $\mathcal{L}(V^n)$, or one first uses
this correspondence to pass to $\mathcal{L}(V^n)$, and then takes the
double commutant there.

	To prove the proposition, let us look first at what
$(\mathcal{A}_n)'$ is.  For $\ell = 1, 2, \ldots, n$, let $P_\ell$
denote the natural projection of $V^n$ onto the $\ell$th copy of $V$
inside $V^n$, defined by
\begin{eqnarray}
	P_\ell(f)(j) & = & f(j) \quad\hbox{when } j = \ell,		\\
		     & = & 0    \ \qquad\hbox{when } j \ne \ell.    \nonumber
\end{eqnarray}
Thus $P_h \circ P_\ell = 0$ when $h \ne \ell$, and
\begin{equation}
\label{sum_{ell = 1}^n P_ell = Identity transformation on V^n}
	\sum_{\ell = 1}^n P_\ell = \hbox{ identity transformation on $V^n$.}
\end{equation}
For any $T \in \mathcal{L}(V)$, 
\begin{equation}
\label{widehat{T} circ P_ell = P_ell circ widehat{T}}
	\widehat{T} \circ P_\ell = P_\ell \circ \widehat{T},
\end{equation}
by the way that $\widehat{T}$ is defined.  Now let $M$ be any linear
transformation on $V^n$, and set $M_{h,\ell} = P_h \circ M \circ P_\ell$.
To say that $M$ lies in $(\mathcal{A}_n)'$ means that 
\begin{equation}
	M \circ \widehat{T} = \widehat{T} \circ M
\end{equation}
for all $T \in \mathcal{A}$.  Hence $M \in (\mathcal{A}_n)'$ if and only
if
\begin{equation}
\label{M_{h, ell} circ widehat{T} = widehat{T} circ M_{h, ell}}
	M_{h, \ell} \circ \widehat{T} = \widehat{T} \circ M_{h, \ell}
\end{equation}
for all $T \in \mathcal{A}$ and $h, \ell = 1, 2, \ldots, n$, because
of the preceding observations.

	For $\ell = 1, 2, \ldots, n$, let $\theta_\ell : V \to V^n$
denote the linear mapping defined by
\begin{eqnarray}
	\theta_\ell(v)(j) & = & v \quad\hbox{when } j = \ell		\\
			  & = & 0 \quad\hbox{when } j \ne \ell.	  \nonumber
\end{eqnarray}
In other words, $\theta_\ell$ is the obvious identitification between
$V$ and the $\ell$th copy of $V$ inside $V^n$.  For each linear
transformation $M$ on $V^n$ and each $h, \ell = 1, 2, \ldots, n$, let
$\widetilde{M}_{h, \ell}$ be the linear transformation on $V$ such that
\begin{equation}
	M_{h, \ell} \circ \theta_\ell = \theta_h \circ \widetilde{M}_{h, \ell}.
\end{equation}
By construction, $M_{h, \ell}$ essentially corresponds to a linear
mapping from the $\ell$th copy of $V$ in $V^n$ to the $h$th copy of
$V$ in $V^n$, and with $\widetilde{M}_{h, \ell}$ we rewrite this as
a mapping on $V$ in the obvious way.

	We also have that
\begin{equation}
	\widehat{T} \circ \theta_\ell = \theta_\ell \circ T
\end{equation}
for all $T \in \mathcal{L}(V)$ and $\ell = 1, 2, \ldots, n$, again
just by the way that everything is defined here.  The bottom line
is that (\ref{M_{h, ell} circ widehat{T} = widehat{T} circ M_{h, ell}})
is equivalent to
\begin{equation}
	\widetilde{M}_{h, \ell} \circ T = T \circ \widetilde{M}_{h, \ell}
\end{equation}
for all $T \in \mathcal{A}$ and $h, \ell = 1, 2, \ldots, n$.  To
summarize, $M \in \mathcal{L}(V^n)$ lies in $(\mathcal{A}_n)'$ if and
only if the ``blocks'' $\widetilde{M}_{h, \ell}$, $1 \le h, \ell \le
n$, in $\mathcal{L}(V)$ all lie in $\mathcal{A}'$.

	Now suppose that $S \in \mathcal{L}(V)$ lies in $\mathcal{A}''$.
Thus $S$ commutes with all linear operators on $V$ in $\mathcal{A}'$.
Using this, it is not hard to show that
\begin{equation}
\label{M circ widehat{S} = widehat{S} circ M}
	M \circ \widehat{S} = \widehat{S} \circ M
\end{equation}
for all $M \in \mathcal{L}(V^n)$ which lie in $(\mathcal{A}_n)'$.
To be more precise, (\ref{M circ widehat{S} = widehat{S} circ M})
holds if and only if
\begin{equation}
	M_{h, \ell} \circ \widehat{S} = \widehat{S} \circ M_{h, \ell}
\end{equation}
for all $h, \ell = 1, 2, \ldots, n$, as before.  This is in turn equivalent
to
\begin{equation}
	\widetilde{M}_{h, \ell} \circ S = S \circ \widetilde{M}_{h, \ell}
\end{equation}
for $h, \ell = 1, 2, \ldots, n$, as operators now on $V$ instead of
$V^n$.  This condition holds for $S \in \mathcal{A}''$ because
$\widetilde{M}_{h, \ell}$ lies in $\mathcal{A}'$ for $h, \ell = 1, 2,
\ldots, n$.

	From (\ref{M circ widehat{S} = widehat{S} circ M}) we obtain
the inclusion
\begin{equation}
	(\mathcal{A}'')_n \subseteq (\mathcal{A}_n)''.
\end{equation}
Next we wish to show that
\begin{equation}
	(\mathcal{A}_n)'' \subseteq (\mathcal{A}'')_n.
\end{equation}

	Let $M \in \mathcal{L}(V^n)$ be any element of $(\mathcal{A}_n)''$.
The first point is that each of the $n$ copies of $V$ inside $V^n$ are
invariant under $M$, which is the same as saying that
\begin{equation}
	M_{h, \ell} = 0 \quad\hbox{when } h \ne \ell.
\end{equation}
This follows from the fact that each $P_\ell$ lies in
$(\mathcal{A}_n)'$, so that $M$ commutes with each $P_\ell$.

	If $\pi$ is any permutation on $\{1, 2, \ldots, n\}$, consider
the linear transformation on $V^n$ which takes $f(j)$ in $V^n$ to
$f(\pi(j))$.  This transformation commutes with $\widehat{T}$ on $V^n$
for any $T \in \mathcal{L}(V)$, and thus this transformation lies in
$(\mathcal{A}_n)'$.

	As a result, $M$ in $(\mathcal{A}_n)''$ commutes with all of
these linear transformations on $\mathcal{L}(V^n)$ obtained from
permutations on $\{1, 2, \ldots, n\}$.  Using this it is not hard to
verify that $M = \widehat{R}$ for some linear operator $R$ on $V$.

	The remaining observation is that $R$ lies in $\mathcal{A}''$.
Indeed, if $S$ is any linear operator on $V$ which lies in
$\mathcal{A}'$, then it is easy to see that $\widehat{S}$ lies in
$(\mathcal{A}_n)'$.  This implies that $M$ commutes with
$\widehat{S}$, and hence $R$ commutes with $S$.  Therefore $R \in
\mathcal{A}''$, since this applies to any $S$ in $\mathcal{A}'$.  In
other words, $M = \widehat{R}$ lies in $(\mathcal{A}'')_n$, as
desired.  This completes the proof of Proposition \ref{prop about
(mathcal{A}'')_n = (mathcal{A}_n)''}.

	Suppose now that $k$ is a symmetric field and that $V$ is
equipped with an inner product $\langle \cdot, \cdot \rangle$.  We
define an inner product $\langle \cdot, \cdot \rangle_n$ on $V^n$ by
\begin{equation}
	\langle f, h \rangle_n = \sum_{j=1}^n \langle f(j), h(j) \rangle
\end{equation}
for all $f, h \in V^n$.  In this way $V^n$ is an orthogonal direct
sum of $n$ copies of $V$, all with the same inner product as on $V$.

\beginlemma
For each linear operator $T$ on $V$, $\widehat{(T^*)} =
(\widehat{T})^*$ (where $T^*$ is defined in terms of the
inner product on $V$, while $(\widehat{T})^*$ is defined using
the inner product on $V^n$).
\end{lemma}

	This is easy to check.

\begincorollary
\label{mathcal{A}_n is a *-algebra, and hence nice}
Under the conditions above, if $\mathcal{A}$ is a $*$-algebra of
operators on $V$, then $\mathcal{A}_n$ is a $*$-algebra of operators
on $V^n$, and hence is nice.
\end{corollary}

	Now assume that $k$ is a field, $V$ is a vector space over
$k$, and $\mathcal{A}$ is the algebra of operators on $V$ generated by
a representation $\rho$ of a finite group $G$ on $V$.  We can
associate to $\rho$ a representation $\rho^n$ of $G$ on $V^n$ in a
simple way, namely by setting
\begin{equation}
	\rho^n_x = \widehat{\rho_x}
\end{equation}
for each $x$ in $G$.

\beginlemma
\label{mathcal{A}_n in the case of a group representation}
Under the conditions just described, $\mathcal{A}_n$ is the
algebra of operators generated by the representation $\rho^n$
of $G$ on $V^n$.  
\end{lemma}

	This follows easily from the definitions.

\begincorollary
\label{mathcal{A}_n in the case of a group representation, char(k) right}
If we also assume that either $k$ has characteristic $0$, or $k$ has
positive characteristic and the number of elements of the group $G$ is
not divisible by the characteristic of $k$, then $\mathcal{A}_n$ is
nice.
\end{corollary}

	This uses Lemma \ref{mathcal{A} nice if from a group rep and
char(k) right} in Subsection \ref{Algebras from group representations}.

\subsection{Very nice algebras of operators}
\label{Very nice algebras of operators}

\begindefinition
\label{def of very nice algebras of operators}
\index{very nice algebras of operators}
Let $k$ be a field and $V$ a vector space over $k$.  An algebra of
operators $\mathcal{A}$ on $V$ is said to be \emph{very nice} if
$\mathcal{A}_n \subseteq \mathcal{L}(V^n)$ (defined in Subsection
\ref{Expanding an algebra}) is a nice algebra of operators for
each positive integer $n$.
\end{definition}

	The $n = 1$ case corresponds to nice algebras of operators.
From Corollaries \ref{mathcal{A}_n is a *-algebra, and hence nice} and
\ref{mathcal{A}_n in the case of a group representation, char(k)
right} we know that $\mathcal{A}$ is very nice if either it is a
$*$-algebra, assuming that $k$ is a symmetric field and that $V$ is
equipped with an inner product, or if $\mathcal{A}$ is generated by a
representation of a finite group on $V$ and either the characteristic
of $k$ is $0$ or the number of elements of $G$ is not divisible by the
characteristic of $k$ when it is positive.

\begintheorem
\label{mathcal{A}'' = mathcal{A} if mathcal{A} is very nice}
Let $k$ be a field, $V$ a vector space over $k$, and
$\mathcal{A}$ an algebra of operators on $V$.  If $\mathcal{A}$
is very nice, then $\mathcal{A}'' = \mathcal{A}$.
\end{theorem}

	This is a relative of the ``double commutant theorem'', and
for the proof (including some of the preliminary steps in the previous
subsections) we are essentially following the argument indicated on p118
of \cite{Douglas}.  As in \cite{Douglas}, this argument applies to
``Von Neumann algebras'' (with suitable adjustments).

	Let $k$, $V$, and $\mathcal{A}$ be as in the statement of the
theorem.  As in (\ref{mathcal{A} subseteq mathcal{A}''}) in Subsection
\ref{Basic notions (Algebras of linear operators)}, $\mathcal{A}
\subseteq \mathcal{A}''$ holds automatically, and so we need only show
that $\mathcal{A}'' \subseteq \mathcal{A}$.

	 Let us first check that for each $v \in V$ and $S \in
\mathcal{A}''$ there is a $T \in \mathcal{A}$ such that $S(v) = T(v)$,
assuming only that $\mathcal{A}$ is nice.  Fix $v \in V$, and set
\begin{equation}
	W = \{T(v) : T \in \mathcal{A}\}.
\end{equation}
It is easy to see that $W$ is a vector subspace of $V$ which is
invariant under $\mathcal{A}$, because $\mathcal{A}$ is an algebra of
operators.  If $\mathcal{A}$ is nice, then $W$ is also invariant under
$\mathcal{A}''$.  Of course $v \in W$, since $I \in \mathcal{A}$, and
therefore $S(v) \in W$ for any $S \in \mathcal{A}''$, because $W$ is
invariant under $\mathcal{A}''$.  This shows that for each $S \in
\mathcal{A}''$ there is a $T \in \mathcal{A}$ such that $T(v) = S(v)$,
by the definition of $W$.

	Now let $n$ be a positive integer, and let $V^n$,
$\mathcal{A}_n$, etc., be as in Subsection \ref{Expanding an algebra}.
The preceding observation can be applied to $V^n$ and $\mathcal{A}_n$
(assuming that $\mathcal{A}$ is very nice, so that $\mathcal{A}_n$ is
nice) to obtain that for each $f \in V^n$ and each $M \in
(\mathcal{A}_n)''$ there is an $N \in \mathcal{A}_n$ such that $N(f) =
M(f)$.

	From Proposition \ref{prop about (mathcal{A}'')_n =
(mathcal{A}_n)''} we know that $(\mathcal{A}_n)'' =
(\mathcal{A}'')_n$.  Thus the previous assertion can be rephrased as
saying that for each $f \in V^n$ and each $M \in (\mathcal{A}'')_n$
there is an $N \in \mathcal{A}_n$ such that $N(f) = M(f)$.  This can
be rephrased again as saying that for each $f \in V^n$ and each
$S \in \mathcal{A}''$ there is a $T \in \mathcal{A}$ such that
\begin{equation}
	\widehat{T}(f) = \widehat{S}(f).
\end{equation}

	This is in turn equivalent to the following statement.  If
$v_1, v_2, \ldots, v_n$ are elements of $V$, and if $S \in
\mathcal{A}''$, then there is a $T \in \mathcal{A}$ such that $T(v_j)
= S(v_j)$ for $j = 1, 2, \ldots, n$.

	If $n$ is the dimension of $V$, then we can choose $v_1, v_2,
\ldots, v_n$ to be a basis of $V$.  In this case the preceding
statement reduces to saying that if $S \in \mathcal{A}''$, then there
is a $T \in \mathcal{A}$ such that $T = S$.  This implies that
$\mathcal{A}'' \subseteq \mathcal{A}$, and hence $\mathcal{A}'' =
\mathcal{A}$, since the other inclusion is automatic.  This completes
the proof of the theorem.

\subsection{Some other approaches}
\label{Some other approaches}

	Let $k$ be a field, $V$ a vector space over $k$, and
$\mathcal{A}$ an algebra of operators on $V$.  Since $\mathcal{A}
\subseteq \mathcal{A}''$ automatically, another way to try to show
that $\mathcal{A} = \mathcal{A}''$ is to show that the dimension of
$\mathcal{A}$ is equal to the dimension of $\mathcal{A}''$, as vector
spaces over $k$.  Of course this is easier to do if $\mathcal{A}''$ is
known and reasonably simple.

	If $\mathcal{A}'$ consists of only scalar multiples of the
identity operator on $V$, then $\mathcal{A}'' = \mathcal{L}(V)$, and
the dimension of $\mathcal{A}''$ is equal to the square of the
dimension of $V$.  In this case we can just say that $\mathcal{A} =
\mathcal{L}(V)$ if $\mathcal{A}$ has dimension equal to the square of
the dimension of $V$.  Let us give a couple of other criteria for
$\mathcal{A}$ to be equal to $\mathcal{L}(V)$.

\beginlemma
\label{criterion for mathcal{A} = mathcal{L}(V), 1}
Let $\mathcal{A}$ be an algebra of linear operators on $V$ which
satisfies the following two conditions: (a) for any vectors $v$, $w$
in $V$ with $v \ne 0$ there is a linear operator $T$ in $\mathcal{A}$
such that $T(v) = w$; (b) there is a nonzero linear operator on $V$
with rank $1$ which lies in $\mathcal{A}$.  Then $\mathcal{A} =
\mathcal{L}(V)$.
\end{lemma}

	Condition (a) in the hypothesis is equivalent to asking that
$\mathcal{A}$ be irreducible in the sense of Definition
\ref{definition of an irreducible algebra of operators}.

	To show that $\mathcal{A} = \mathcal{L}(V)$, it is enough to
show that every linear operator on $V$ with rank $1$ lies in
$\mathcal{A}$.  By assumption, there is a nonzero operator $R$ of rank
$1$ which lies in $\mathcal{A}$, and we can write $R$ as $R(u) = f(u)
\, z$, where $z$ is a nonzero element of $V$, and $f$ is a nonzero
linear mapping from $V$ to $k$, i.e., a nonzero linear functional on
$V$.  If $w$ is any other nonzero vector in $v$, then there is an
operator $T$ in $\mathcal{A}$ such that $T(z) = w$, and hence $(T \circ
R)(u) = f(u) \, T(z) = f(u) \, w$ and $T \circ R$ lies in
$\mathcal{A}$.

	Let $h$ be any nonzero linear functional on $V$.  We would
like to show that there is a linear operator $S$ in $\mathcal{A}$ such
that $h = f \circ S$.  If we can do that, then $(T \circ R \circ S)(u)
= h(u) \, w$ and $T \circ R \circ S$ lies in $\mathcal{A}$, and we get
that any linear operator on $V$ of rank $1$ lies in $\mathcal{A}$.

	Consider the set 
\begin{equation}
	Z = \{y \in V : f \circ S(y) = 0 
			\hbox{ for all } S \in \mathcal{A}\}.
\end{equation}
The irreducibility assumption for $\mathcal{A}$ implies that $Z$
contains only the zero vector in $V$.  On the other hand, the set of
linear functionals on $V$ of the form $f \circ S$ is a vector subspace
of the vector space of all linear functionals on $V$.  The set $Z$ is
exactly the intersection of the kernels of these linear functionals,
and if they did not span the whole space of linear functionals on $V$,
then $Z$ would contain a nonzero vector.  This is a basic result of
linear algebra.  Since $Z$ does contain only the zero vector, we
obtain that the set of linear functionals on $V$ of the form $f \circ
S$, $S \in \mathcal{A}$, is equal to the set of all linear functionals
on $V$, as desired.  This completes the proof of Lemma \ref{criterion
for mathcal{A} = mathcal{L}(V), 1}.

\beginlemma
\label{criterion for mathcal{A} = mathcal{L}(V), 2}
Let $\mathcal{A}$ be an algebra of linear operators on $V$ such that
for any pair $v_1$, $v_2$ of linearly independent vectors in $V$ and
any arbitrary pair of vectors $w_1$, $w_2$ in $V$ there is a $T$
in $\mathcal{A}$ such that $T(v_1) = w_1$, $T(v_2) = w_2$.  Then
$\mathcal{A} = \mathcal{L}(V)$.
\end{lemma}

	We may as well assume that $V$ has dimension at least $2$,
since the conclusion of the lemma is automatic when the dimension of
$V$ is $1$.  The hypothesis of the lemma is also vacuous in this
case.  Because of Lemma \ref{criterion for mathcal{A} =
mathcal{L}(V), 1}, it suffices to show that $\mathcal{A}$ contains a
nonzero operator of rank $1$.  

	Let $R$ be any nonzero operator in $\mathcal{A}$.  If $R$ has
rank $1$, then we are finished, and so we assume that $R$ has rank at
least $2$.  This means that there are vectors $u_1$, $u_2$ in $V$ such
that $R(u_1)$, $R(u_2)$ are linearly independent.  By assumption,
there is an operator $T$ in $\mathcal{A}$ such that $T(R(u_1)) = 0$
and $T(R(u_2)) \ne 0$.  Thus $T \circ R$ is nonzero and has rank
smaller than the rank of $R$, because the dimension of the image has
been reduced, or, equivalently, the dimension of the kernel has been
increased.  We also have that $T \circ R$ is in $\mathcal{A}$, since
$R$ and $T$ are elements of $\mathcal{A}$.  By repeating this process,
as needed, we can get a nonzero element of $\mathcal{A}$ of rank $1$.

\beginlemma
\label{criterion for mathcal{A} = mathcal{L}(V), 3}
Suppose that $\mathcal{A}$ is an algebra of linear operators on $V$
such that $\mathcal{A}'$ consists of only the scalar multiples of the
identity and $\mathcal{A}_n \subseteq \mathcal{L}(V^n)$ (defined in
Subsection \ref{Expanding an algebra}) is a nice algebra of operators
when $n = 2$.  Then $\mathcal{A} = \mathcal{L}(V) = \mathcal{A}''$.
\end{lemma}

	The assumptions in this lemma imply the hypothesis of Lemma
\ref{criterion for mathcal{A} = mathcal{L}(V), 2}.

\section{Vector spaces with definite scalar product}
\label{Vector spaces with definite scalar product}
\setcounter{equation}{0}

	Let $k$ be a field and let $V$ be a vector space.  Suppose
that $(v,w)$ is a \emph{definite scalar product}\index{definite scalar
product}\index{scalar product} on $V$, by which we mean the following:
(a) $(v,w)$ is a bilinear form\index{bilinear form} on $V$, i.e., it
is a function from $V \times V$ into $k$ such that $v \mapsto (v,w)$
is linear for each $w$ in $V$, and $w \mapsto (v,w)$ is linear for
each $v$ in $V$; (b) $(v,w)$ is symmetric\index{symmetric bilinear
form} in $v$ and $w$, so that $(v,w) = (w,v)$ for all $v$ and $w$ in
$V$; and (c) $(v,v) = 0$ if and only if $v = 0$.  This last can be a
somewhat complicated condition, which does not necessarily entail
positivity.  See \cite{Serre-arith}, for instance.  For finite fields
any homogeneous polynomial of degree 2 in at least $3$ variables has a
nontrivial solution, as on p338 of \cite{Cassels} and Corollary 2 on p6
of \cite{Serre-arith}.

	Two vectors $v$, $w$ in $V$ are said to be
\emph{orthogonal}\index{orthogonal vectors} if $(v,w) = 0$.
A collection of vectors $v_1, \ldots, v_k$ is said to be
\emph{orthogonal} if $(v_j, v_l) = 0$ for all $j$, $l$ such
that $1 \le j, l \le k$ and $j \ne l$.  As usual, nonzero
orthogonal vectors are linearly independent.

	If $W$ is a linear subspace of $V$, then we define the
\emph{orthogonal complement}\index{orthogonal complement}
$U^\perp$\index{$U^\perp$} of $U$ by
\begin{equation}
	U^\perp = \{v \in V : (v,u) = 0 \hbox{ for all } u \in U\}.
\end{equation}
Thus $U^\perp$ is also a linear subspace of $V$, and
\begin{equation}
	U \cap U^\perp = \{0\},
\end{equation}
since our scalar product is assumed to be definite.

	Suppose that $u_1, \ldots, u_n$ is a collection of
nonzero orthogonal vectors in $V$, and let $U$ denote the
span of $u_1, \ldots, u_n$.  Define a linear operator $P$
on $V$ by
\begin{equation}
	P(v) = \sum_{j=1}^n \frac{(v,u_j)}{(u_j,u_j)} \, u_j.
\end{equation}
We have that $P(u) = u$ when $u \in U$, $P(v) = 0$ when $v \in
U^\perp$, and 
\begin{equation}
\label{splitting property for P}
	P(w) \in U, \ w - P(w) \in U^\perp \quad\hbox{for all } w \in V.
\end{equation}
For any $w$ in $V$, there is at most one vector $y$ in $V$ such that
$y \in U$ and $w - y \in U^\perp$, because if $y'$ is another such
vector, then $y - y' \in U \cap U^\perp$.  Hence $P$ is characterized
by (\ref{splitting property for P}), and depends only on $U$, and not
the choice of $u_1, \ldots, u_n$.  We say that $P$ is the
\emph{orthogonal projection}\index{orthogonal projections} of $V$ onto
$U$.

\beginlemma
\label{orthogonalizing vectors}
If $x_1, \ldots, x_n$ are linearly independent vectors in $V$, then
there are nonzero orthogonal vectors $u_1, \ldots, u_n$ in $V$
such that $\Span (x_1, \ldots, x_n) = \Span (u_1, \ldots, u_n)$.
\end{lemma}

	As usual, this can be proved using induction on $n$.  The
orthogonalization of the last vector is obtained with the help of an
orthogonal projection onto the span of the previous vectors, where the
orthogonal projection just mentioned is derived from the induction
hypothesis.

\begincorollary
\label{orthogonal projections onto any subspace}
If $U$ is any linear subspace of $V$, then there is an orthogonal
projection of $V$ onto $U$.  In particular, $V$ is spanned by $U$
and $U^\perp$, and $(U^\perp)^\perp = U$.
\end{corollary}

	Now let $T$ be a linear operator on $V$.  The
\emph{transpose}\index{transpose of a linear operator} of $T$ is the
linear operator $T^t$\index{$T^t$ (transpose of an operator $T$)} on
$V$ such that
\begin{equation}
\label{(T(v),w) = (v, T^t(w))}
	(T(v),w) = (v, T^t(w))
\end{equation}
for all $v$, $w$ in $V$.  It is not difficult to establish the
existence of such an operator $T^t$, and uniqueness is an easy
consequence of (\ref{(T(v),w) = (v, T^t(w))}).  

	Note that
\begin{equation}
	(\alpha \, S + \beta \, T) = \alpha \, S^t + \beta \, T^t
\end{equation}
and
\begin{equation}
	(S \circ T)^t = T^t \circ S^t
\end{equation}
for any scalars $\alpha$, $\beta$ and linear operators $S$, $T$ on $V$.
The transpose of the identity operator $I$ on $V$ is itself, and
\begin{equation}
	(T^t)^t = T
\end{equation}
for any linear operator $T$ on $V$.

\begindefinition
\label{def of t-algebras}
Let $\mathcal{A}$ be an algebra of linear operators on $V$, as in
Definition \ref{def of algebras of linear operators}.  We say that
$\mathcal{A}$ is a $t$-algebra\index{algebrast@$t$-algebras of linear
operators} if $T^t$ lies in $\mathcal{A}$ whenever $T$ does.
\end{definition}

	Of course the algebra $\mathcal{L}(V)$ of all operators on $V$
is a $t$-algebra, as is the algebra consisting only of scalar multiples
of the identity.

\beginlemma
\label{commutant of a t-algebra is a t-algebra}
If $\mathcal{A}$ is a $t$-algebra of linear operators on $V$,
then the commutant $\mathcal{A}'$ is also a $t$-algebra of
linear operators on $V$.
\end{lemma}

	This is a simple exercise.

\beginlemma
\label{W^perp invariant under mathcal{A} if W is and mathcal{A} a t-algebra}
Let $\mathcal{A}$ be an algebra of operators on $V$, and let $U$ be
a linear subspace of $V$ which is invariant under $\mathcal{A}$.
If $\mathcal{A}$ is a $t$-algebra, then $U^\perp$ is invariant under
$\mathcal{A}$, and the orthogonal projection $P_U$ of $V$ onto $U$ lies
in $\mathcal{A}'$.
\end{lemma}

	The statement that $U$ is invariant under $\mathcal{A}$
is equivalent to
\begin{equation}
	(T(u), v) = 0 
\end{equation}
for all $u \in U$, $v \in U^\perp$, and $T \in \mathcal{A}$.  This
implies that
\begin{equation}
	(u, T^t(v)) = 0
\end{equation}
for all $u \in U$, $v \in U^\perp$, and $T \in \mathcal{A}$, so that
$U^\perp$ is invariant under $T^t$ for all $T$ in $\mathcal{A}$.  If
$\mathcal{A}$ is a $t$-algebra, then it follows that $U^\perp$ is
invariant under $\mathcal{A}$.  The information that $U$ and $U^\perp$
are both invariant under $\mathcal{A}$ implies that $P_U \in
\mathcal{A}'$, as in Lemma \ref{complementary subspaces and commuting
operators}.

	As a result, if $\mathcal{A}$ is a $t$-algebra of operators
on $V$, then $\mathcal{A}$ is a nice algebra of operators, in the
sense of Definition \ref{def of nice algebras of operators}.

	Let $n$ be a positive integer, and define $V^n$,
$\mathcal{A}_n$ as in Subsection \ref{Expanding an algebra}.  Fix nonzero
elements $\lambda_1, \ldots, \lambda_n$ in $k$, and consider 
\begin{equation}
	(f, h)_n = \sum_{j=1}^n \lambda_j \, (f(j), h(j)), \quad f, h \in V^n.
\end{equation}
In general, this will not be definite, but if it is, then it is easy
to see that $\mathcal{A}_n$ is a nice algebra of operators on $V^n$,
because it is a $t$-algebra with respect to $(f,h)_n$, as in Corollary
\ref{mathcal{A}_n is a *-algebra, and hence nice}.

\section{Irreducibility}
\label{section on irreducibility}
\setcounter{equation}{0}

\begindefinition
\label{definition of an irreducible algebra of operators}
\index{irreducible algebras of operators} 
Let $V$ be a vector space over a field $k$.  An algebra of operators
$\mathcal{A}$ on $V$ is said to be \emph{irreducible} if there are no
vector subspaces $W$ of $V$ which are invariant under $\mathcal{A}$
(Definition \ref{invariant subspaces of an algebra of operators})
except for $W = \{0\}$ and $W = V$.
\end{definition}

\beginlemma
\label{reformulation of irreducibility}
An algebra of operators $\mathcal{A}$ on a vector space $V$ over a
field $k$ is irreducible if and only if for every pair of vectors $v$,
$w$ in $V$ with $v \ne 0$ there is an operator $T$ in $\mathcal{A}$
such that $T(v) = w$.
\end{lemma}

	The ``if'' part is a simple consequence of the definition.
For the ``only if'' part, let $v$ be a nonzero vector in $V$, and
consider the set of vectors in $V$ of the form $T(v)$, $T \in
\mathcal{A}$.  This set of vectors is a vector subspace of $V$ which
is invariant under $\mathcal{A}$, because $\mathcal{A}$ is an algebra
of operators, and it is nonzero because the identity operator lies in
$\mathcal{A}$, so that $v$ lies in the subspace.  Thus irreducibility
of $\mathcal{A}$ implies that this subspace is all of $V$, which means
exactly that for every vector $w$ in $V$ there is an operator $T$ in
$\mathcal{A}$ such that $T(v) = w$.

\beginlemma
\label{irreducibility for group rep. and associated algebra}
If $V$ is a vector space, $G$ is a finite group, $\rho$ is a
representation of $G$ on $V$, and $\mathcal{A}$ is the algebra of
operators on $V$ generated by $\rho_x$ for $x$ in $G$ as in Subsection
\ref{Algebras from group representations}, then $\rho$ is an
irreducible representation of $G$ if and only if $\mathcal{A}$ is
irreducible.
\end{lemma}

	This follows from Lemma \ref{W inv under mathcal{A} iff
W inv under rho} in Subsection \ref{Algebras from group representations}.

\beginproposition
\label{irreducibility of mathcal{A} implies ....}
Let $V$ be a vector space over a field $k$, and let $\mathcal{A}$ be
an algebra of operators on $V$.  If $\mathcal{A}$ is irreducible, then
for each $T$ in the commutant $\mathcal{A}'$, either $T = 0$ or $T$ is
invertible as an operator on $V$.  In the latter event $T^{-1}$ also
lies in $\mathcal{A}'$.
\end{proposition}

	If $T$ lies in $\mathcal{A}'$, then the kernel and image of
$T$ are invariant subspaces for $\mathcal{A}$, as in Lemma
\ref{invariant subspaces from elements of the commutant}.  From this
it is easy to see that either $T$ is the zero operator or $T$ is
invertible when $\mathcal{A}$ is irreducible.  Proposition
\ref{inverse in the algebra when it exists} implies that $T^{-1}$
lies in $\mathcal{A}'$ when $T$ is invertible, since $\mathcal{A}'$
is an algebra of operators.

\beginproposition
\label{k alg closed, irreducibility, imply mathcal{A}' = scalars I}
Let $V$ be a vector space over an algebraically closed field $k$.  If
$\mathcal{A}$ is an irreducible algebra of operators on $V$, then the
commutant $\mathcal{A}'$ is equal to the set of multiples of the
identity operator by scalars (elements of $k$).
\end{proposition}

	The commutant of any algebra always contains the scalar
multiples of the identity, since the identity operator commutes with
all other linear mappings, and so it suffices to show that every
linear transformation in $\mathcal{A}'$ is a multiple of the identity.
This is in essence the same as one-half of the classical \emph{Schur's
lemma}, which can be proved as follows.

	Let $T$ be any linear operator in $\mathcal{A}'$.  Because $k$
is algebraically closed, there is at least one nontrivial eigenvalue
$\lambda \in k$ of $T$.  In other words, there is at least one
$\lambda$ in $k$ for which the eigenspace
\begin{equation}
	E(T, \lambda) = \{v \in V : T(v) = \lambda \, v\}
\end{equation}
contains nonzero vectors.  This is because $\lambda$ is an eigenvalue
of $T$ if and only if $\det (T - \lambda \, I) = 0$, and $\det (T -
\lambda \, I)$ is a polynomial in $\lambda$ (of degree equal to the
dimension of $V$, which is positive).  

	If $E(T, \lambda) = V$, then $T = \lambda \, I$, which is what
we want.  Thus we suppose that $E(T, \lambda)$ is not all of $V$.  It
is not hard to check that if $S$ is a linear operator on $V$ which
commutes with $T$, then
\begin{equation}
\label{S(E(T, lambda)) subseteq E(T, lambda)}
	S(E(T, \lambda)) \subseteq E(T, \lambda).
\end{equation}
In particular, this holds for all $S \in \mathcal{A}$, since $T \in
\mathcal{A}'$ by hypothesis.  Hence $E(T, \lambda)$ is an invariant
subspace of $\mathcal{A}$, in contradiction to the assumption that
$\mathcal{A}$ is irreducible.  We conclude that $E(T, \lambda) = V$.

	An alternative approach is to say that if $T$ lies in
$\mathcal{A}'$, then $T$ generates a field extension of $k$ of finite
degree (since $\mathcal{A}' \subseteq \mathcal{L}(V)$ has finite
dimension over $k$), and that this field extension should be equal
to $k$ when $k$ is algebraically closed.

	Now let us consider the case of a vector space $V$ over the
real numbers ${\bf R}$.  Suppose as before that $\mathcal{A}$ is an
algebra of operators on $V$ which is irreducible.  One can start with
Proposition \ref{irreducibility of mathcal{A} implies ....} and apply
a famous result of Frobenius (as in \cite{Albert}) to conclude that
$\mathcal{A}'$ either consists of real multiples of the identity, or
that $\mathcal{A}'$ is isomorphic to the complex numbers or to the
quaternions, as algebras over the real numbers.  Here we shall
describe a different (also classical) analysis, under the additional
assumption that $V$ is equipped with an inner product $\langle \cdot,
\cdot \rangle$ and that $\mathcal{A}$ is a $*$-algebra.  A broader
discussion can be found in \cite{Curtis1, Harvey}, which also allows
for nonassociative algebras.

\beginlemma
\label{T = T_1 + T_2, T_1^* = T_1, T_2^* = -T_2, T_i in mathcal{A}'}
If $T \in \mathcal{A}'$, then $T = T_1 + T_2$, where $T_1^* = T_1$,
$T_2^* = - T_2$, and $T_1, T_2 \in \mathcal{A}'$.
\end{lemma}

	As usual, one can take $T_1 = (T + T^*)/2$ and $(T - T^*)/2$
to get $T_1$, $T_2$ so that $T = T_1 + T_2$, $T_1^* = T_1$, and
$T_2^* = - T_2$.  Because $\mathcal{A}'$ is a $*$-algebra, we also
have that $T_1, T_2 \in \mathcal{A}'$ when $T \in \mathcal{A}_1$.

\beginlemma
\label{T in mathcal{A}', mathcal{A} irreducible, T^* = T}
If $T \in \mathcal{A}'$ is symmetric, $T^* = T$, then $T = \alpha \, I$
for some real number $\alpha$.
\end{lemma}

	For this we use the fact that $T$ is diagonalizable, since it
is symmetric.  We really only need that $T$ has an eigenvalue, as in
the setting of Proposition \ref{k alg closed, irreducibility, imply
mathcal{A}' = scalars I}.  The eigenspace is an invariant subspace for
$\mathcal{A}$, for the same reason as before, and hence the eigenspace
is all of $V$, by the irreducibility of $V$.  This says exactly that
$T$ is a scalar multiple of the identity.

	If $R$ and $T$ are elements of $\mathcal{A}'$, then $(R^* \, T
+ T^* \, R)/2$ is an element of $\mathcal{A}'$, since $\mathcal{A}'$
is a $*$-algebra, and it is also symmetric as an operator on $V$.
Hence it is a real number times the identity operator, by Lemma \ref{T
in mathcal{A}', mathcal{A} irreducible, T^* = T}.  Let $(R,T)$ denote
this real number, so that
\begin{equation}
\label{(R^* T + T^* R)/2 = (R, T) I}
	(R^* \, T + T^* \, R)/2 = (R, T) \, I.
\end{equation}
It is easy to see that $(R,T)$ is linear in $R$ and in $T$, and
symmetric in $R$ and $T$.

	If $R = T$, then we can rewrite (\ref{(R^* T + T^* R)/2 = (R,
T) I}) as
\begin{equation}
\label{R^* R = (R, R) I}
	R^* \, R = (R, R) \, I.
\end{equation}
If $v$ is any vector in $V$, then 
\begin{equation}
\label{langle R(v), R(v) rangle = ... = (R, R) langle v, v rangle} 
	\langle R(v), R(v) \rangle = 
	\langle (R^* \, R)(v), v \rangle = (R, R) \, \langle v, v \rangle.
\end{equation}
This implies that $(R,R)$ is always nonnegative, and that it is equal
to $0$ exactly when $R = 0$.  Thus $(R,T)$ defines an inner product
on $\mathcal{A}'$.  If $(R, R) > 0$, we also get that $R$ is
invertible.

	Next, if $R_1, R_2 \in \mathcal{A}'$, then
\begin{equation}
	(R_1 \, R_2, R_1 \, R_2) = (R_1, R_1) \, (R_2, R_2).
\end{equation}
This is because
\begin{eqnarray}
	(R_1 \, R_2)^* \, (R_1 \, R_2) = R_2^* \, R_1^* \, R_1 \, R_2
  		& = & (R_1, R_1) \, R_2^* \, R_2		\\
 	 	& = & (R_1, R_1) \, (R_2, R_2) \, I.	\nonumber
\end{eqnarray}

	Now let us focus somewhat on antisymmetric operators in
$\mathcal{A}'$.  If $R \in \mathcal{A}'$ is antisymmetric, then
\begin{equation}
\label{R^2 = -(R,R) I}
	R^2 = -(R,R) \, I.
\end{equation}
If $R_1, R_2 \in \mathcal{A}'$ are antisymmetric and $(R_1, R_2) = 0$,
then $R_1^* \, R_2 + R_2^* \, R_1 = 0$, which reduces to
\begin{equation}
\label{R_1 R_2 = - R_2 R_1}
	R_1 \, R_2 = - R_2 \, R_1.
\end{equation}
In this case we also have that
\begin{equation}
	(R_1 \, R_2)^* = R_2^* \, R_1^* = R_2 \, R_1 = - R_1 \, R_2,
\end{equation}
so that $R_1 \, R_2$ is also antisymmetric.

\beginlemma
\label{R_1, R_2, and products thereof}
Suppose that $R_1$ and $R_2$ are antisymmetric elements of
$\mathcal{A}'$ such that $(R_1, R_2) = 0$.  Then
\begin{equation}
	R_1 \, (R_1 \, R_2) = - (R_1 \, R_2) \, R_1, \quad
	R_2 \, (R_1 \, R_2) = - (R_1 \, R_2) \, R_2,
\end{equation}
these elements of $\mathcal{A}'$ are antisymmetric, and
\begin{equation}
	(R_1, R_1 \, R_2) = (R_2, R_1 \, R_2) = 0.
\end{equation}
\end{lemma}

	This is not hard to check.

	Assume for the moment that $T_1$ and $T_2$ are nonzero
antisymmetric elements of $\mathcal{A}'$ such that $(T_1, T_2) = 0$,
and that the set of antisymmetric elements of $\mathcal{A}'$ is not
spanned by $T_1$, $T_2$, and $T_1 \, T_2$.  Then there is another
nonzero antisymmetric element $U$ of $\mathcal{A}'$ such that
\begin{equation}
	(U, T_1) = (U, T_2) = (U, T_1 \, T_2) = 0.
\end{equation}
We know from the observations above that $T_1 \, T_2$ is
antisymmetric, and hence $(U, T_1 \, T_2) = 0$ implies that
\begin{equation}
\label{U (T_1 T_2) = - (T_1 T_2) U}
	U \, (T_1 \, T_2) = - (T_1 \, T_2) \, U,
\end{equation}
as in (\ref{R_1 R_2 = - R_2 R_1}).  Similarly,
\begin{equation}
	U \, T_1 = - T_1 \, U \hbox{ and }
	U \, T_2 = - T_2 \, U,
\end{equation}
since $(U, T_1) = (U, T_2) = 0$.  These two equations imply
that
\begin{eqnarray}
	U \, (T_1 \, T_2) & = & (U \, T_1) \, T_2
	   =  - (T_1 \, U) \, T_2 				\\
	  & = & - T_1 \, (U \, T_2)
	   = - T_1 \, (- T_2 \, U) = (T_1 \, T_2) \, U.	     \nonumber
\end{eqnarray}
This together with (\ref{U (T_1 T_2) = - (T_1 T_2) U})
lead to
\begin{equation}
	U \, (T_1 \, T_2) = 0.
\end{equation}
This is a contradiction, because $T_1$, $T_2$, and $U$ are assumed to
be nonzero elements of $\mathcal{A}'$, and hence are invertible, so
that the product cannot be $0$.

	The conclusion of this is that if $T_1$ and $T_2$ are nonzero
antisymmetric elements of $\mathcal{A}'$, then the set of
antisymmetric elements of $\mathcal{A}'$ is spanned by $T_1$, $T_2$,
and $T_3$.  

	Thus we obtain three possibilities for $\mathcal{A}'$.  The
first is that there are no nonzero antisymmetric elements of
$\mathcal{A}'$.  In other words, all elements of $\mathcal{A}'$ are
symmetric, because of Lemma \ref{T = T_1 + T_2, T_1^* = T_1, T_2^* =
-T_2, T_i in mathcal{A}'}.  From Lemma \ref{T in mathcal{A}',
mathcal{A} irreducible, T^* = T} it follows that $\mathcal{A}'$
consists exactly of multiples of the identity operator by real
numbers.

	The second possibility is that there is a nonzero
antisymmetric element $J$ of $\mathcal{A}'$, and that all other
antisymmetric elements of $\mathcal{A}'$ are multiples of $J$ by
real numbers.  

	The third possibility is that $\mathcal{A}'$ contains nonzero
antisymmetric elements $T_1$, $T_2$ such that $(T_1, T_2) = 0$.  In
this case $T_1 \, T_2$ is also a nonzero antisymmetric element of
$\mathcal{A}'$, $T_1 \, T_2$ is orthogonal to $T_1$ and $T_2$ relative
to the inner product $(\cdot, \cdot)$ on $\mathcal{A}'$, and all
antisymmetric elements of $\mathcal{A}'$ are linear combinations of
$T_1$, $T_2$, and $T_1 \, T_2$.

	These are the only possibilities for $\mathcal{A}'$.  In other
words, the dimension of the real vector space of antisymmetric
elements of $\mathcal{A}'$ is either $0$, $1$, or greater than or
equal to $2$, and these cases correspond exactly to the possibilities
just described.  In particular, if the vector space of antisymmetric
elements of $\mathcal{A}'$ has dimension greater than or equal to
$2$, then the dimension is in fact equal to $3$.

	In the first situation, where $\mathcal{A}'$ consists of only
multiples of the identity, it follows from Theorem \ref{mathcal{A}'' =
mathcal{A} if mathcal{A} is very nice} that $\mathcal{A}$ is equal to
the algebra $\mathcal{L}(V)$ of all linear transformations on $V$.
Note that if $\mathcal{A} = \mathcal{L}(V)$, then $\mathcal{A}'$ is
equal to the set of multiples of the identity operator by real
numbers.

	In the second situation, there is a nonzero antisymmetric
element $J$ of $\mathcal{A}'$ such that all other antisymmetric
elements of $\mathcal{A}'$ are equal to real numbers times $J$.  We
may as well assume that $(J, J) = 1$, since we can multiply $J$ by a
positive real number to obtain this.  From (\ref{R^2 = -(R,R) I}) we
obtain that
\begin{equation}
\label{J^2 = -I}
	J^2 = -I.
\end{equation}
We also have that
\begin{equation}
\label{langle J(v), J(v) rangle = langle v, v rangle}
	\langle J(v), J(v) \rangle = \langle v, v \rangle
\end{equation}
for all $v \in V$, as in (\ref{langle R(v), R(v) rangle = ... = (R, R)
langle v, v rangle}).  (Notice that any two of the conditions $J^* =
-J$, (\ref{J^2 = -I}), and (\ref{langle J(v), J(v) rangle = langle v,
v rangle}) inplies the third.)

	In general, if $J$ is a linear operator on a real vector space
that satisfies $J^2 = -I$, then $J$ defines a \emph{complex structure}
on the vector space.  That is, one can use $J$ as a definition for
multiplication by $i$ on the vector space (and this choice satisfies
the requirements of a complex vector space).  A real linear operator
on the vector space that commutes with $J$ is the same as a
complex-linear operator on the vector space, with respect to this
complex structure.  If the vector space is equipped with a (real)
inner product, then the additional condition of antisymmetry, or
equivalently (\ref{langle J(v), J(v) rangle = langle v, v rangle}), is
a compatibility property of the complex structure with the inner
product.

	Here we have such an operator $J$, and $\mathcal{A}'$ is
spanned by the identity operator $I$ and $J$.  By Theorem
\ref{mathcal{A}'' = mathcal{A} if mathcal{A} is very nice},
$\mathcal{A} = \mathcal{A}''$ is the set of linear operators on $V$
which commute with $J$, i.e., which are complex-linear with respect to
the complex structure defined by $J$.

	Conversely, one could start with a linear operator $J$ on $V$
such that (\ref{J^2 = -I}) holds, and then define $\mathcal{A}$ to be
the set of linear operators on $V$ which are complex-linear with
respect to the complex structure defined by $J$.  It is easy to see
that $\mathcal{A}$ is then an algebra.  If $J$ is compatible with
the inner product $\langle \cdot, \cdot \rangle$, in the sense of
(\ref{langle J(v), J(v) rangle = langle v, v rangle}) (or,
equivalently, by being antisymmetric), then one can check that
$\mathcal{A}$ is a $*$-algebra.  Also, $\mathcal{A}'$ is equal to
the span of $J$ and the identity operator in this case.

	The third situation is analogous to the second one, except for
having a ``quaternionic structure'' on $V$ rather than a complex
structure, corresponding to the larger family of operators in
$\mathcal{A}'$.  Let us briefly review some aspects of the
quaternions.

	Just as one might think of a complex number as being a
number of the form $a + b \, i$, where $i^2 = -1$ and $i$ commutes
with all real numbers, one can think of a quaternion as being
something which can be expressed in a unique way as
\begin{equation}
	a + b \, i + c \, j + d \, k,
\end{equation}
where $a$, $b$, $c$, and $d$ are real numbers, and $i$, $j$, and $k$
are special quaternions analogous to $i$ in the complex numbers.
Specifically, $i$, $j$, and $k$ commute with real numbers and satisfy
\begin{equation}
	i^2 = j^2 = k^2 = -1, \quad i \, j = - j \, i = k.
\end{equation}
As a result of the latter equations, one also has 
\begin{equation}
	i \, k = - k \, i, \quad j \, k = - k \, j.
\end{equation}

	If $x = a + b \, i + c \, j + d \, k$ is a quaternion,
then its \emph{conjugate} $\overline{x}$ is defined by
\begin{equation}
	\overline{x} = a - b \, i - c \, j - d \, k,
\end{equation}
in analogy with complex conjugation.  As in the setting of complex
numbers, one can check that
\begin{equation}
	x \, \overline{x} = \overline{x} \, x = a^2 + b^2 + c^2 + d^2,
\end{equation}
where the right side is always a real number.  One defines the
``norm'' or ``modulus'' $|x|$ of $x$ to be $\sqrt{x \, \overline{x}}$.
This is the same as the usual Euclidean norm of $(a, b, c, d) \in {\bf
R}^4$.  Clearly the conjugate of the conjugate of $x$ is $x$ again,
and one can verify that
\begin{equation}
	\overline{(xy)} = \overline{y} \, \overline{x}
\end{equation}
for all quaternions $x$ and $y$.  As a result, 
\begin{equation}
	|x y|^2 =  (x y) \, \overline{(xy)}
			= x \, y \, \overline{y} \, \overline{x}
			= |x|^2 \, |y|^2
\end{equation}
for all quaternions $x$, $y$, i.e., $|x y| = |x| \, |y|$.

	In the third situation above, $\mathcal{A}'$ is isomorphic to
the quaternions in a natural sense.  One can think of the identity
operator in $\mathcal{A}'$ as corresponding to the real number $1$ in
the quaternions, and in general the real multiples of the identity
operator in $\mathcal{A}'$ correspond to the real numbers inside the
quaternions.  The antisymmetric operators in $\mathcal{A}'$
correspond to the ``imaginary'' quaternions, which are the
quaternions of the form $b \, i + c \, j + d \, k$.  There are not
necessarily specific counterparts for $i$, $j$, and $k$ in
$\mathcal{A}'$, but one can use any antisymmetric operators $R_1$,
$R_2$, and $R_3$ such that $R_3 = R_1 \, R_2$, $(R_1, R_1) = (R_2,
R_2) = 1$, and $(R_1, R_2) = 0$ (as in the earlier discussion).

	This leads to a one-to-one correspondence between
$\mathcal{A}'$ and the quaternions that preserves sums and products.
Similarly, the conjugation operation $x \mapsto \overline{x}$ on the
quaternions matches with the adjoint operation $T \mapsto T^*$ on
$\mathcal{A}'$.  The norm on $\mathcal{A}'$ coming from the inner
product $(\cdot,\cdot)$ matches with the norm on quaternions.

	In this manner the operators in $\mathcal{A}'$ define a
quaternionic structure on $V$.  That is, a complex structure on a
real vector space gives a way to have complex numbers operate on
vectors in the vector space, and now we have a way for the
quaternions to operate on a vector space, in a way that extends the
usual scalar multiplication by real numbers.  This quaternionic
structure is also compatible with the inner product on $V$, in the
sense that a quaternion with norm $1$ corresponds to an operator on
$V$ which preserves the inner product there.  This uses (\ref{R^* R =
(R, R) I}).

	Theorem \ref{mathcal{A}'' = mathcal{A} if mathcal{A} is very
nice} says that $\mathcal{A} = \mathcal{A}''$, so that $\mathcal{A}$
consists exactly of the real-linear operators on $V$ that commute with
the operators that come from the quaternionic structure (i.e., the
operators in $\mathcal{A}'$ in this case).  In other words,
$\mathcal{A}$ consists exactly of the operators which are
``quaternionic-linear'' on $V$ with respect to this quaternionic
structure coming from $\mathcal{A}'$, analogous to the second
situation and complex-linear operators with respect to the complex
structure on $V$ that arose there.

\section{Division algebras of operators}
\label{Division algebras of operators}
\setcounter{equation}{0}

\begindefinition
\label{def of division algebra (of operators)}
Let $\mathcal{B}$ be an algebra of operators on a vector space $V$
over a field $k$.  Suppose that for each $T$ in $\mathcal{B}$, either
$T = 0$ or $T$ is an invertible operator on $V$ (so that $T^{-1}$ lies
in $\mathcal{B}$, by Proposition \ref{inverse in the algebra when it
exists}).  In this case $\mathcal{B}$ is said to be a \emph{division
algebra of operators}.\index{division algebra of operators}
\end{definition}

	For the rest of this section we assume that $k$ is a field,
$V$ is a vector space over $k$, and $\mathcal{B}$ is a division
algebra of operators on $V$.

	If $v_1, \ldots, v_m$ are vectors in $V$, then we say that
they are
\emph{$\mathcal{B}$-independent}\index{B-independent@$\mathcal{B}$-independent}
if
\begin{equation}
	T_1(v_1) + \cdots + T_m(v_m) = 0
\end{equation}
implies that 
\begin{equation}
	T_1 = \cdots = T_m = 0
\end{equation}
for all $T_1, \ldots, T_m$ in $\mathcal{B}$.  The
\emph{$\mathcal{B}$-span}\index{B-span@$\mathcal{B}$-span} of $v_1,
\ldots, v_m$ is defined to be the subspace of $V$ consisting of
vectors of the form
\begin{equation}
\label{T_1(v_1) + cdots + T_m(v_m)}
	T_1(v_1) + \cdots + T_m(v_m),
\end{equation}
where $T_1, \ldots, T_m$ lie in $\mathcal{B}$.  Thus
$\mathcal{B}$-independence of $v_1, \ldots, v_m$ is equivalent to
saying that if a vector can be represented as (\ref{T_1(v_1) + cdots +
T_m(v_m)}) with $T_1, \ldots, T_m \in \mathcal{B}$, then this
representation is unique.

	Note that a single nonzero vector in $V$ is
$\mathcal{B}$-independent.  This uses the assumption that an operator
in $\mathcal{B}$ is either $0$ or invertible.

\beginlemma
\label{independent vectors, adding one more}
Assume that $v_1, \ldots, v_m$ are $\mathcal{B}$-independent
vectors in $V$.  If $w$ is a vector in $V$ that does not lie in
the $\mathcal{B}$-span of $v_1, \ldots, v_m$, then the vectors
$v_1, \ldots, v_m, w$ are $\mathcal{B}$-independent.
\end{lemma}

	Indeed, suppose that $T_1, \ldots, T_m, U$ are operators
in $\mathcal{B}$ such that
\begin{equation}
\label{T_1(v_1) + cdots + T_m(v_m) + U(w) = 0}
	T_1(v_1) + \cdots + T_m(v_m) + U(w) = 0.
\end{equation}
We would like to show that $T_1 = \cdots = T_m = U = 0$.  If $U = 0$,
then this reduces to the $\mathcal{B}$-independence of $v_1, \ldots,
v_m$.  If $U \ne 0$, then $U$ is invertible, by our assumptions on
$\mathcal{B}$, and we can rewrite (\ref{T_1(v_1) + cdots + T_m(v_m) +
U(w) = 0}) as
\begin{equation}
	w = - (U^{-1} \circ T_1)(v_1) - \cdots - (U^{-1} \circ T_m)(v_m).
\end{equation}
This implies that $w$ is in the $\mathcal{B}$-span of $v_1, \ldots,
v_m$, contrary to hypothesis, and Lemma \ref{independent vectors,
adding one more} follows.

\begincorollary
\label{mathcal{B}-independent vectors which mathcal{B}-span V}
There are vectors $u_1, \ldots, u_\ell$ in $V$ which are
$\mathcal{B}$-independent and whose $\mathcal{B}$-span is equal to
$V$.
\end{corollary}

	To see this, one can start with a single nonzero vector $u_1$
in $V$.  If the $\mathcal{B}$-span of this vector is all of $V$, then
we stop.  Otherwise, we choose a vector $u_2$ which is not in the
$\mathcal{B}$-span of $u_1$.  Lemma \ref{independent vectors, adding
one more} implies that $u_1, u_2$ are $\mathcal{B}$-independent.  If
the $\mathcal{B}$-span of $u_1, u_2$ is all of $V$, then we stop, and
otherwise we choose $u_3$ to be a vector in $V$ not in the
$\mathcal{B}$-span of $u_1, u_2$.  Proceeding in this manner, one
eventually gets a collection of vectors $u_1, \ldots, u_\ell$ in $V$
as in the corollary.

\begincorollary
\label{dim V is divisible by dim mathcal{B}}
The dimension of $V$ is divisible by the dimension of $\mathcal{B}$
(as vector spaces over $k$).
\end{corollary}

	Specifically, the dimension of $V$ is equal to $\ell$ times
the dimension of $\mathcal{B}$, if $\ell$ is as in the previous
corollary.

\begincorollary
\label{bases adapted to invariant subspaces}
Suppose that $W$ is a nonzero vector subspace of $V$ which is
invariant under $\mathcal{B}$.  There are vectors $u_1, \ldots,
u_\ell$ in $V$ and a positive integer $m$ such that $u_1, \ldots,
u_\ell$ are $\mathcal{B}$-independent, the $\mathcal{B}$-span of $u_1,
\ldots, u_\ell$ is equal to $V$, $u_j$ lies in $W$ exactly when $j \le
m$, and $W$ is equal to the $\mathcal{B}$-span of the $u_j$'s, $j \le
m$.
\end{corollary}

	This can be established in the same manner as the Corollary
\ref{mathcal{B}-independent vectors which mathcal{B}-span V}.  The
main difference is that in the beginning of the construction one
chooses the $u_j$'s in $W$, for as long as one can.  One stops doing
this exactly when there are enough $u_j$'s so that their
$\mathcal{B}$-span is $W$.  One then continues as before, with
additional $u_j$'s as needed so that the $\mathcal{B}$-span of the
total collection of vectors is $V$.

\begincorollary
\label{invariant subspaces have invariant complements}
If $W$ is a nonzero vector subspace of $V$ which is invariant under
$\mathcal{B}$ and which is not equal to $V$, then there is a vector
subspace $Z$ of $V$ which is invariant under $\mathcal{B}$ and which
is a complement of $W$ (so that $W \cap Z = \{0\}$ and $W + Z = V$).
\end{corollary}

	For this one can take $Z$ to be the $\mathcal{B}$-span of
$u_{m+1}, \ldots, u_\ell$, where $u_1, \ldots, u_\ell$ is as in 
Corollary \ref{bases adapted to invariant subspaces}.

	Now let us consider $\mathcal{B}'$.  Let $u_1, \ldots, u_\ell$
be vectors in $V$ which are $\mathcal{B}$-independent and whose
$\mathcal{B}$-span is equal to $V$, as in Corollary
\ref{mathcal{B}-independent vectors which mathcal{B}-span V}.  If $R
\in \mathcal{B}'$, then $R$ is determined uniquely by $R(u_1), \ldots,
R(u_\ell)$.  For if $T \in \mathcal{B}$, then $R(T(u_j)) = T(R(u_j))$
for each $j$, and this and linearity determine $R$ on all of $V$.

	Conversely, if $w_1, \ldots, w_\ell$ are arbitrary vectors in
$V$, then there is an $R \in \mathcal{B}'$ so that $R(u_j) = w_j$ for
each $j$.  One can simply define $R$ on $V$ by
\begin{equation}
	R(T_1(u_1) + \cdots + T_\ell(u_\ell))
		= T_1(w_1) + \cdots + T_\ell(w_\ell),
\end{equation}
for all $T_1, \ldots, T_\ell \in \mathcal{B}$.  It is easy to
check that this gives a linear operator $R$ on $V$ which lies
in $\mathcal{B}'$.

\beginlemma
\label{dimension of mathcal{B}' as a vector space}
The dimension of $\mathcal{B}'$ is equal to the square of the
dimension of $V$ divided by the dimension of $\mathcal{B}$, as vector
spaces over $k$.
\end{lemma}

	This follows from the preceding description of $\mathcal{B}'$.

\beginproposition
\label{mathcal{B}'' = mathcal{B} when mathcal{B} is a division alg of ops}
$\mathcal{B}'' = \mathcal{B}$.
\end{proposition}

	This can be verified directly, using our characterization of
$\mathcal{B}'$.  Let us mention another argument, in terms of the
methods of Section \ref{Algebras of linear operators}.  It is enough
to show that $\mathcal{B}$ is a very nice algebra of operators, by
Theorem \ref{mathcal{A}'' = mathcal{A} if mathcal{A} is very nice} in
Subsection \ref{Very nice algebras of operators}.  Using Lemma
\ref{complementary subspaces and commuting operators} in Subsection
\ref{Nice algebras of operators} and Corollary \ref{invariant
subspaces have invariant complements} one obtains that $\mathcal{B}$
is a nice algebra of operators (Definition \ref{def of nice algebras
of operators}).  To say that $\mathcal{B}$ is very nice means
(Definition \ref{def of very nice algebras of operators}) that the
expanded algebras $\mathcal{B}_n$ defined in Subsection \ref{Expanding an
algebra} are all nice.  By construction, the expanded algebras
$\mathcal{B}_n$ are all division algebras of operators, and hence they
are nice for the same reason that $\mathcal{B}$ is.

	Next we discuss some analogues of the results in Subsection
\ref{Some other approaches}.  Let us make the standing assumption
for the rest of the section that 
\begin{equation}
\label{standing assumption about mathcal{A}}
	\hbox{$\mathcal{A}$ is an algebra of operators on $V$ such that }
		\mathcal{A} \subseteq \mathcal{B}'.
\end{equation}

\beginlemma
\label{criterion for mathcal{A} = mathcal{B}', 0}
If the dimension of $V$ is equal to the dimension of $\mathcal{B}$,
as vector spaces over $k$, then $\mathcal{A} = \mathcal{B}'$ if and
only if $\mathcal{A}$ is irreducible.
\end{lemma}

	The assumption that the dimension of $V$ is equal to the
dimension of $\mathcal{B}$ implies that $V$ is equal to the
$\mathcal{B}$-span of any nonzero vector in $V$.  Using this and the
reformulation of irreducibility in Lemma \ref{reformulation of
irreducibility} it is easy to check that $\mathcal{A} = \mathcal{B}'$
when $\mathcal{A}$ is irreducible.  For the converse, notice that
$\mathcal{B}'$ is always irreducible, without any restriction on the
dimension of $V$.

\beginlemma
\label{criterion for mathcal{A} = mathcal{B}', 1}
If $\mathcal{A}$ is irreducible and there is a nonzero operator in
$\mathcal{A}$ whose image is contained in the $\mathcal{B}$-span of a
single vector, then $\mathcal{A} = \mathcal{B}'$.
\end{lemma}

	The assumptions in this lemma are also necessary for
$\mathcal{A}$ to be equal to $\mathcal{B}'$.

	For each nonzero vector $u$ in $V$, consider the vector space
of linear operators on $V$ which lie in $\mathcal{B}'$ and have image
contained in the $\mathcal{B}$-span of $u$.  We would like to show
that this vector space of linear operators is contained in
$\mathcal{A}$.  If we can do this for any nonzero $u$ in $V$, then it
follows that $\mathcal{A} = \mathcal{B}'$.

	Fix $u \ne 0$ in $V$.  Observe first that there is a nonzero
linear operator on $V$ which lies in $\mathcal{A}$ and has image
contained in the $\mathcal{B}$-span of $u$.  The hypothesis of the
lemma says that this is true for at least one $u_1$, and one can get
any other $u$ by composing with an operator in $\mathcal{A}$ that
sends $u_1$ to $u$.  Such an operator exists, by irreducibility.

	As a vector space over $k$, the collection of operators in
$\mathcal{B}'$ whose image is contained in the $\mathcal{B}$-span of
$u$ has dimension equal to the dimension of $V$.  This can be derived
from the earlier description of $\mathcal{B}'$.  Thus it is enough to
show that the collection of operators in $\mathcal{A}$ with image
contained in the $\mathcal{B}$-span of $u$ has dimension which is at
least the dimension of $v$.  Let $R$ be a nonzero operator in
$\mathcal{A}$ whose image is contained in, and hence equal to, the
$\mathcal{B}$-span of $u$, and let $z$ be a vector in $V$ such that
$R(z) = u$.  For each vector $v$ in $V$ there is a linear operator $T$
in $\mathcal{A}$ such that $T(v) = z$.  For this $T$ we have that $R
\circ T$ lies in $\mathcal{A}$, the image of $R \circ T$ is contained
in the $\mathcal{B}$-span of $u$, and $(R \circ T)(v) = u$.  Because
we can do this for each $v$ in $V$, it is not hard to see that the
dimension of the collection of linear operators in $\mathcal{A}$ with
image contained in the $\mathcal{B}$-span of $u$ has dimension at
least the dimension of $V$, as desired.

\beginlemma
\label{criterion for mathcal{A} = mathcal{B}', 2}
Suppose that $\mathcal{A}$ is irreducible, and that for every four
vectors $v_1$, $v_2$, $w_1$, $w_2$ in $V$ such that $v_1$, $v_2$
are $\mathcal{B}$-independent there is a $T$ in $\mathcal{A}$
which satisfies $T(v_1) = w_1$, $T(v_2) = w_2$.  Then $\mathcal{A}
= \mathcal{B}'$.
\end{lemma}

	The assumption of irreducibility of $\mathcal{A}$ follows from
the second condition when the dimension of $V$ is larger than the
dimension of $\mathcal{B}$.  When the two dimensions are equal, the
second condition is vacuous, and one could just as well apply Lemma
\ref{criterion for mathcal{A} = mathcal{B}', 0}.

	To prove Lemma \ref{criterion for mathcal{A} = mathcal{B}',
2}, it suffices to show that there is a nonzero operator in
$\mathcal{A}$ whose image is contained in the $\mathcal{B}$-span of a
single vector, by Lemma \ref{criterion for mathcal{A} = mathcal{B}',
1}.  Let $R$ be any nonzero operator in $\mathcal{A}$.  If the image
of $R$ is not contained in the $\mathcal{B}$-span of a single vector,
then there exist $x_1$, $x_2$ in $V$ such that $R(x_1)$, $R(x_2)$ are
$\mathcal{B}$-independent.  By hypothesis, there is an operator $T$ in
$\mathcal{A}$ such that $T(R(x_1)) \ne 0$ and $T(R(x_2)) = 0$.  Thus
$T \circ R$ lies in $\mathcal{A}$, is nonzero, and the dimension of
its image strictly less than that of $R$.  If the image of $T \circ R$
is contained in the $\mathcal{B}$-span of a single vector, then we are
finished, and otherwise we can repeat this process until such an
operator is obtained.  This proves the lemma.

\beginlemma
\label{criterion for mathcal{A} = mathcal{B}', 3}
Suppose that $\mathcal{A}' = \mathcal{B}$ and $\mathcal{A}_n \subseteq
\mathcal{L}(V^n)$ (defined in Subsection \ref{Expanding an algebra}) is a
nice algebra of operators when $n = 2$.  Then $\mathcal{A} =
\mathcal{B}' = \mathcal{A}''$.
\end{lemma}

	The hypotheses of this lemma imply those of Lemma
\ref{criterion for mathcal{A} = mathcal{B}', 2}.

\section{Group representations and their algebras}
\label{Group representations and their algebras}
\setcounter{equation}{0}

	Throughout this section, $G$ will be a finite group, and
$k$ will be a field.

\beginproposition
\label{one half of Schur's lemma}
Let $V_1$ and $V_2$ be vector spaces over $k$, and let $\rho^1$,
$\rho^2$ be representations of $G$ on $V_1$, $V_2$.  Suppose that $T :
V_1 \to V_2$ is a linear mapping which
\emph{intertwines}\index{interwining mapping (with respect to two
representations of a group)} the representations $\rho^1$, $\rho^2$,
in the sense that
\begin{equation}
	T \circ \rho^1_x = \rho^2_x \circ T
\end{equation}
for all $x$ in $G$.  If the representations $\rho^1$, $\rho^2$ are
\emph{irreducible}, then either $T$ is $0$, or $T$ is a one-to-one
mapping from $V_1$ onto $V_2$, and the representations $\rho^1$,
$\rho^2$ are isomorphic to each other.
\end{proposition}

	This is one half of \emph{Schur's lemma}.  (The other half
came up in Section \ref{section on irreducibility}.)

	Let $V_1$, $V_2$, etc., be as in the statement above.
Consider the kernel of $T$, which is a vector subspace of $V_1$.  From
the intertwining property it is easy to see that the kernel of $T$ is
invariant under the representation $\rho^1$.  The irreducibility of
$\rho^1$ then implies that the kernel is either all of $V_1$ or the
subspace consisting of only the zero vector.  In the first case $T =
0$, and we are finished.  In the second case we get that $T$ is
one-to-one.

	Now consider the image of $T$ in $V_2$.  This is a vector
subspace that is invariant under $\rho^2$, because of the intertwining
property.  The irreducibility of $\rho^2$ implies that the image of
$T$ is either the subspace of $V_2$ consisting of only the zero
vector or all of $V_2$.  This is the same as saying that $T$ is
either equal to $0$ or it maps $V_1$ onto $V_2$.  Combining this
with the conclusion of the preceding paragraph, we obtain that $T$
is either equal to $0$, or that it is one-to-one and maps $V_1$
onto $V_2$.  This proves the proposition.

	For the rest of this section, let us assume that
\begin{eqnarray}
  && \hbox{either the characteristic of $k$ is $0$, or it is positive and}
								\\
  && \hbox{does not divide the number of elements of the group $G$.}
								\nonumber
\end{eqnarray}

	Let $V$ be a vector space over $k$, and let $\rho$ be a
representation of $G$ on $V$.  As in Lemma \ref{decomp into
irreducible pieces} in Subsection \ref{Reducibility, continued}, there is
an independent system of subspaces $W_1, \ldots, W_h$ of $V$ such that
the span of the $W_j$'s is equal to $V$, each $W_j$ is invariant under
$\rho$, and the restriction of $\rho$ to each $W_j$ is irreducible.
Let us make the following additional assumption:
\begin{eqnarray}
  && \hbox{for each $j$, $l$ such that $1 \le j, l \le h$}
		\hbox{ and $j \ne l$, the}			\\
  && \hbox{restriction of $\rho$ to $W_j$ is not isomorphic to the}
							\nonumber \\
  && \hbox{restriction of $\rho$ to $W_l$, as a representation of $G$.}
							\nonumber
\end{eqnarray}

\beginnotation
\label{notation about mathcal{A}(U)}
{\rm Suppose that $\mathcal{A}$ is an algebra of operators on a vector
space $V$.  If $U$ is a vector subspace of $V$ which is invariant
under $\mathcal{A}$ (Definition \ref{invariant subspaces of an algebra
of operators}), then we let
$\mathcal{A}(U)$\index{$A(U)$@$\mathcal{A}(U)$} denote the set of
operators on $U$ which are the restrictions to $U$ of the operators in
$\mathcal{A}$.  It is easy to see that $\mathcal{A}(U)$ is an algebra
of operators on $U$.  The commutant and double commutant of
$\mathcal{A}(U)$ in $\mathcal{L}(U)$ are denoted $\mathcal{A}(U)'$ and
$\mathcal{A}(U)''$.  }
\end{notation}

	In our present circumstances, we let $\mathcal{A}$ denote the
algebra of operators on $V$ generated by $\rho$ as in Subsection
\ref{Algebras from group representations}.  A subspace $U$ of $V$ is
invariant under $\mathcal{A}$ if and only if it is invariant under
$\rho$, and in this case $\mathcal{A}(U)$ is the same as the algebra
generated by the restrictions of the $\rho_x$'s to $U$, $x$ in $G$.

\begintheorem
\label{mathcal{A} from a group rep., no irred. comp's repeated} 
Under the preceding conditions, $\mathcal{A}$ consists of the linear
operators $T$ on $V$ such that $T(W_j) \subseteq W_j$ and the
restriction of $T$ to $W_j$ lies in $\mathcal{A}(W_j)$ for each
$j$, $1 \le j \le h$.
\end{theorem}

	It is clear from the definitions that each $T \in \mathcal{A}$
satisfies the properties described in the theorem.  The key point is
that in the converse, the restrictions of $T$ to the various $W_j$'s
can be chosen independently of each other.  To prove this, we use
Theorem \ref{mathcal{A}'' = mathcal{A} if mathcal{A} is very nice}.

	Because $W_1, \ldots, W_h$ is an independent system of
subspaces of $V$ which spans $V$, there is for each $j = 1, \ldots, h$
a natural linear operator $P_j$ on $V$ which projects $V$ onto $W_j$.
Specifically, $P_j$ is the linear operator which satisfies $P_j(u) =
u$ when $u \in W_j$ and $P_j(z) = 0$ when $z \in W_l$, $l \ne j$.
Thus $P_j \circ P_l = 0$ when $j \ne l$, and $\sum_{j=1}^h P_j = I$,
the identity operator on $V$.  We also have that each $P_j$ lies in
$\mathcal{A}'$, because the $W_j$'s are invariant subspaces.

\beginlemma
\label{S in mathcal{A}', j ne l, imply P_l circ S circ P_j = 0}
If $S \in \mathcal{A}'$ and $1 \le j, l \le h$, $j \ne l$, then $P_l
\circ S \circ P_j = 0$.
\end{lemma}

	Let $S$, $j$, and $l$ be given as in the lemma, and
set $S_{j,l} = P_l \circ S \circ P_j$.  Note that $S_{j,l}$
lies in $\mathcal{A}'$, since $S$, $P_j$, and $P_l$ do.  Thus
\begin{equation}
\label{S_{j,l} circ rho_x = rho_x circ S_{j,l}}
	S_{j,l} \circ \rho_x = \rho_x \circ S_{j,l}
\end{equation}
for all $x$ in $G$ (and this is equivalent to the statement that
$S_{j,l} \in \mathcal{A}'$, since $\mathcal{A}$ is generated by the
$\rho_x$'s).

	Let $R_{j,l}$ be the linear mapping from $W_j$ to $W_l$ which
is the restriction of $S_{j,l}$ to $W_j$.  Thus $R_{j,l}$ is
essentially the same as $S_{j,l}$, except that it is formally viewed
as a mapping from $W_j$ to $W_l$ rather than as an operator on $V$.
On $W_j$ and $W_l$ we have representations of $G$, namely, the
restrictions of $\rho$ to $W_j$ and $W_l$.  By hypothesis, these
representations are irreducible and not isomorphically equivalent.
The mapping $R_{j,l}$ intertwines these representations, by
(\ref{S_{j,l} circ rho_x = rho_x circ S_{j,l}}).  Therefore, $R_{j,l}
= 0$, by Proposition \ref{one half of Schur's lemma}.  This is
equivalent to $S_{j,l} = 0$, and the lemma follows.

\beginlemma
\label{description of mathcal{A}' in this case}
A linear operator $S$ on $V$ lies in $\mathcal{A}'$ if and
only if $S(W_j) \subseteq W_j$ and the restriction of $S$
to $W_j$ lies in $\mathcal{A}(W_j)'$ for each $j$, $1 \le j \le h$.
\end{lemma}

	If $S \in \mathcal{A}'$, then $S(W_j) \subseteq W_j$ for each
$j$ by Lemma \ref{S in mathcal{A}', j ne l, imply P_l circ S circ P_j
= 0}.  To say that $S \in \mathcal{A}'$ means that $Z \circ S = S
\circ Z$ for all $Z \in \mathcal{A}$, or
\begin{equation}
\label{(Z circ S)(v) = (S circ Z)(v)}
	(Z \circ S)(v) = (S \circ Z)(v)
\end{equation}
for all $Z \in \mathcal{A}$ and all vectors $v$ in $V$.  Because $V$
is spanned by $W_1, \ldots, W_h$, (\ref{(Z circ S)(v) = (S circ
Z)(v)}) holds for all $v$ in $V$ if and only if it holds for all
$v$ in $W_j$ for each $j = 1, \ldots, h$.  If $v \in W_j$, then
$S(v) \in W_j$ and $Z(v) \in W_j$ for all $Z \in \mathcal{A}$.
Let $S_j$ and $Z_j$ denote the restrictions of $S$ and $Z$ to $W_j$.
Then we get that (\ref{(Z circ S)(v) = (S circ Z)(v)}) holds for
all $v$ in $V$ and all $Z$ in $\mathcal{A}$ if and only if
\begin{equation}
	(Z_j \circ S_j)(v) = (S_j \circ Z_j)(v)
\end{equation}
for all $v$ in $W_j$, all $Z$ in $\mathcal{A}$, and all $j = 1,
\ldots, h$.  This is the same as saying that for every $j$ the
restriction of $S$ to $W_j$ commutes with the restriction of every
element of $\mathcal{A}$ to $W_j$, or, equivalently, with every
element of $\mathcal{A}(W_j)$.  This completes the proof of the lemma.

\beginlemma
\label{description of mathcal{A}'' in this situation}
A linear operator $T$ on $V$ lies in $\mathcal{A}''$ if and
only if $T(W_j) \subseteq W_j$ and the restriction of $T$
to $W_j$ lies in $\mathcal{A}(W_j)''$ for each $j$.
\end{lemma}

	We know that each $W_j$ is an invariant subspace of $T \in
\mathcal{A}''$ because each $P_j$ lies in $\mathcal{A}'$.  This is
essentially the same as applying Lemma \ref{invariant subspaces from
elements of the commutant} in Subsection \ref{Nice algebras of
operators}.  One can finish the argument in essentially the same way
as for Lemma \ref{description of mathcal{A}' in this case}.

	Theorem \ref{mathcal{A} from a group rep., no irred. comp's
repeated} can be derived from Lemma \ref{description of mathcal{A}''
in this situation}.  Alternatively, it is enough to show that the
operators $T$ on $V$ which satisfy the conditions in Theorem
\ref{mathcal{A} from a group rep., no irred. comp's repeated} should
lie in $\mathcal{A}$, and for this it is enough to show that they
commute with every element of $\mathcal{A}'$ (so that they then lie in
$\mathcal{A}''$, which is equal to $\mathcal{A}$ by the Theorem
\ref{mathcal{A}'' = mathcal{A} if mathcal{A} is very nice} in Subsection
\ref{Very nice algebras of operators}).  This is an easy consequence
of Lemma \ref{description of mathcal{A}' in this case}.

\section{Basic facts about decompositions}
\label{Basic facts about decompositions}
\setcounter{equation}{0}

	Throughout this section $G$ is a finite group, and $k$ is a
field which is assumed to have characteristic $0$ or positive
characteristic which does not divide the number of elements in $G$.

	Suppose that $\rho$ is a representation of $G$ on a vector
space $V$ over $k$.  As in Lemma \ref{decomp into irreducible pieces}
in Subsection \ref{Reducibility, continued}, there is an independent
system of subspaces $W_1, \ldots, W_h$ of $V$ such that $V =
\Span(W_1, \ldots, W_h)$, each $W_j$ is invariant under $\rho$, and
the restriction of $\rho$ to each $W_j$ is irreducible.  We do not ask
that the restrictions of $\rho$ to the $W_j$'s be isomorphically
distinct.

\beginlemma
\label{U in relation to the W_j's}
In addition to the preceding conditions, assume that $U$ is a nonzero
vector subspace of $V$ which is invariant under $\rho$ and for which
the restriction of $\rho$ to $U$ is irreducible.  Let $I$ denote the
set of integers $j$, $1 \le j \le h$, such that the restriction of
$\rho$ to $W_j$ is isomorphic to the restriction of $\rho$ to $U$.
Then $I$ is not empty, and
\begin{equation}
	U \subseteq \Span \{W_j : j \in I\}.
\end{equation}
\end{lemma}

	For each $j = 1, \ldots, h$, let $P_j : V \to V$ be the
projection onto $W_j$ that comes naturally from the independent system
of subspaces $W_1, \ldots, W_h$, i.e., $P_j(u) = u$ when $u \in W_j$
and $P_j(z) = 0$ when $z \in W_l$, $l \ne j$.  Thus $\sum_{j=1}^h P_j
= I$ and each $P_j$ commutes with the operators $\rho_x$, $x$ in $G$,
from the representation, because the $W_j$'s are invariant under the
representation.

	Define $T_j : U \to W_j$ for $j = 1, \ldots, h$ to be the
restriction of $P_j$ to $U$, now viewed as a mapping into $W_j$.
Notice that $T_j$ intertwines the representations on $U$, $W_j$ which
are the restrictions of $\rho$ to $U$, $W_j$, since $P_j$ commutes
with the operators $\rho_x$, $x$ in $G$.  Because the restrictions of
$\rho$ to $U$, $W_j$ are irreducible, we have that each $T_j$ is
either equal to $0$ or is a one-to-one linear mapping onto $W_j$, by
Proposition \ref{one half of Schur's lemma}.

	If $J$ denotes the set of $j$'s such that $T_j \ne 0$, then
\begin{equation}
	U \subseteq \Span \{W_j : j \in J\},
\end{equation}
simply because $P_j$ is equal to $0$ on $U$ for $j \not\in J$.  In
particular, $J$ is not equal to the empty set, since $U \ne \{0\}$.
We also have that $J \subseteq I$, since $T_j : U \to W_j$ defines an
isomorphism between the restriction of $\rho$ to $U$ and the
restriction of $\rho$ to $W_j$ when $j \in J$.  This proves Lemma
\ref{U in relation to the W_j's}.

	Now suppose that $Y_1, \ldots, Y_m$ are nonzero subspaces of
$V$ with properties analogous to those of $W_1, \ldots, W_h$, namely,
$Y_1, \ldots, Y_m$ is an independent system of subspaces of $V$ such
that $V = \Span(Y_1, \ldots, Y_m)$, each $Y_j$ is invariant under
$\rho$, and the restriction of $\rho$ to each $Y_j$ is irreducible.
In general it is not true that the $Y_j$'s have to be the same
subspaces of $V$ as the $W_l$'s, after a permutation of the indices.
For instance, if $\rho$ is the trivial representation, so that each
$\rho_x$ is equal to the identity mapping on $V$, then any subspace of
$V$ is invariant under $\rho$, and any $1$-dimensional subspace has
the property that the restriction of $\rho$ to it is irreducible.  If
the dimension of $V$ is strictly larger than $1$, then there are
numerous ways in which to decompose $V$ into $1$-dimensional
subspaces.

	However, it is true that the two decompositions of $V$ have to
be nearly the same in some respects.  Let $J$ be a nonempty subset of
$\{1, \ldots, h\}$ such that the restriction of $\rho$ to $W_j$ is
isomorphic to the restriction of $\rho$ to $W_l$ when $j, l \in J$,
and such that the restriction of $\rho$ to $W_j$ is not isomorphic to
the restriction of $\rho$ to $W_l$ when $j \in J$ and $l \not\in J$.
Let $L$ be the set of $l \in \{1, \ldots, m\}$ such that the
restriction of $\rho$ to $W_j$ is isomorphic to the restriction of
$\rho$ to $Y_l$ when $j \in J$.  Then
\begin{equation}
\label{Span {W_j : j in J} = Span {Y_l : l in L}}
	\Span \{W_j : j \in J\} = \Span \{Y_l : l \in L\}.
\end{equation}
This is because $Y_l \subseteq \Span \{W_j : j \in J\}$ when $l \in
L$, by Lemma \ref{U in relation to the W_j's}, and similarly
$W_j \subseteq \Span \{Y_l : l \in L\}$ when $j \in J$.

	Because of the isomorphisms, the $W_j$'s for $j \in J$ and
the $Y_l$'s for $l \in L$ all have the same dimension.  We also
have that
\begin{equation}
	\dim \Span \{W_j : j \in J\} =  \sum_{j \in J} \dim W_j
\end{equation}
and
\begin{equation}
	\dim \Span \{Y_l : l \in L\}  =  \sum_{l \in L} \dim Y_l,	
\end{equation}
since the $W_j$'s and $Y_l$'s form independent systems of subspaces.
Using (\ref{Span {W_j : j in J} = Span {Y_l : l in L}}) we obtain that
$J$ and $L$ have the same number of elements

	We can do this for all subsets $J$ of $\{1, \ldots, h\}$ with
the properties described above, and these subsets exhaust all of $\{1,
\ldots, h\}$.  The corresponding sets $L \subseteq \{1, \ldots, m\}$
exhaust all of $\{1, \ldots, m\}$ as well, because of Lemma \ref{U in
relation to the W_j's}.  To summarize, the same irreducible
representations of $G$ occur as restrictions of $\rho$ to the $W_j$'s
as occur as restrictions of $\rho$ to the $Y_l$'s, up to isomorphic
equivalence, and they occur the same number of times.  In particular,
we get that $h = m$.  Although the various subspaces do not have to
match up exactly, we do have (\ref{Span {W_j : j in J} = Span {Y_l : l
in L}}).

	Now let $Z$ be a vector space over $k$, and $\sigma$ be an
irreducible representation of $G$ on $Z$.  Let $F(G)$ denote the
vector space of $k$-valued functions on $G$.

\beginlemma
\label{irred. rep.'s isomorphic to a subrep. of the reg. rep.}
Suppose that $\lambda$ is a nonzero linear functional on $Z$, i.e., a
nonzero linear mapping from $Z$ into $k$.  For each $v$ in $Z$,
consider the function $f_v(y)$ on $G$ defined by
\begin{equation}
	f_v(y) = \lambda(\sigma_{y^{-1}}(v)).
\end{equation}
Define $U \subseteq F(G)$ by
\begin{equation}
	U = \{f_v(y) : v \in Z\}.
\end{equation}
Then $U$ is a vector subspace of $F(G)$ which is invariant under the
left regular representation (Subsection \ref{Representations of finite
groups}), and the mapping $v \mapsto f_v$ is a one-to-one linear
mapping from $Z$ onto $U$ which intertwines the representations
$\sigma$ on $Z$ and the restriction of the left regular representation
to $U$.  In particular, these two representations are isomorphic.
\end{lemma}

	Clearly $v \mapsto f_v$ is a linear mapping from $Z$ onto
$U$, so that $U$ is a vector subspace of $F(G)$.  Let us check that
the kernel of this mapping is the zero subspace of $Z$.  

	Suppose that $v$ is a vector in $Z$ such that $f_v(y)$ is the
zero function on $G$.  This is the same as saying that
$\lambda(\sigma_{y^{-1}}(v)) = 0$ for all $y$ in $G$.  Let $W$ be the
subspace of $Z$ which is spanned by $\sigma_y(v)$, $y \in G$.  Then
$W$ is invariant under the representation $\sigma$, because of the way
that it is defined.  If $v \ne 0$, then $W$ is not the zero subspace,
and indeed $v$ lies in $W$ since $\sigma_e(v) = v$.  The
irreducibility of $\sigma$ then implies that $W = Z$.  On the other
hand, $\lambda$ is equal to $0$ on $W$, and we are assuming that
$\lambda$ is a nonzero linear functional on $Z$.  Thus we conclude
that $v = 0$.  This shows that the kernel of $v \mapsto f_v$ is the
zero subspace of $Z$, and hence that this mapping is one-to-one.

	Let $L_x$, $x \in G$, denote the left regular representation
of $G$, as in Subsection \ref{Representations of finite groups}.
For each $x$ in $G$ and $v$ in $Z$ we have that
\begin{eqnarray}
	L_x(f_v)(y) = f_v(x^{-1}y) 
		& = & \lambda(\sigma_{y^{-1}x}(v))			\\
		& = & \lambda(\sigma_{y^{-1}}(\sigma_x(v)))
			= f_{\sigma_x(v)}(y).			\nonumber
\end{eqnarray}
This says exactly that the mapping $v \mapsto f_v$ intertwines the
representations $\sigma$ on $Z$ and the left regular representation on
$F(G)$.  In particular, $U$ is invariant under the left regular
representation.  Lemma \ref{irred. rep.'s isomorphic to a subrep. of
the reg. rep.} follows easily from these observations.

\beginremark
\label{lemma works for any field k, no assumption on char(k)}
{\rm Lemma \ref{irred. rep.'s isomorphic to a subrep. of the
reg. rep.} and its proof work for any field $k$, without any
assumption on the characteristic of $k$.}
\end{remark}

\beginproposition
\label{coeff.s not all 0 lead to nonzero operators for some irr. rep.}
Suppose that for each $x$ in $G$ an element $a_x$ of the field $k$ is
chosen, where $a_x \ne 0$ for at least one $x$.  Then there is a
vector space $Y$ over $k$ and an irreducible representation $\tau$ on
$Y$ such that the linear operator
\begin{equation}
\label{sum_{x in G} a_x tau_x}
	\sum_{x \in G} a_x \, \tau_x
\end{equation}
on $Y$ is not the zero operator.  
\end{proposition}

	Suppose first that we do not worry about having $\tau$ be an
irreducible representation.  Then we can use the left regular
representation on the space of $k$-valued functions on $G$.  The
operator in question is
\begin{equation}
\label{sum_{x in G} a_x L_x}
	\sum_{x \in G} a_x \, L_x.
\end{equation}
Let $\phi_e(z)$ denote the function on $G$ defined by $\phi_e(z) = 1$
when $z = e$ and $\phi_e(z) = 0$ when $z \ne e$.  Then
\begin{equation}
	\sum_{x \in G} a_x \, L_x(\phi_e)(w) = a_w
\end{equation}
for all $w$ in $G$.  The assumption that the $a_x$'s are not all
$0$ says exactly that this is not the zero function on $G$, and
hence that (\ref{sum_{x in G} a_x L_x}) is not the zero operator
on the space of functions on $G$.

	From Lemma \ref{decomp into irreducible pieces} in Subsection
\ref{Reducibility, continued} we know that the vector space of
$k$-valued functions on $G$ is spanned by subspaces which are
invariant under the left regular representation, and to which the
restriction of the left regular representation is irreducible.  It
follows easily that the restriction of the operator (\ref{sum_{x in G}
a_x L_x}) to at least one of these subspaces is nonzero, since it is
not the zero operator on the whole space of functions.  To get
Proposition \ref{coeff.s not all 0 lead to nonzero operators for some
irr. rep.}, we take $\tau$ to be the restriction of the left regular
representation to such a subspace.

	Now let $V$ be a vector space over $k$, and let $\rho$ be a
representation of $G$ on $V$ of the following special type (which is
unique up to isomorphism).  We assume that $V$ has an independent
system of subspaces $W_1, \ldots, W_m$ such that $V = \Span (W_1,
\ldots, W_m)$, each $W_j$ is invariant under $\rho$, the restriction
of $\rho$ to each $W_j$ is an irreducible representation of $G$, the
restriction of $\rho$ to $W_j$ is not isomorphic to the restriction of
$\rho$ to $W_l$ when $j \ne l$, and every irreducible representation
of $G$ on a vector space over $k$ is isomorphic to the restriction of
$\rho$ to some $W_j$, $1 \le j \le m$.  In other words, this
representation is isomorphic to a direct sum of irreducible
representations, in which every irreducible representation of $G$ on a
vector space over $k$ appears exactly once, up to isomorphism.  One
can also think in terms of starting with the left regular
representation of $G$ (over $k$) and passing to a suitable invariant
subspace, to avoid repetitions of irreducible representations (up to
isomorphism).

	Let $\mathcal{A}$ denote the set of linear operators on $V$
associated to the representation $\rho$, as in Subsection \ref{Algebras
from group representations}.

\beginlemma
\label{mathcal{A} corresponds to ``free'' linear sums of the rho_x's}
A linear operator $T \in \mathcal{L}(V)$ lies in $\mathcal{A}$ if and
only if it can be written as
\begin{equation}
\label{T = sum_{x in G} a_x rho_x}
	T = \sum_{x \in G} a_x \, \rho_x,
\end{equation}
where each $a_x$ lies in $k$.  Each operator $T \in \mathcal{A}$
can be written as (\ref{T = sum_{x in G} a_x rho_x}) in a unique
way.  In particular, the dimension of $\mathcal{A}$ as a vector space
is equal to the number of elements of $G$.
\end{lemma}

	The first statement is a rephrasal of the definition of
$\mathcal{A}$ given at the beginning of Subsection \ref{Algebras from
group representations}.  The uniqueness of the expression for $T$ in
(\ref{T = sum_{x in G} a_x rho_x}) follows from Proposition
\ref{coeff.s not all 0 lead to nonzero operators for some irr. rep.},
since all irreducible representations of $G$ occur as restrictions of
$\rho$ to a subspace of $V$.  The last assertion in the lemma is an
immediate consequence of the first two.

	Because the restrictions of $\rho$ to the $W_j$'s are
isomorphically distinct, Theorem \ref{mathcal{A} from a group rep., no
irred. comp's repeated} leads to rather precise information about how
$\mathcal{A}$ looks in this situation.

	If the field $k$ is algebraically closed, such as the field
${\bf C}$ of complex numbers, then the algebra of operators generated
by an irreducible representation is the algebra of all linear
operators on the corresponding vector space.  For this we use
Proposition \ref{k alg closed, irreducibility, imply mathcal{A}' = scalars I}
to say the commutant of the algebra consists of only the scalar
multiples of the identity, and then Theorem \ref{mathcal{A}'' =
mathcal{A} if mathcal{A} is very nice} in Subsection \ref{Very nice
algebras of operators}.  It follows that the number of elements of $G$
is equal to the sum of the squares of the degrees of the
isomorphically-distinct irreducible representations of $G$ in this
case.  Note that there is always an irreducible representation of
degree $1$, namely the representation of $G$ on a vector space of
dimension $1$ in which each $x$ in $G$ is associated to the identity
operator.

	Let us look now at the \emph{center} of $\mathcal{A}$, which
is the set of operators in $\mathcal{A}$ which commute with all other
elements of $\mathcal{A}$.

\beginlemma
An operator $T = \sum_{x \in G} a_x \, \rho_x$ in $\mathcal{A}$ lies
in the center of $\mathcal{A}$ if and only if $a_x = a_y$ whenever $x$
and $y$ are \emph{conjugate} inside the group $G$, i.e., whenever
there is a $w$ in $G$ such that $y = w x w^{-1}$.
\end{lemma}

	Indeed, $T$ commutes with all elements of $\mathcal{A}$
if and only if $T$ commutes with $\rho_z$ for all $z$ in the group.
This is the same as saying that
\begin{equation}
	\rho_z^{-1} \Bigl(\sum_{x \in G} a_x \, \rho_x \Bigr) \rho_z =
		\sum_{x \in G} a_x \, x
\end{equation}
for all $z$ in $G$.  This can be rewritten as
\begin{equation}
	\sum_{x \in G} a_x \, \rho_{z^{-1} x z}
			 = \sum_{x \in G} a_x \, \rho_x,
\end{equation}
or as
\begin{equation}
	\sum_{x \in G} a_{z x z^{-1}} \, \rho_x 
			= \sum_{x \in G} a_x \, \rho_x.
\end{equation}
In other words, the mapping from $x$ to $z x z^{-1}$ permutes the
elements of the group, and in the last step we made a change of
variables in the sum using this permutation.  We conclude that $T$
lies in the center of $\mathcal{A}$ if and only if $a_{z x z^{-1}} =
a_x$ for all $x$ and $z$ in the group.  This proves the lemma.

	If we write $x \sim y$ when $x$ and $y$ are conjugate elements
of the group, then this defines an equivalence relation on the group,
as is well-known and easy to verify.  The equivalence classes for this
relation are called \emph{conjugacy classes}.  The lemma can be
rephrased as saying that $T = \sum_{x \in G} a_x \, \rho_x$ lies in
the center of $\mathcal{A}$ if and only if the coefficients $a_x$ are
constant on the conjugacy classes of $G$.

	Note that the center of $\mathcal{A}$ is automatically a
subalgebra, and a vector subspace in particular.

\begincorollary
The dimension of the center of $\mathcal{A}$ is equal to the
number of conjugacy classes in the group.
\end{corollary}

	This is an easy consequence of the previous remarks.

	If the field $k$ is algebraically closed, then it follows that
the number of conjugacy classes in the group is equal to the number of
isomorphically-distinct irreducible representations of the group.
This uses Proposition \ref{k alg closed, irreducibility, imply
mathcal{A}' = scalars I}.

\section{$p$-adic numbers}
\label{$p$-adic numbers}
\index{p-adic numbers@$p$-adic numbers}
\setcounter{equation}{0}

	Let ${\bf Q}$\index{$Q$@${\bf Q}$} denote the field of
rational numbers, and fix a prime number $p$.

	The \emph{$p$-adic absolute value}\index{p-adic absolute
value@$p$-adic absolute value} $|\cdot|_p$\index{$"|\cdot"|_p$} on ${\bf
Q}$ is defined as follows.  If $x = \frac{a}{b} p^k$, where $a$, $b$,
and $k$ are integers, with $a, b \ne 0$ and neither $a$ nor $b$
divisible by $p$, then we set
\begin{equation}
	|x|_p = p^{-k}.
\end{equation}
If $x = 0$, then we set $|x|_p = 0$.

	The $p$-adic absolute value $|\cdot|_p$ enjoys many of
the same features as the classical absolute value $|\cdot|$,
defined by $|x| = x$ when $x \ge 0$ and $|x| = -x$ when $x < 0$.
In particular,
\begin{equation}
\label{|x + y|_p le |x|_p + |y|_p}
	|x + y|_p \le |x|_p + |y|_p
\end{equation}
and 
\begin{equation}
\label{|x y|_p = |x|_p |y|_p}
	|x y|_p = |x|_p \, |y|_p
\end{equation}
for all $x, y \in {\bf Q}$.  In fact, instead of (\ref{|x + y|_p le
|x|_p + |y|_p}) there is a stronger ``ultrametric''\index{ultrametric}
version of the triangle inequality, which states that
\begin{equation}
\label{|x + y|_p le max(|x|_p, |y|_p)}
	|x + y|_p \le \max(|x|_p, |y|_p),
\end{equation}
whose analogue for the standard absolute value $|\cdot|$ is not true.
It is not difficult to verify (\ref{|x y|_p = |x|_p |y|_p}) and
(\ref{|x + y|_p le max(|x|_p, |y|_p)}), just from the definitions.

	Just as for the standard absolute value, we can define a
$p$-adic distance function\index{p-adic distance function@$p$-adic
distance function} $d_p(x,y)$\index{$d_p(\cdot,\cdot)$} on ${\bf Q}$
using $|\cdot|_p$, by
\begin{equation}
	d_p(x,y) = |x - y|_p.
\end{equation}
This satisfies the usual requirements for a metric (as in \cite{Rudin}),
which is to say that $d_p(x,y)$ is a nonnegative real number which is
equal to $0$ if and only if $x = y$, $d_p(x,y) = d_p(y,x)$, and 
$d_p(x,y)$ satisfies the triangle inequality.  In this case we have
the stronger condition
\begin{equation}
\label{d_p(x,z) le max(d_p(x,y), d_p(y,z))}
	d_p(x,z) \le \max(d_p(x,y), d_p(y,z))
\end{equation}
for all $x, y, z \in {\bf Q}$, because of (\ref{|x + y|_p le
max(|x|_p, |y|_p)}).  Thus one says that $d_p(x,y)$ is an
\emph{ultrametric},\index{ultrametric} referring to this stronger
version of the triangle inequality.

	Recall that if $(M, D(u,v))$ is any metric space, then a
sequence $\{u\}_{j=1}^\infty$ of points in $M$ is said to
\emph{converge} to a point $u$ in $M$ if for every positive real
number $\epsilon$ there is a positive integer $L$ such that
\begin{equation}
	D(u_j, u) < \epsilon  \qquad\hbox{for all } j \ge L.
\end{equation}
Similarly, a sequence $\{v_j\}_{j=1}^\infty$ in $M$ is said to be
a \emph{Cauchy sequence} if for every $\epsilon > 0$ there is
an $L > 0$ so that
\begin{equation}
	D(v_j, v_k) < \epsilon \qquad\hbox{for all } j, k \ge L.
\end{equation}
It is easy to see that every convergent sequence is also a Cauchy
sequence.  It is not true in general that every Cauchy sequence
converges to some point in the metric space, but this is true
in some situations.  

	A metric space $(M, D(u,v))$ is said to be \emph{complete} if
every Cauchy sequence in $M$ also converges in $M$.  A classical
example is the real line ${\bf R}$, equipped with the standard metric
$|x - y|$, as in \cite{Rudin}.  The set ${\bf Q}$ of rational numbers
equipped with the standard metric $|x - y|$ is not complete, because a
sequence in ${\bf Q}$ which converges as a sequence in ${\bf R}$ is
always a Cauchy sequence, whether or not the limit lies in ${\bf Q}$.

	By a \emph{completion} of a metric space $(M, D(u,v))$ we mean
a metric space $(M_1, D_1(u,v))$ together with an embedding $\theta_1 :
M \to M_1$ with the following three properties: (a) the embedding is
isometric, in the sense that
\begin{equation}
	D_1(\theta_1(u), \theta_1(v)) = D(u,v)
			\qquad\hbox{for all } u, v \in M;
\end{equation}
(b) the image of $M$ in $M_1$ under $\theta_1$ is a dense subset of
$M_1$, so that
\begin{eqnarray}
	&& \hbox{for every $w \in M_1$ and $\epsilon > 0$ there is}	\\
	&& \hbox{a $u \in M$ such that $D_1(w, \theta_1(u)) < \epsilon$;}
								\nonumber
\end{eqnarray}
and (c) the metric space $(M_1, D_1(u,v))$ is itself complete.  A
completion of a metric space always exists, as in Problem 24 on p82
and Problem 24 on p170 of \cite{Rudin}.  If $(M_1, D_1(u,v))$,
$\theta_1 : M \to M_1$ and $(M_2, D_2(z,w))$, $\theta_2 : M \to M_2$
are two completions of the same metric space $(M, D(x,y))$, then there
is an isometry $\phi$ from $M_1$ onto $M_2$ such that $\theta_2 = \phi
\circ \theta_1$.  Indeed, one can define $\phi$ on the dense subset
$\theta_1(M)$ of $M_1$ through this equation, and then extend $\phi$
to all of $M_1$ using uniform continuity and completeness (as in
Problem 13 on p99 of \cite{Rudin}).	

	The standard embedding of ${\bf Q}$ in ${\bf R}$ provides
a completion of ${\bf Q}$ equipped with the standard metric $|x - y|$.
What about ${\bf Q}$ equipped with the $p$-adic metric $d_p(x,y)$?
Let us begin with a representation for rational numbers connected
to the $p$-adic geometry.

\beginlemma
\label{x a nonzero rational number, expansion in p's}
Let $x$ be a nonzero rational number, and let $l$ be the nonzero
integer such that $|x|_p = p^{-l}$.  There is a sequence
$\{\alpha_j\}_{j=l}^\infty$ of nonnegative integers such that
$\alpha_l \ne 0$, $\alpha_j \le p-1$ for all $j$, and the series
\begin{equation}
	\sum_{j=l}^\infty \alpha_j \, p^j
\end{equation}
converges to $x$ in ${\bf Q}$ with respect to the $p$-adic metric.
\end{lemma}

	In other words, if
\begin{equation}
\label{s_n = sum_{j=l}^n alpha_j p^j}
	s_n = \sum_{j=l}^n \alpha_j \, p^j
\end{equation}
for $n \ge l$, then $d_p(s_n, x) \to 0$ as $n \to \infty$.  In fact,
\begin{equation}
\label{d_p(s_n, x) le p^{-n-1} for each n ge l}
	d_p(s_n, x) \le p^{-n-1} \qquad\hbox{for each } n \ge l.
\end{equation}

	We can use (\ref{d_p(s_n, x) le p^{-n-1} for each n ge l}) to
choose the $\alpha_j$'s one after the other.  Let us start with
$\alpha_l$.  Since $|x|_p = p^{-l}$, we can write $x$ as $(a/b) \,
p^l$, where $a$ and $b$ are nonzero integers which are not divisible
by $p$.  We choose $\alpha_l$, $1 \le \alpha_l \le p-1$, so that $a
\equiv \alpha_l \, b$ modulo $p$.  This ensures that $d_p(s_l, x) \le
p^{-l-1}$, where $s_l = \alpha_l \, p^l$.

	Now suppose that $\alpha_j$ has been chosen for $j = l, l+1,
\ldots, n$ for some integer $n$, in such a way that (\ref{d_p(s_n, x)
le p^{-n-1} for each n ge l}) holds, and let us choose $\alpha_{n+1}$.
Since $d_p(s_n, x) \le p^{-n-1}$, we can write $x - s_n$ as $(a'/b')
\, p^{n+1}$, where $a'$ and $b'$ are integers, and $b'$ is not zero
and not divisible by $p$.  We choose $\alpha_{n+1}$, $0 \le
\alpha_{n+1} \le p-1$, so that $a' \equiv \alpha_{n+1} \, b'$ modulo
$p$.  This ensures that $d_p(s_{n+1}, x) \le p^{-n-2}$, with $s_{n+1}
= s_n + \alpha_{n+1} \, p^{n+1}$, as in (\ref{s_n = sum_{j=l}^n
alpha_j p^j}).

	This gives a sequence $\{\alpha_j\}_{j=l}^\infty$ as required
in Lemma \ref{x a nonzero rational number, expansion in p's}.  We shall
explain another way to get such a sequence in a moment, using the
next lemma.

\beginlemma
\label{series with integers to ones with integers in [0, p-1]}
Let $l$ be an integer, and let $\{\beta_j\}_{j=l}^\infty$ be a 
sequence of arbitrary nonnegative integers.  There is a
sequence $\{\alpha_j\}_{j=l}^\infty$ of nonnegative integers
such that $\alpha_j \le p-1$ for all $j$ and
\begin{equation}
	\sum_{j=l}^n \beta_j \, p^j - \sum_{j=l}^n \alpha_j \, p^j
\end{equation}
can be written as $c_n \, p^{n+1}$ with $c_n$ a nonnegative integer
for all integers $n \ge l$.
\end{lemma}

	This can be derived in much the same manner as before.  We
choose $\alpha_l$ so that it is equal to $\beta_l$ modulo $p$.
Suppose that $\alpha_j$, $l \le j \le n$, have the properties
described in the lemma, and let us choose $\alpha_{n+1}$.  By
assumption,
\begin{equation}
	c_n = p^{-n-1} \,
  \Bigl(\sum_{j=l}^n \beta_j \, p^j - \sum_{j=l}^n \alpha_j \, p^j \Bigr)
\end{equation}
is a nonnegative integer, and we choose $\alpha_{n+1}$, $0 \le
\alpha_{n+1} \le p-1$, so that it is equal to $c_n + \beta_{n+1}$
modulo $p$.  This implies that
\begin{equation}
	c_{n+1} = p^{-n-2} \,
\Bigl(\sum_{j=l}^{n+1} \beta_j \, p^j - \sum_{j=l}^{n+1} \alpha_j \, p^j \Bigr)
\end{equation}
is also a nonnegative integer, as desired.

	In the context of Lemma \ref{series with integers to ones with
integers in [0, p-1]}, if $\sum_{j=l}^\infty \beta_j \, p^j$ converges
to a rational number in the $p$-adic metric, then $\sum_{j=l}^\infty
\alpha_j \, p^j$ converges to the same rational number in the
$p$-adic metric, since the difference between the partial sums tends
to $0$ in the $p$-adic absolute value.

	Let us return to the setting of Lemma \ref{x a nonzero
rational number, expansion in p's}.  Let $x$ be a nonzero rational
number, and assume first that $x$ is a positive integer.
Then we can write $x$ as
\begin{equation}
	x = \sum_{j=0}^n \alpha_j \, p^j,
\end{equation}
where $n$ is a nonnegative integer and each $\alpha_j$ is a
nonnegative integer such that $0 \le \alpha_j \le p - 1$.  One
can begin by choosing $n$ as large as possible so that $p^n \le x$,
and then take $\alpha_n$ so that $0 \le x - \alpha_n \, p^n < p^n$.
Afterwards the remaining part $x - \alpha_n \, p_n$ can be treated in
a similar manner, and so on.  Alternatively, $\alpha_0$ is uniquely
determined by the requirement that $x \equiv \alpha_0$ modulo $p$,
$\alpha_1$ is determined by $x - \alpha_0 \equiv \alpha_1 \, p$ modulo
$p^2$, etc.  

	More generally, if $x$ is a positive integer times $p^l$ for
some $l \in {\bf Z}$, then there is a finite expansion for $x$ as in
Lemma \ref{x a nonzero rational number, expansion in p's}.  Of course
we can write $x$ in this way where the integer factor is not divisible
by $p$, so that the leading term in the expansion is nonzero.

	What about negative numbers?  Consider $x = -1$, for instance.
We can write $-1$ as $(p-1)/(1-p)$, and this leads to the expansion
\begin{equation}
	-1 = \sum_{j=0}^\infty (p-1) \, p^j.
\end{equation}
This is a kind of geometric series, which would not converge in the
classical sense, because $p > 1$.  However, this series does converge
to $-1$ with respect to the $p$-adic metric.  Indeed, by the usual
formula, 
\begin{equation}
	\sum_{j=0}^n (p-1) \, p^j = (p-1) \sum_{j=0}^n p^j
		= (p-1) \frac{1 - p^{n+1}}{1-p} = - 1 - p^{n+1},
\end{equation}
so that
\begin{equation}
	d_p\Bigl(\sum_{j=0}^n (p-1) \, p^j, -1 \Bigr)
		= |(-1 - p^{n+1}) - 1|_p = |-p^{n+1}|_p = p^{-n-1}.
\end{equation}

	More generally, suppose that $x = -1/b$, where $b$ is a positive
integer which is not divisible by $p$.  Let $c$ be a positive integer
such that $c \le p-1$ and $-b \, c \equiv 1$ modulo $p$.  Thus 
$x = -1/b = c/(1 - b_1 \, p)$, where $b_1$ is a positive integer.
As in the preceding situation, we get that
\begin{equation}
	x = \sum_{j=0}^\infty c \, b_1^j \, p^j,
\end{equation}
with convergence in the $p$-adic metric, because
\begin{equation}
	\sum_{j=0}^n c \, b_1^j \, p^j = c \sum_{j=0}^n (b_1 \, p)^j
		= c \frac{1 - (b_1 \, p)^{n+1}}{1 - b_1 \, p}
		= x \, (1 - (b_1 \, p)^{n+1}).
\end{equation}

	This expansion does not quite fit the conditions in Lemma
\ref{x a nonzero rational number, expansion in p's}, because the
coefficients $c \, b_1^j$, which are positive integers, are not
bounded by $p-1$ in general.  However, one can use Lemma \ref{series
with integers to ones with integers in [0, p-1]} to convert this
expansion into one that does satisfy the additional restriction on the
coefficients.

	Before proceeding, let us note the following.

\beginlemma
\label{products of expansions}
Let $l$, $l'$ be integers, and let $\{\beta_j\}_{j=l}^\infty$ and
$\{\beta'_j\}_{j=l'}^\infty$ be sequences of nonnegative integers.
Define a sequence $\{\gamma_i\}_{i=l+l'}^\infty$ by
\begin{equation}
	\gamma_i = \sum \{\beta_j \, \beta'_k : 
				i = j + k, j \ge l, k \ge l'\}
\end{equation}
(i.e., $\{\gamma_i\}_{i=l+l'}^\infty$ is the ``Cauchy product'' of
$\{\beta_j\}_{j=l}^\infty$ and $\{\beta'_j\}_{j=l'}^\infty$).
Then, for each integer $n \ge 0$,
\begin{equation}
\label{product of sums minus sum of Cauchy product}
	\Bigl(\sum_{j=l}^{l + n} \beta_j \, p^j \Bigr) \,
		\Bigl(\sum_{k=l'}^{l' + n} \beta'_k \, p^k \Bigr)
	- \sum_{i = l + l'}^{l + l' + n} \gamma_i \, p^i
\end{equation}
is a nonnegative integer which is divisible by $p^{l+l'+n+1}$.
\end{lemma}

	To see this, one might prefer to think about the case where $l
= l' = 0$, to which one can easily reduce.  The main point is that 
every term occurring in the sum of the $\gamma_i$'s in (\ref{product
of sums minus sum of Cauchy product}) also occurs when one expands the
product of the two sums in (\ref{product of sums minus sum of Cauchy
product}).  There are additional terms in the product of the sums, 
which are all of the form $c \, p^m$, where $c$ is a nonnegative integer
and $m \ge l+l'+n+1$.

\begincorollary
\label{products of sums}
In the situation of Lemma \ref{products of expansions}, suppose that
the series $\sum_{j=l}^\infty \beta_j \, p^j$ and $\sum_{k=l'}^\infty
\beta'_k \, p^k$ converge to rational numbers $u$, $u'$, respectively,
in the $p$-adic metric.  Then the series $\sum_{i = l + l'}^\infty
\gamma_i \, p^i$ converges to the product $u \, u'$ in the $p$-adic
metric.
\end{corollary}

	This is an easy consequence of the lemma.  Note that this
result does not hold for infinite series in the classical sense
without some additional hypotheses.  See \cite{Rudin} for more
information.

	To finish our alternate discussion of Lemma \ref{x a nonzero
rational number, expansion in p's}, suppose that $x$ is a nonzero
rational number.  For special cases we have given expansions for $x$
in terms of sums of products of nonnegative integers with powers of
$p$.  One can get arbitrary $x$ by taking products, using Corollary
\ref{products of sums} to ensure that this leads to the correct
limits.  One also employs Lemma \ref{series with integers to ones with
integers in [0, p-1]} to obtain an expansion in which the coefficients
are bounded by $p-1$.

	Let us record some more lemmas.

\beginlemma
\label{bound for p-adic absolute value of sums}
Suppose that $l$ is an integer, $\{\beta_j\}_{j=l}^\infty$ is a
sequence of integers (of arbitrary sign), and that the sum
$\sum_{j=l}^\infty \beta_j \, p^j$ converges in ${\bf Q}$ in the
$p$-adic metric.  Then
\begin{equation}
	\biggl|\sum_{j=l}^\infty \beta_j \, p^j \biggr|_p \le p^{-l}.
\end{equation}
\end{lemma}

	This is not hard to verify, from the definitions.  Of
course
\begin{equation}
	\biggl|\sum_{j=m}^n \beta_j \, p^j \biggr|_p \le p^{-m}
\end{equation}
for all integers $m$ and $n$ such that $l \le m \le n$.

\beginlemma
\label{p-adic absolute value in terms of l under suitable conditions}
Suppose that $l$ is an integer, $\{\beta_j\}_{j=l}^\infty$ is a
sequence of integers (of arbitrary sign), and that the sum
$\sum_{j=l}^\infty \beta_j \, p^j$ converges in ${\bf Q}$ in the
$p$-adic metric.  If $\beta_l \not\equiv 0$ modulo $p$, then
\begin{equation}
	\biggl|\sum_{j=l}^\infty \beta_j \, p^j \biggr|_p = p^{-l}.
\end{equation}
\end{lemma}

	Again, this is not hard to verify.  By assumption, the leading
term in the sum has $p$-adic absolute value $p^{-l}$, while the rest
has $p$-adic absolute value $\le p^{-l-1}$.

\beginlemma
\label{uniqueness of the expansion for rational numbers}
The expansion described in Lemma \ref{x a nonzero rational number,
expansion in p's} is unique.  In other words, suppose that $x$ is a
nonzero rational number, $l$, $l'$ are integers, and
$\{\alpha_j\}_{j=l}^\infty$, $\{\alpha'_j\}_{j=l'}^\infty$ are
sequences of nonnegative integers such that $\alpha_l \ne 0$,
$\alpha'_{l'} \ne 0$, $0 \le \alpha_j \le p - 1$ for each $j \ge l$,
$0 \le \alpha'_j \le p - 1$ for each $j \ge l'$, and the sums
$\sum_{j=l}^\infty \alpha_j \, p^j$, $\sum_{j=l'}^\infty \alpha'_j \,
p^j$ both converge in ${\bf Q}$ to $x$ in the $p$-adic metric.  Then
$l = l'$ and $\alpha_j = \alpha'_j$ for all $j \ge l$.
\end{lemma}

	Notice first that $l = l'$ under these conditions, since
$|x|_p = p^{-l}$ and $|x|_p = p^{-l'}$, as in Lemma \ref{p-adic
absolute value in terms of l under suitable conditions}.  Next,
$\alpha_l = \alpha'_l$, because otherwise the difference between the
two series has $p$-adic absolute value $p^{-l}$.  One can subtract the
leading terms and repeat the argument to get that $\alpha_j =
\alpha'_j$ for all $j$, as desired.

\beginlemma
\label{uniqueness for 0}
Suppose that $l$ is an integer and that $\{\alpha_j\}_{j=l}^\infty$
is a sequence of nonnegative integers such that $\alpha_j \le p-1$
for all $j$.  If $\sum_{j=l}^\infty \alpha_j \, p^j$ converges
in ${\bf Q}$ to $0$ in the $p$-adic metric, then $\alpha_j = 0$
for all $j$.
\end{lemma}

	This is easy to check.

	Let us now define ${\bf Q}_p$,\index{$Q_p$@${\bf Q}_p$} the
set of $p$-adic numbers, to consist of the formal sums of the form
\begin{equation}
\label{sum_{j = l}^infty alpha_j p^j}
	\sum_{j = l}^\infty \alpha_j \, p^j,
\end{equation}
where $l$ is an integer, and each $\alpha_j$ is an integer such that
$0 \le \alpha_j \le p - 1$ for each $j \ge l$.  If
\begin{equation}
	\sum_{j = l'}^\infty \alpha'_j \, p^j
\end{equation}
is another such formal sum, with $l' < l$, say, then the two sums are
viewed as representing the same element of ${\bf Q}_p$ if $\alpha'_j =
0$ when $l' \le j < l$ and $\alpha'_j = \alpha_j$ when $j \ge l$.  

	If a formal sum (\ref{sum_{j = l}^infty alpha_j p^j}) actually
converges to a rational number $x$ with respect to the $p$-adic
metric, then let us identify that element of ${\bf Q}_p$ with $x \in
{\bf Q}$.  The preceding lemmas imply that no more than one element of
${\bf Q}_p$ is identified in this manned with a single rational number
$x$, and Lemma \ref{x a nonzero rational number, expansion in p's}
says that every rational number is included in ${\bf Q}_p$ through
this recipe.  Thus we can think of ${\bf Q}$ as a subset of ${\bf
Q}_p$.

	The $p$-adic absolute value\index{p-adic absolute
value@$p$-adic absolute value} $|\cdot|_p$\index{$"|\cdot"|_p$} can be
extended to ${\bf Q}_p$ by setting
\begin{equation}
	\biggl|\sum_{j = l}^\infty \alpha_j \, p^j \biggr|_p
		= p^{-l}
\end{equation}
when $\alpha_l \ne 0$ and $0 \le \alpha_j \le p-1$, $\alpha_j \in {\bf
Z}$ for all $j$ (as usual).  The zero element of ${\bf Q}_p$, which
corresponds to the series with all coefficients equal to $0$, and to
the rational number $0$ under the identification discussed in the
previous paragraph, has $p$-adic absolute value equal to the number
$0$.  For the elements of ${\bf Q}_p$ which correspond to rational
numbers, this definition of the $p$-adic absolute value agrees with
the original one, as in the earlier lemmas.

	Next, let us extend the $p$-adic distance
function\index{p-adic distance function@$p$-adic distance function}
$d_p(\cdot,\cdot)$\index{$d_p(\cdot,\cdot)$} to ${\bf Q}_p$.
Suppose that
\begin{equation}
	\sum_{j = l}^\infty \alpha_j \, p^j, \quad
		\sum_{j = l'}^\infty \alpha'_j \, p^j
\end{equation}
are two elements of ${\bf Q}_p$, where, as usual, $l$ and $l'$ are
integers, $\alpha_j$ is an integer such that $0 \le \alpha_j \le p-1$
for all $j \ge l$, and $\alpha'_j$ is an integer such that $0 \le
\alpha'_j \le p-1$ for all $j \ge l'$.  If either of these two sums is
the zero element of ${\bf Q}_p$, so that all of the coefficients are
$0$, then the $p$-adic distance between the two sums is defined to be
the $p$-adic absolute value of the other sum.  If they are both the
zero element of ${\bf Q}_p$, then the $p$-adic distance is equal to
$0$.  Assume now that both sums correspond to nonzero elements of
${\bf Q}_p$.  In this event we ask that $\alpha_l \ne 0$ and
$\alpha'_{l'} \ne 0$, which can always be arranged by dropping initial
terms with coefficient $0$, if necessary.  If $l \ne l'$, then the
$p$-adic distance between the two sums is defined to be $p^{- \min(l,
l')}$.  Suppose instead that $l = l'$.  If $\alpha_j = \alpha'_j$ for
all $j \ge l$, then the two sums are the same, and the $p$-adic
distance between them is defined to be $0$.  Otherwise, let $j_0$ be
the smallest integer $\ge l$ such that $\alpha_{j_0} \ne
\alpha'_{j_0}$.  In this case the distance between the two sums is
defined to be $p^{-j_0}$.  It is not hard to check that this definition
of the $p$-adic distance agrees with the earlier one when both sums
correspond to rational numbers.

	By definition, the $p$-adic distance on ${\bf Q}_p$ is
nonnegative and equal to $0$ exactly when the two elements of ${\bf
Q}_p$ are the same.  The distance is also symmetric in the two
elements of ${\bf Q}_p$.  As before, the $p$-adic distance
satisfies the ultrametric\index{ultrametric} condition
\begin{equation}
	d_p(w_1, w_3) \le \max(d_p(w_1, w_2), d_p(w_2, w_3))
\end{equation}
for all $w_1$, $w_2$, and $w_3$ in ${\bf Q}_p$.  This is not
difficult to verify.

	In this way ${\bf Q}_p$ becomes a metric space.  It is easy to
see that ${\bf Q}$ defines a dense subset of ${\bf Q}_p$, since sums
(\ref{sum_{j = l}^infty alpha_j p^j}) in which all but finitely many
coefficients are $0$ correspond to rational numbers, and arbitrary
sums can be approximated by these.  One can also show that ${\bf Q}_p$
is complete as a metric space, so that it is indeed a completion of
${\bf Q}$.

	Like the real line with its standard metric, ${\bf Q}_p$
enjoys the property that closed and bounded subsets of it are 
compact.  For this it is enough to know that closed balls around
the origin are compact.  Fix an integer $l$, and consider the
set of $w$ in ${\bf Q}_p$ such that $|w|_p \le p^{-l}$.  This
can be described exactly as the set of sums
\begin{equation}
	\sum_{j=l}^\infty \alpha_j \, p^j,
\end{equation}
where each $\alpha_j$ is an integer such that $0 \le \alpha_j \le
p-1$.  Now all of the sums have the same starting point (at $l$), but
it is important to allow the initial coefficients to be $0$, or even
all of the coefficients to be $0$, to get the right subset of ${\bf
Q}_p$.

	Topologically, this set is equivalent to a Cantor set (with
$p$ pieces at each stage).  It is convenient to look at the topology
in terms of sequences in the set, and convergence of sequences.
Namely, a sequence $\{w_t\}_{t=1}^\infty$ of elements of this set
converges to another element $w$ of this set exactly if, for each $j
\ge l$, the $j$th coefficient in the sum associated to the $w_t$'s is
equal to the $j$th coefficient of the sum associated to $w$ for all
sufficiently large $t$.  (How large $t$ should be is permitted to
depend on $j$.)  The statement that the set is compact means that for
any sequence $\{w_t\}_{t=1}^\infty$ of elements of the set there is a
subsequence that converges to an element of the set.  This can be
established through classical arguments.

	Suppose that $l$ is an integer, and that $\{\beta_j\}_{j=l}^\infty$
is a sequence of nonnegative integers, without the restriction
$\beta_j \le p-1$.  Consider the formal sum
\begin{equation}
\label{sum_{j=l}^infty beta_j p^j}
	\sum_{j=l}^\infty \beta_j \, p^j.
\end{equation}
We can view this as still giving rise to an element of ${\bf Q}_p$,
using Lemma \ref{series with integers to ones with integers in [0, p-1]}
to convert this to a sum $\sum_{j=l}^\infty \alpha_j \, p^j$ where
each $\alpha_j$ is a nonnegative integer such that $\alpha_j \le p-1$.

	This remark permits us to define operations of addition and
multiplication on ${\bf Q}_p$.  Specifically, one first adds or
multiplies two sums in the obvious manner, to get a sum in the form
(\ref{sum_{j=l}^infty beta_j p^j}).  One then gets an element of
${\bf Q}_p$ as described in the previous paragraph.

	For elements of ${\bf Q}_p$ that correspond to rational
numbers, these operations of addition and multiplication give the same
result as the usual ones.  The elements of ${\bf Q}_p$ corresponding
to the rational numbers $0$ and $1$ are additive and multiplicative
identity elements for all of ${\bf Q}_p$.  One also has the usual
commutative, associative, and distributive laws on ${\bf Q}_p$.
Furthermore, these operations define continuous mappings from ${\bf
Q}_p \times {\bf Q}_p$ into ${\bf Q}_p$, with respect to the $p$-adic
metric.

	As before, $-1$ can be written as $\sum_{j=0}^\infty (p-1) \,
p^j$.  Multiplication by $-1$ gives additive inverses in ${\bf Q}_p$,
just as in ${\bf Q}$.  If $w$ is a nonzero element of ${\bf Q}_p$,
then $w$ has a multiplicative inverse in ${\bf Q}_p$.  To see this, it
is convenient to find a representation for $-1/w$, from which one can
get $1/w$ by multiplication by $-1$.  Assume that $w$ is given as
$\sum_{j=l}^\infty \alpha_j \, p^j$, where $l$ is an integer,
$\alpha_j$ is a nonnegative integer such that $\alpha_j \le p-1$ for
all $j$, and $\alpha_l \ne 0$.  Let $c$ be an integer such that $1 \le
c \le p-1$ and $c \, \alpha_l \equiv -1$ modulo $p$, and let $\beta$ be
the nonnegative integer such that $c \, \alpha_l + 1 = \beta \, p$.
Consider the expression
\begin{equation}
\label{big sum for -1/w}
	\sum_{j=0}^\infty c \, p^{-l} \, 
   \Bigl(\beta \, p + \sum_{i=1}^\infty c \, \alpha_{l+i} \, p^i \Bigr)^j.
\end{equation}
It is easy to expand this out to get a sum of the form
(\ref{sum_{j=l}^infty beta_j p^j}), and hence an element of ${\bf
Q}_p$.  One can think of (\ref{big sum for -1/w}) as corresponding to
\begin{equation}
	\frac{-1}{w} = \frac{c \, p^{-l}}{1 - (c \, w \, p^{-l} - 1)}
		= \frac{c \, p^{-l}}{1 - \Bigl(\beta \, p + 
		\sum_{i=1}^\infty c \, \alpha_{l+i} \, p^i \Bigr)}.
\end{equation}
This works, since the product of
\begin{equation}
	\sum_{j=0}^\infty 
   \Bigl(\beta \, p + \sum_{i=1}^\infty c \, \alpha_{l+i} \, p^i \Bigr)^j
\end{equation}
and 
\begin{equation}
	(\beta \, p - 1) + \sum_{i=1}^\infty c \, \alpha_{l+i} \, p^i 
	= \sum_{i=0}^\infty c \, \alpha_{l+i} \, p^i
\end{equation}
is equal to $1$.

	Thus ${\bf Q}_p$ is a field.  With subtraction defined, it
makes sense to say that $d_p(w,z) = |w-z|_p$ on ${\bf Q}_p$.  It can
be somewhat more convenient to write this as $d_p(w,w+u) = |u|_p$ for
$w$, $u$ in ${\bf Q}_p$.  Notice too that
\begin{equation}
	|w + z|_p \le \max(|w|_p, |z|_p)
\end{equation}
and
\begin{equation}
	|w z|_p = |w|_p \, |z|_p
\end{equation}
for $w$, $z$ in ${\bf Q}_p$.

\section{Absolute values on fields}
\label{Absolute values on fields}
\index{absolute values (on a field)}
\setcounter{equation}{0}

	Let $k$ be a field.  A function $|\cdot |_*$ on $k$ is called
an \emph{absolute value function} on $k$, or a choice of
\emph{absolute values}, if $|x|_*$ is a nonnegative real number for
all $x$ in $k$, $|x|_* = 0$ if and only if $x = 0$, and
\begin{equation}
\label{|x y|_* = |x|_* |y|_*}
	|x \, y|_* = |x|_* \, |y|_*
\end{equation}
and
\begin{equation}
\label{|x + y|_* le |x|_* + |y|_*}
	|x + y|_* \le |x|_* + |y|_*
\end{equation}
for all $x$, $y$ in $k$.  

	Notice that (\ref{|x y|_* = |x|_* |y|_*}) yields
\begin{equation}
	|1|_* = 1,
\end{equation}
where the $1$ on the left side is the multiplicative identity
element in $k$, and the $1$ on the right side is the real number $1$.
If $x$ is a nonzero element of $k$, then we obtain that
\begin{equation}
	|x^{-1}|_* = |x|_*^{-1},
\end{equation}
because $1 = |1|_* = |x \, x^{-1}|_* = |x|_* \, |x^{-1}|_*$.  Also,
\begin{equation}
	|-1|_* = 1,
\end{equation}
since $|-1|_*^2 = |(-1)^2|_* = |1|_* = 1$.  As a result,
$|-x|_* = |x|_*$ for all $x$ in $k$.

	For example, the usual absolute values for real numbers
defines an absolute value function on ${\bf R}$ and ${\bf Q}$.  The
usual modulus of a complex number defines an absolute value function
on ${\bf C}$.

	An absolute value function $|\cdot |_*$ on $k$ is said to be
\emph{non-Archimedian},\index{non-Archimedian} or
\emph{ultrametric},\index{ultrametric} if it satisfies the stronger
condition
\begin{equation}
\label{|x + y|_* le max(|x|_*, |y|_*)}
	|x + y|_* \le \max(|x|_*, |y|_*)
\end{equation}
for all $x$, $y$ in $k$, instead of (\ref{|x + y|_* le |x|_* + |y|_*}).
The $p$-adic absolute values on ${\bf Q}$ or on ${\bf Q}_p$ have this
property.

	On any field one can define the trivial\index{trivial absolute
value function} absolute value function, which is equal to $0$ at the
zero element of the field and to $1$ at all nonzero elements in the
field.  Note that this is an ultrametric absolute value function.
If $k$ is a finite field, then it is not hard to see that the only
absolute value function on $k$ is the trivial one.

	If $k_0$ is a field, let $k_0(t)$ denote the field of rational
functions in one variable $t$ over $k_0$.  Thus every element of
$k_0(t)$ can be written as $P(t)/Q(t)$, where $P(t)$ and $Q(t)$ are
polynomials in $t$ with coefficients in $k_0$, and $Q(t)$ has at least
one nonzero coefficient.  Two such representations
\begin{equation}
	P_1(t)/Q_1(t), \ P_2(t)/Q_2(t)
\end{equation}
are considered to define the same element of $k_0(t)$ when 
\begin{equation}
	P_1(t) \, Q_2(t) = P_2(t) \, Q_1(t).
\end{equation}
If $P(t)$ is the zero polynomial, so that all of its coefficients are
$0$, then $P(t)/Q(t)$ is the zero rational function in $k_0(t)$ for
any $Q(t)$ which is not the zero polynomial.

	Fix a real number $A$ with $A > 1$.  We can define a
nontrivial absolute value function on $k_0(t)$ by setting it equal to
$0$ for the zero rational function, and to $A^{-l}$ when the rational
function can be expressed as $t^l \, P_0(t)/Q_0(t)$, where $P_0(t)$
and $Q_0(t)$ are polynomials in $t$ with nonzero constant terms.  It
is easy to see that this satisfies the conditions described before,
including the ultrametric version of the triangle inequality.  If we
restrict this absolute value function to the subfield of $k_0(t)$
consisting of constant functions, then we get the trivial absolute
value function on $k_0$.

	Let $k_0((t))$ denote the field of formal Laurent expansions
of finite order over $k_0$ in $t$.  In other words, an element of $k_0((t))$
is given by a formal series
\begin{equation}
\label{sum_{j=n}^infty a_j t^j}
	\sum_{j=n}^\infty a_j \, t^j,
\end{equation}
where $n$ is an integer and the $a_j$'s are arbitrary elements of
$k_0$.  If $m$ is an integer with $m \le n$ and
\begin{equation}
	\sum_{l=m}^\infty b_j \, t^j
\end{equation}
is another such series, then the two series are viewed as defining the
same element of $k_0((t))$ if and only if $b_j = 0$ when $m \le j < n$
and $a_j = b_j$ when $j \ge n$.  The series with all coefficients
equal to $0$ and whatever starting point $t^n$ correspond to the zero
element of $k_0((t))$.  Elements of $k_0((t))$ can be added and
multiplied in the usual manner, term by term.  The series with
constant term equal to $1$ and all others equal to $0$ is the
multiplicative identity element in $k_0((t))$.  If $R(t)$ is a nonzero
element of $k_0((t))$, then $R(t)$ can be written as
\begin{equation}
	R(t) = a_n \, t^n (1 - R_0(t)),
\end{equation}
where $n$ is an integer, $a_n$ is a nonzero element of $k_0$, and
$R_0(t) \in k_0((t))$ is a given by a sum with only positive powers of
$t$.  One can check that the standard formula
\begin{equation}
	(1 - R_0(t))^{-1} = \sum_{j=0}^\infty R_0(t)^j,
\end{equation}
with $R_0(t)^0$ interpreted as being $1$, makes sense and works in
$k_0((t))$.  Thus nonzero elements of $k_0((t))$ have multiplicative
inverses, so that $k_0((t))$ is a field.  There is a natural embedding
of $k_0(t)$ in $k_0((t))$, by associating to a rational function over
$k_0$ its Laurent expansion around $0$.  Of course polynomials are
contained in $k_0((t))$ as finite sums of multiples of nonnegative
powers of $t$, and rational functions can be obtained from this using
products and multiplicative inverses.

	Let $A$ be a real number with $A > 1$, as before.  One can
define an absolute value function on $k_0((t))$ by taking the absolute
value of the zero element of $k_0((t))$ to be $0$, and the absolute
value of (\ref{sum_{j=n}^infty a_j t^j}) to be $A^{-n}$ when $a_n \ne
0$.  It is not difficult to check that this defines an absolute value
function with the ultrametric version of the triangle inequality, and
that it agrees with the one described earlier for $k_0(t)$ when
applied to Laurent series coming from rational functions.

	Suppose that $k$ is a field and that $|\cdot |_*$ is an
absolute value function on $k$ which satisfies the ultrametric version
of the triangle inequality.  We shall say that $|\cdot |_*$ is
\emph{nice}\index{nice absolute value function (on a field)} if there
is a subset $E$ of the set of nonnegative real numbers such that
$|x|_* \in E$ for all $x$ in $k$ and $E$ has no limit point in the
real line except possibly at $0$.  This is equivalent to saying that
for any real numbers $a$, $b$ such that $0 < a < b$, the set $E \cap
[a,b]$ is finite.

\beginlemma
\label{characterization of | |_* being nice}
An ultrametric absolute value function $|\cdot |_*$ on a field $k$
is nice if and only if there is a real number $r$ such that $0 \le r <
1$ and $|x|_* \le r$ for all $x \in k$ such that $|x|_* < 1$.
\end{lemma}

	The ``only if'' part of this statement is immediate from the
definition.  Conversely, suppose that there is a real number $r$ as in
the lemma.  If $x$, $y$ are elements of $k$ such that $|x|_* < |y|_*$,
then $|x \, y^{-1}|_* < 1$, so that $|x \, y^{-1}|_* \le r$ and $|x|_*
\le r \, |y|_*$.  One can use this to show that the absolute value
function takes values in a set $E$ as above.

\beginlemma
\label{``local total boundedness'' implies nice}
An ultrametric absolute value function $|\cdot |_*$ on a field $k$ is
nice if there is a positive real number $s < 1$ and a finite
collection $x_1, \ldots, x_m$ of elements of $k$ such that for each
$y$ in $k$ with $|y|_* < 1$ there is an $x_j$, $1 \le j \le m$, such
that $|y - x_j|_* \le s$.
\end{lemma}

	To prove this, suppose that the hypothesis of Lemma
\ref{``local total boundedness'' implies nice} holds, and let us check
that the condition in Lemma \ref{characterization of | |_* being nice}
is satisfied.  If $1 \le j \le m$ and $y \in k$ satisfy $|y|_* < 1$ and
$|y - x_j|_* \le s$, then
\begin{equation}
	|x_j|_* \le \max(|x_j - y|_*, |y|_*) < 1.
\end{equation}
We may as well assume that $|x_j|_* < 1$ for all $j$, since otherwise
we can reduce to a smaller collection of $x_j$'s with the same
property as in the lemma.  Take $r$ to be the maximum of $s$ and the
numbers $|x_j|_*$, $1 \le j \le m$, so that $r < 1$.  If $y$ is an
element of $k$ and $|y|_* < 1$, then $|y - x_j|_* \le s$ for some $j$,
and hence $|y|_* \le \max(|y-x_j|_*, |x_j|_*) \le r$, as desired.

\section{Norms on vector spaces}
\label{Norms on vector spaces}
\index{norms}
\setcounter{equation}{0}

	Fix a field $k$ and an absolute value function $|\cdot |_*$ on
$k$, and let $V$ be a vector space over $k$.  A \emph{norm} on $V$
with respect to this choice of absolute value function on $k$ is a
real-valued function $N(\cdot)$ on $V$ such that the following three
properties are satisfied: (a) $N(v) \ge 0$ for all $v$ in $V$, with
$N(v) = 0$ if and only if $v = 0$; (b) $N(\alpha \, v) = |\alpha|_* \,
N(v)$ for all $\alpha$ in $k$ and $v$ in $V$; (c) $N(v + w) \le N(v) +
N(w)$ for all $v$, $w$ in $V$.

	In this section we make the standing assumption that $|\cdot
|_*$ is an ultrametric absolute value function on $k$, and we shall
restrict our attention to norms $N$ on vector spaces $V$ over $k$ with
respect to $|\cdot |_*$ that are \emph{ultrametric
norms},\index{ultrametric} in the sense that
\begin{equation}
	N(v+w) \le \max(N(v), N(w))
\end{equation}
for all $v$, $w$ in $V$.  Observe that if $N(\cdot)$ is an ultrametric
norm on $V$, then $d(v,w) = N(v-w)$ is an ultrametric on $V$, so that
\begin{equation}
	d(u,w) \le \max(d(u,v), d(v,w))
\end{equation}
for all $u$, $v$, and $w$ in $V$.

	One can think of $k$ as a $1$-dimensional vector space over
itself, and then the absolute value function $|\cdot |_*$ defines an
ultrametric norm on this vector space.  If $n$ is a positive integer,
then $k^n$, the space of $n$-tuples of elements of $k$, is an
$n$-dimensional vector space over $k$, with respect to coordinatewise
addition and scalar multiplication.  Consider the expression
\begin{equation}
	\max_{1 \le j \le n} |x_j|_*
\end{equation}
for each $x = (x_1, \ldots, x_n)$ in $k^n$.  It is easy to check that
this defines an ultrametric norm on $k^n$.

\beginlemma
\label{N(v+w) = max(N(v), N(w)) when N(v) ne N(w)}
Let $V$ be a vector space over $k$, and let $N$ be an ultrametric norm
on $V$.  Suppose that $v$, $w$ are elements of $V$, and that $N(v) \ne
N(w)$.  Then $N(v + w) = \max(N(v), N(w))$.
\end{lemma}

	We have that $N(v + w) \le \max(N(v), N(w))$ for all $v$, $w$
in $V$, and so we want to show that the reverse inequality holds when
$N(v) \ne N(w)$.  Assume for the sake of definiteness that
$N(v) < N(w)$.  Notice that
\begin{eqnarray}
	N(w) = N((v+w) + (-v)) & \le & \max(N(v+w), N(-v)) 		\\
		& = & \max(N(v+w), N(v)).			\nonumber
\end{eqnarray}
Because $N(v) < N(w)$, this implies that $N(w) \le N(v+w)$.  Hence
\begin{equation}
	\max(N(v), N(w)) \le N(v+w),
\end{equation}
as desired.

\begincorollary
\label{N(sum v_j) = max(N(v_j)) when the N(v_j)'s are distinct}
Let $V$ be a vector space over $k$, and let $N$ be a ultrametric norm
on $V$.  Suppose that $m$ is a positive integer, and that $v_1,
\ldots, v_m$ are elements of $V$ such that $N(v_j) = N(v_l)$ only when
either $j = l$ or $v_j = v_l = 0$.  Then
\begin{equation}
	N \Bigl(\sum_{j=1}^m v_j \Bigr) = \max_{1 \le j \le m} N(v_j).
\end{equation}
\end{corollary}

	This can be derived from Lemma \ref{N(v+w) = max(N(v), N(w)) 
when N(v) ne N(w)} using induction.

\beginlemma
\label{values of an ultrametric norm N on V}
Suppose that $V$ is a vector space over $k$ of dimension $n$, and that
$N$ is an ultrametric norm on $V$.  Let $E$ be a subset of the set of
nonnegative real numbers such that $|x|_* \in E$ for all $x$ in $k$.
If $v_1, \ldots, v_{n+1}$ are nonzero elements of $V$, then at
least one of the ratios $N(v_j)/N(v_l)$, $1 \le j < l \le n+1$,
lies in $E$.
\end{lemma}

	Let $v_1, \ldots, v_{n+1}$ be nonzero vectors in $V$.  If the
ratios $N(v_j)/N(v_l)$ do not lie in $E$ for any $j$, $l$ with $j \ne
l$, then one can check that $v_1, \ldots, v_{n+1}$ are linearly
independent in $V$, using Corollary \ref{N(sum v_j) = max(N(v_j)) when
the N(v_j)'s are distinct}.  This contradicts the assumption that
$V$ has dimension $n$, and the lemma follows.

	In analogy with the definition of a nice ultrametric absolute
value function in Section \ref{Absolute values on fields}, let us say
that an ultrametric norm $N$ on a vector space $V$ over $k$ is
\emph{nice}\index{nice ultrametric norm} if there is a subset $E_1$ of
the set of nonnegative real numbers such that $N(v)$ lies in $E_1$ for
every vector $v$ in $V$ and $E_1$ has no limit point in the real line
except possibly for $0$.  The latter is equivalent to asking that $E_1
\cap [a,b]$ be finite for every pair $a$, $b$ of positive real
numbers.

\begincorollary
\label{N nice if abs. val. fcn. nice}
Let $V$ be a vector space over $k$, and let $N$ be an ultrametric norm
on $V$.  If the absolute value function $|\cdot |_*$ on $k$ is
nice, then $N$ is a nice ultrametric norm on $V$.
\end{corollary}

	More precisely, if $V$ has dimension $n$, and if $E$ is a
subset of the set of nonnegative real numbers such that the absolute
value function $|\cdot |_*$ on $k$ takes values in $E$, then there are
positive real numbers $a_1, \ldots, a_n$ such that $N$ takes values in
the set $E_1 = \bigcup_{j=1}^n a_j \, E$.  Here $a \, E = \{a \, s : s
\in E\}$.  This can be derived from Lemma \ref{values of an
ultrametric norm N on V}.  (Note that the $a_j$'s need not be distinct.)

	Fix a positive integer $n$, and let us take our vector space
to be $k^n$.  We shall say that a norm $N$ on $k^n$ is
\emph{nondegenerate} if there is a positive real number $c$ such that
\begin{equation}
\label{nondegeneracy condition}
	c \, \max_{1 \le j \le n} |x_j|_* \le N(x)
\end{equation}
for all $x = (x_1, \ldots, x_n)$ in $k^n$.  This condition is
automatically satisfied if $\{y \in k : |y|_* \le 1\}$ is a compact
subset of $k$, using the metric $|u-v|_*$ on $k$.  To see this,
observe that it is enough to check that $N$ has a positive lower bound
on the set of $x$ in $k^n$ such that $\max_{1 \le j \le n} |x_j|_* =
1$, because of homogeneity.  By assumption, this is a compact subset
of $k^n$ in the product topology, and $N$ is positive at each element
of this set since $N$ is a norm.  One can verify that $N$ is continuous,
and in fact $N$ is locally constant away from $0$ in $k^n$.

\beginremark
\label{bound for a norm}
{\rm For any norm $N$ on $k^n$ there is a positive real number $C$
so that
\begin{equation}
\label{general upper bound for N(x)}
	N(x) \le C \, \max_{1 \le j \le n} |x_j|_*
\end{equation}
for all $x$ in $k^n$.  This is easy to see, by writing $x$ as a linear
combination of standard basis vectors, and it applies to norms in general,
whether or not they are ultrametric norms.  For that matter, the notion
of nondegeneracy makes sense for norms in general.
}
\end{remark}

\beginlemma
\label{point in W subseteq k^n as close as possible to z}
Suppose that the absolute value function $|\cdot |_*$ on $k$ is nice.
Let $N$ be a nondegenerate ultrametric norm on $k^n$, let $W$ be a
vector subspace of $k^n$, and let $z$ be an element of $k^n$ which
does not lie in $W$.  There exists an element $x_0$ of $W$ such that
the distance $N(z-x_0)$ is as small as possible.
\end{lemma}

	Notice first that there is a $\delta > 0$ such that $N(z-x)
\ge \delta$ for all $x \in W$.  This uses the fact that $W$ is a
closed subset of $k^n$ with respect to the product topology.  For
instance, one can describe $W$ as the set of vectors where finitely
many linear functions vanish, and these linear functions are
continuous.  Because $|\cdot |_*$ is nice, $N$ is nice, as in
Corollary \ref{N nice if abs. val. fcn. nice}, and hence the infimum
of $N(z-x)$ over $x$ in $W$ is attained, since it is positive.

\beginremark
\label{points in W of minimal distance to z}
{\rm If $w \in W$ satisfies $N(w) \le N(z-x_0)$, then $N(z-(x_0+w))
\le N(z-x_0)$, so that $N(z-(x_0+w)) = N(z-x_0)$ if $N(z-x_0)$
is as small as possible.  Conversely, if $N(z-(x_0+w)) = N(z-x_0)$,
then $N(w) \le N(z-x_0)$.  }
\end{remark}

\beginlemma
\label{lemma about extending linear maps, keeping norm fixed}
Assume that $|\cdot |_*$ is a nice absolute value function on $k$.
Let $N$ be a nondegenerate ultrametric norm on $k^n$, let $V_1$ be a
vector space over $k$, and let $N_1$ be an ultrametric norm on $V_1$.
Suppose that $W$ is a vector subspace of $k^n$, and that $T$ is a
linear mapping from $W$ to $V_1$.  Assume also that $A$ is a
nonnegative real number such that
\begin{equation}
	N_1(T(v)) \le A \, N(v)
\end{equation}
for all $v$ in $W$.  Then there is a linear mapping $T_1$ from $k^n$
to $V_1$ such that $T_1(v) = T(v)$ when $v$ lies in $W$, and
\begin{equation}
\label{N_1(T_1(v)) le A N(v)}
	N_1(T_1(v)) \le A \, N(v)
\end{equation}
for all $v$ in $k^n$.
\end{lemma}

	To prove this, let $V_1$, $W$, $T$, etc., be given as above.
We may as well assume that $W$ is not all of $k^n$, since otherwise
there is nothing to do.  Let $z$ be an element of $k^n$ which does not
lie in $W$, and set $W_1 = \Span (W,z)$.  We would like to extend $T$
first to $W_1$.

	Let $x_0$ be an element of $W$ such that $N(z-x_0)$ is as
small as possible, as in Lemma \ref{point in W subseteq k^n as close as
possible to z}, and set $y_0 = T(x_0)$.  Define $T_1$ on $W_1$ by
\begin{equation}
	T_1(v + \alpha \, z) = T(v) + \alpha \, y_0
\end{equation}
for all $v$ in $W$ and $\alpha$ in $k$.  With this definition
$T_1$ is clearly linear on $W_1$ and agrees with $T$ on $W$.  We want
to check that
\begin{equation}
	N_1(T_1(v + \alpha \, z)) \le A \, N(v + \alpha \, z)
\end{equation}
for all $v$ in $W$ and $\alpha$ in $k$.  This inequality holds
when $\alpha = 0$ by assumption, and so we may restrict ourselves to
$\alpha \ne 0$.  By the homogeneity of the norms, we are reduced to
showing that
\begin{equation}
	N_1(T_1(v + z)) \le A \, N(v + z)
\end{equation}
for all $v$ in $W$, which is the same as
\begin{equation}
	N_1(T(v) + y_0) \le A \, N(v + z).
\end{equation}
We can rewrite this further as
\begin{equation}
	N_1(T(v + x_0)) \le A \, N(v + z),
\end{equation}
by the definition of $y_0$.  

	Since
\begin{equation}
	N_1(T(v + x_0)) \le A \, N(v + x_0)
\end{equation}
by hypothesis (because $x_0$ lies in $W$), it suffices to show
that
\begin{equation}
\label{N(v + x_0) le N(v + z)}
	N(v + x_0) \le N(v + z)
\end{equation}
for all $v$ in $W$.  Notice that
\begin{eqnarray}
	N(v + x_0) & = & N((v+z) + (-z+x_0)) 			\\
		& \le & \max(N(v+z), N(-z+x_0))		\nonumber \\
		& =   & \max(N(v+z), N(z-x_0)).		\nonumber
\end{eqnarray}
By the manner in which $x_0$ was chosen, $N(v+z) \ge N(z-x_0)$.
This yields (\ref{N(v + x_0) le N(v + z)}).

	Thus $T$ can be extended to a linear mapping $T_1$ from the
larger subspace $W_1$ into $V_1$ with the inequality (\ref{N_1(T_1(v))
le A N(v)}) holding there.  By repeating the process, we can get an
extension to all of $k^n$.  This proves Lemma \ref{lemma about
extending linear maps, keeping norm fixed}.

\begincorollary
\label{projection onto a subspace W}
Assume that $|\cdot |_*$ is a nice absolute value function on $k$.
Let $N$ be a nondegenerate ultrametric norm on $k^n$, and let $W$ be a
vector subspace of $k^n$.  There is a linear mapping $P : k^n \to W$
which is a projection, so that $P(w) = w$ when $w \in W$ and $P(v)$
lies in $W$ for all $v$ in $k^n$, and which satisfies
\begin{equation}
	N(P(v)) \le N(v)
\end{equation}
for all $v$ in $k^n$.
\end{corollary}

	This follows by defining $P$ first on $W$ to be the identity
mapping, and then extending this to a mapping from $k^n$ to $W$ using
Lemma \ref{lemma about extending linear maps, keeping norm fixed}.

	Suppose that $W$ and $P$ are as in the corollary, and let
$W_0$ denote the kernel of $P$.  We may as well assume that $W$ is
neither the zero subspace of $k^n$ nor all of $k^n$, so that the same is
true of $W_0$.  As usual, if $I$ denotes the identity mapping on $k^n$,
then $I-P$ is a projection of $k^n$ onto $W_0$, i.e., $(I-P)(v)$ lies in
$W_0$ for all $v$ in $k^n$ and $(I-P)(v) = v$ when $v$ lies in $W_0$.
Here we also have that
\begin{equation}
	N((I-P)(v)) \le N(v)
\end{equation}
for all $v$ in $k^n$, because
\begin{eqnarray}
	N((I-P)(v)) & = & N(v - P(v)) 				\\
		    & \le & \max(N(v), N(P(v))) = N(v).		\nonumber
\end{eqnarray}

	We can go further and say that
\begin{equation}
	N(v) = \max(N(P(v)), N((I-P)(v)))
\end{equation}
for all $v$ in $k^n$.  That is,
\begin{eqnarray}
	N(v) & = & N(P(v) + (I-P)(v)) 				\\
	     & \le & \max(N(P(v)), N((I-P)(v))),		\nonumber
\end{eqnarray}
while the reverse inequality follows from the ones that have already
been derived.

	Let us pause a moment for some terminology.  Suppose that
$V$ is a vector space over $k$ of dimension $n$, and that
$x_1, \ldots, x_n$ is a basis for $V$.  This determines a dual
family $f_1, \ldots, f_n$ of linear functionals on $V$, i.e.,
each $f_j$ is a linear mapping from $V$ into $k$ such
that $f_j(x_j) = 1$ and $f_j(x_l) = 0$ when $j \ne l$.

	Now assume that $V$ is equipped with an ultrametric norm $N$.
We say that the basis $x_1, \ldots, x_n$ is
\emph{normalized}\index{normalized basis} if
\begin{equation}
\label{normalization condition}
	N(x_j) \cdot |f_j(v)|_* \le N(v)
\end{equation}
for all vectors $v$ in $V$.  Note that equality occurs when $v = x_j$.

\beginlemma
\label{normalized bases exist}
Assume that $|\cdot |_*$ is a nice absolute value function on $k$.
If $N$ is a nondegenerate ultrametric norm on $k^n$, then there is
a normalized basis for $k^n$ with respect to $N$.
\end{lemma}

	When $n = 1$ there is nothing to do.  For $n > 1$ one can
use induction, with the projection operators discussed above
permitting one to go from $n$ to $n+1$.

\beginlemma
\label{formula for ultrametric norm N in terms of a normalized basis}
Let $V$ be a vector space over $k$ with dimension $n$, equipped with
an ultrametric norm $N$, and let $x_1, \ldots, x_n$ be a normalized
basis for $V$, with dual linear functionals $f_1, \ldots, f_n$.  Then
\begin{equation}
\label{N(v) = max {N(x_j) |f_j(v)|_* : 1 le j le n }}
	N(v) = \max \{N(x_j) \cdot |f_j(v)|_* : 1 \le j \le n \}
\end{equation}
for all vectors $v$ in $V$.
\end{lemma}

	The normalization condition (\ref{normalization condition})
gives one inequality in (\ref{N(v) = max {N(x_j) |f_j(v)|_* : 1 le j
le n }}).  For the opposite inequality, we write
\begin{equation}
	v = f_1(v) \, x_1 + \cdots + f_n(v) \, x_n,
\end{equation}
so that 
\begin{equation}
	N(v) \le \max \{N(f_j(v) \, x_j) : 1 \le j \le n \}.
\end{equation}

\beginremark
\label{further normalization}
{\rm If $x_1, \ldots, x_n$ is a normalized basis for $V$, then we can
multiply the $x_j$'s by arbitrary nonzero elements of $k$ and get a
new normalized basis for $V$.  In particular, if $N$ has the feature
that it takes values in the set of nonnegative real numbers which
occur as values of $|\cdot |_*$, then we can multiply the $x_j$'s by
elements of $k$ to get a basis $y_1, \ldots, y_n$ which is 
normalized and satisfies $N(y_j) = 1$ for each $j$.  If $h_1, \ldots,
h_n$ is the corresponding dual basis of linear functionals on $V$,
then (\ref{N(v) = max {N(x_j) |f_j(v)|_* : 1 le j le n }}) becomes
\begin{equation}
\label{N(v) = max {|h_j(v)|_* : 1 le j le n }}
	N(v) = \max \{|h_j(v)|_* : 1 \le j \le n \}
\end{equation}
for all $v$ in $V$.
}
\end{remark}

	Let $G$ be a finite group, and $\rho$ be a representation of
$G$ on $k^n$.  Let $N_0$ be any fixed ultrametric norm on $k^n$.
Define $N(v)$ on $V$ by
\begin{equation}
\label{N(v) = max {N_0(rho_a(v)) : a in G}}
	N(v) = \max \{N_0(\rho_a(v)) : a \in G\}.
\end{equation}
It is not difficult to check that this defines an ultrametric norm on
$k^n$.  Indeed, each $N(\rho_a(v))$ defines an ultrametric norm on $k^n$,
and passing to the maximum yields an ultrametric norm as well.  By
construction, this norm has the property that
\begin{equation}
	N(\rho_b(v)) = N(v)
\end{equation}
for all $b$ in the group $G$ and all vectors $v$ in $V$.  If we choose
$N_0$ so that it is nondegenerate, such as $N_0(x) = \max_{1 \le j \le
n} |x_j|_*$, then $N$ is also nondegenerate.  Similarly, if we choose
$N_0$ so that it takes values in the set of nonnegative real numbers
which occur as values of the absolute value function $|\cdot |_*$ on
$k$, as with $N_0(x) = \max_{1 \le j \le n} |x_j|_*$, then $N$ has the
same feature.

\section{Operator norms}
\label{Operator norms}
\setcounter{equation}{0}

	Let $k$ be a field with an absolute value function $|\cdot
|_*$.  Fix a positive integer $n$, and consider the vector space $k^n$
of $n$-tuples of elements of $k$.  Let $N(\cdot )$ be a norm on $k^n$
with respect to $|\cdot |_*$, which we assume to be nondegenerate, in
the sense of (\ref{nondegeneracy condition}).

	If $T$ is a linear operator from $k^n$ to itself, then the
\emph{operator norm} $\|T\|_{op}$ of $T$ with respect to $N$ can be
defined by
\begin{equation}
	\|T\|_{op} = \sup \{N(T(x)) : x \in k^n, \ N(x) = 1\}.
\end{equation}
This supremum is finite because of the nondegeneracy condition for $N$.
One can reformulate this definition as saying that
\begin{equation}
	N(T(x)) \le \|T\|_{op} \, N(x)
\end{equation}
for all $x$ in $k^n$, and that $\|T\|_{op}$ is the smallest nonnegative
real number with this property.  

	The collection of linear operators from $k^n$ to itself forms
a vector space over $k$ in the usual manner, and $\|T\|_{op}$ defines
a norm on this vector space, i.e., $\|T\|_{op}$ is a nonnegative
real number which is equal to $0$ exactly when $T$ is the zero linear
operator on $k^n$,
\begin{equation}
	\|\alpha \, T\|_{op} = |\alpha |_* \, \|T\|_{op}
\end{equation}
for all $\alpha$ in $k$ and linear operators $T$ on $k^n$, and
\begin{equation}
	\|T_1 + T_2\|_{op} \le \|T_1\|_{op} + \|T_2\|_{op}
\end{equation}
for all linear operators $T_1$, $T_2$ on $k^n$.  These properties
follow from the corresponding features of the norm $N$.  Notice that
\begin{equation}
	\|I\|_{op} = 1,
\end{equation}
where $I$ denotes the identity operator on $k^n$, and that
\begin{equation}
	\|T_1 \circ T_2\|_{op} \le \|T_1\|_{op} \, \|T_2\|_{op}
\end{equation}
for any two linear operators $T_1$, $T_2$.  If $|\cdot |_*$ is an
ultrametric absolute value function and $N(\cdot )$ is an ultrametric norm,
then the operator norm $\|\cdot \|_{op}$ is also an ultrametric norm,
so that
\begin{equation}
	\|T_1 + T_2\|_{op} \le \max(\|T_1\|_{op}, \|T_2\|_{op})
\end{equation}
for all linear operators $T_1$, $T_2$ on $k^n$.

	Suppose that $T$, $A$ are linear operators on $k^n$, $T$ is
invertible, and 
\begin{equation}
\label{||A||_{op} < ||T^{-1}||_{op}^{-1}}
	\|A\|_{op} < \|T^{-1}\|_{op}^{-1}.
\end{equation}
In this case $T + A$ is also an invertible linear operator on $k^n$.
To see this, it suffices to show that the kernel of $T + A$ is
trivial.  Let $x$ be any element of $k^n$.  If $y = T(x)$, then
\begin{equation}
	N(T^{-1}(y)) \le \|T^{-1}\|_{op} \, N(y),
\end{equation}
which is the same as saying that
\begin{equation}
	\|T^{-1}\|_{op}^{-1} \, N(x) \le N(T(x)).
\end{equation}
We also have that
\begin{eqnarray}
\label{N(T(x)) = N((T+A)(x) - A(x)), etc}
	N(T(x)) & = & N((T+A)(x) - A(x)) 			\\
		& \le & N((T+A)(x)) + N(A(x))		\nonumber \\
		& \le &  N((T+A)(x)) + \|A\|_{op} \, N(x),
							\nonumber
\end{eqnarray}
and hence
\begin{equation}
	(\|T^{-1}\|_{op}^{-1} - \|A\|_{op}) \, N(x) \le N((T + A)(x)).
\end{equation}
Thus $(T + A)(x) \ne 0$ when $x \ne 0$, so that $T + A$ is invertible,
and we get the norm estimate
\begin{equation}
	\|(T + A)^{-1}\|_{op} \le (\|T^{-1}\|_{op}^{-1} - \|A\|_{op})^{-1}.
\end{equation}

	If $|\cdot |_*$ is an ultrametric absolute value function and
$N$ is an ultrametric norm, then (\ref{N(T(x)) = N((T+A)(x) - A(x)),
etc}) can be replaced with
\begin{eqnarray}
	N(T(x)) & \le & \max(N((T+A)(x)), N(A(x))) 		\\
		& \le & \max(N((T+A)(x)), \|A\|_{op} \, N(x)),
								\nonumber
\end{eqnarray}
so that
\begin{equation}
	\|T^{-1}\|_{op}^{-1} \, N(x) \le \max(N((T+A)(x)), \|A\|_{op} \, N(x))
\end{equation}
for all $x$ in $k^n$.  Because $\|A\|_{op} < \|T^{-1}\|_{op}^{-1}$, we
obtain that
\begin{equation}
	\|T^{-1}\|_{op}^{-1} \, N(x) \le N((T+A)(x)).
\end{equation}
More precisely, one derives this initially for $x \ne 0$, and then notes
that it holds automatically for $x = 0$, so that the inequality applies
to all $x$ in $k^n$.  This leads to
\begin{equation}
	\|(T + A)^{-1}\|_{op} \le \|T^{-1}\|_{op}
\end{equation}
in the ultrametric case (under the condition (\ref{||A||_{op} <
||T^{-1}||_{op}^{-1}})).

\begindefinition
\label{def of uniform division algebra}
Let $\mathcal{B}$ be an algebra of operators on $k^n$ which is a
division algebra (Definition \ref{def of division algebra (of
operators)}).  We say that $\mathcal{B}$ is a \emph{uniform division
algebra}\index{uniform division algebra} if there is a positive
real number $C$ so that
\begin{equation}
\label{||T^{-1}||_{op} le C ||T||_{op}^{-1}}
	\|T^{-1}\|_{op} \le C \, \|T\|_{op}^{-1} 
\end{equation}
for all nonzero operators $T$ in $\mathcal{B}$.
\end{definition}

\beginremark
\label{||T^{-1}||_{op} ge ||T||_{op}^{-1} automatic}
{\rm The opposite inequality is automatic, i.e.,
\begin{equation}
	1 \le \|T^{-1}\|_{op} \, \|T\|_{op},
\end{equation}
since $T^{-1} \circ T = I$ has norm $1$. }
\end{remark}

	This property does not depend on the choice of a nondegenerate
norm $N$ on $k^n$, although the choice of $N$ can affect the constant
in the definition.  If $k$ is locally compact, so that $\{\alpha \in k
: |\alpha|_* \le 1\}$ is a compact set with respect to the metric
$|\alpha - \beta|_*$ on $k$, then any algebra of operators on $k^n$
which is a division algebra is in fact a uniform division algebra.
Indeed, it suffices to establish (\ref{||T^{-1}||_{op} le C
||T||_{op}^{-1}}) for the $T$'s such that $\|T\|_{op} = 1$, and this
is a compact subset of the vector space of linear operators on $k^n$
under our assumption.  It is not hard to use compactness to get a
uniform bound for $\|T^{-1}\|_{op}$, as desired.

	For the rest of this section, let us assume that $|\cdot |_*$
is an ultrametric absolute value function on $k$, and consider the
ultrametric norm
\begin{equation}
	N_1(x) = \max_{1 \le j \le n} |x_j|_*
\end{equation}
on $k^n$.

	Let $T$ be a linear operator on $k^n$.  Thus $T$ can be described
by an $n \times n$ matrix $(a_{j,l})$ with entries in $k$, so that
\begin{equation}
	(T(x))_j = \sum_{l=1}^n a_{j,l} \, x_l
\end{equation}
for all $x$ in $k^n$, where $(T(x))_j$ denotes the $j$th component of
$T(x)$, just as $x_l$ denotes the $l$th component of $x$.  The operator
norm $\|T\|_{op}$ of $T$ with respect to $N_1$ can be given by
\begin{equation}
\label{||T||_{op} = max_{1 le j, l le n} |a_{j,l}|_*}
	\|T\|_{op} = \max_{1 \le j, l \le n} |a_{j,l}|_*.
\end{equation}
To check that $\|T\|_{op}$ is less than or equal to the right side
of the equation, one can use the ultrametric property of $|\cdot |_*$.
For the opposite inequality, one can look at $T(x)$ in the special case
where $x$ has one component equal to $1$ and the rest equal to $0$.

	When is $T$ an isometry, i.e., when does $T$ satisfy
\begin{equation}
	N_1(T(x)) = N_1(x)
\end{equation}
for all $x$ in $k^n$?  This is clearly equivalent to asking that
$\|T\|_{op} \le 1$ and $\|T^{-1}\|_{op} \le 1$.  Let us check that
this happens if and only if $\|T\|_{op} \le 1$ and $|\det T|_* = 1$,
where $\det T$ denotes the determinant of $T$ (or, equivalently, of
the matrix $(a_{j,l})$).  Notice first that if $\|T\|_{op} \le 1$,
then $|\det T|_* \le 1$.  This uses merely the definition of the
determinant, together with the fact that the entries of the matrix
associated to $T$ have absolute value less than or equal to $1$, and
the properties of the absolute value function.  Similarly, the
determinant of any square submatrix of $T$ has absolute value in $k$
less than or equal to $1$.  By Cramer's rule, $T^{-1}$ is given by
$(\det T)^{-1}$ times the cofactor transpose of $T$.  The matrix
entries for the cofactor transpose have absolute value less than or
equal to $1$ when $\|T\|_{op} \le 1$, since they are given in terms of
determinants of submatrices of the matrix of $T$.  Thus the operator
norm of the cofactor transpose is less than or equal to $1$ when
$\|T\|_{op} \le 1$, because of the formula (\ref{||T||_{op} = max_{1
le j, l le n} |a_{j,l}|_*}) applied to the cofactor transpose (instead
of $T$).  If $|\det T|_* = 1$, then it follows that $\|T^{-1}\|_{op}
\le 1$.  Conversely, if $\|T^{-1}\|_{op} \le 1$, then $|\det T^{-1}|_*
\le 1$, and hence $|\det T|_* \ge 1$, since $\det T^{-1} = (\det
T)^{-1}$.  This implies that $|\det T|_* = 1$ if $\|T\|_{op} \le 1$
too.

	As in Section \ref{Norms on vector spaces}, one way that
isometries on $k^n$ come up is through representations of finite
groups on $k^n$.  That is, we might start with the norm $N_1$, not
necessarily invariant under the group representation, and then obtain
an invariant norm using (\ref{N(v) = max {N_0(rho_a(v)) : a in G}}).
The new norm can be converted back into $N_1$ after a linear change of
variables on $k^n$, if the absolute value function $|\cdot |_*$ is
nice, as in Section \ref{Norms on vector spaces}.  In other words, the
representation of the finite group can be conjugated by an invertible
linear transformation on $k^n$ to get a representation that acts by
isometries with respect to $N_1$.  Compare with Appendix 1 of Chapter
IV in Part II of \cite{Serre-lalg}.

	In addition to norms for operators on a single vector space,
one can consider norms for operators between vector spaces.  Suppose
that $V_1$, $V_2$ are vector spaces over the field $k$, and that
$N_1(\cdot)$, $N_2(\cdot)$ are norms on $V_1$, $V_2$ with respect to
the absolute value function $|\cdot |_*$ on $k$.  We assume that $N_1$
is nondegenerate, as before, with respect to a linear equivalence
between $V_1$ and $k^n$, where $n$ is the dimension of $V_1$.  It is
easy to check that the condition of nondegeneracy does not depend on
the choice of linear isomorphism between $V_1$ and $k^n$.

	If $T$ is a linear mapping from $V_1$ to $V_2$, then the
operator norm $\|T\|_{op, 12}$ of $T$ with respect to $N_1$, $N_2$ on
$V_1$, $V_2$ is defined by
\begin{equation}
	\|T\|_{op, 12} = \sup \{N_2(T(x)) : x \in V_1, N_1(x) = 1\}.
\end{equation}
The supremum makes sense because of the assumption of nondegeneracy for
$N_1$.  As before, the operator norm is characterized by the conditions
that
\begin{equation}
	N_2(T(x)) \le \|T\|_{op, 12} \, N_1(x)
\end{equation}
for all $x$ in $V_1$ and that $\|T\|_{op, 12}$ be the smallest
nonnegative real number for which this property holds.

	The operator norm defines a norm on the vector space of linear
mappings from $V_1$ to $V_2$, as one can easily verify.  If $|\cdot
|_*$ is an ultrametric absolute value function on $k$, and $N_1$,
$N_2$ are ultrametric norms on $V_1$, $V_2$, then $\|\cdot \|_{op,
12}$ is an ultrametric norm on the vector space of linear mappings
from $V_1$ to $V_2$.

	Suppose that $V_3$ is another vector space over $k$, and that
$N_3$ is a norm on $V_3$ with respect to the absolute value function
$|\cdot |_*$ on $k$.  Assume that the norm $N_2$ on $V_2$ is also
nondegenerate in the sense described before.  We can define operator
norms $\|\cdot \|_{op, 13}$ and $\|\cdot \|_{op, 23}$ for linear
mappings from $V_1$ to $V_3$ and from $V_2$ to $V_3$, respectively,
using the norms $N_1$, $N_2$, $N_3$ on $V_1$, $V_2$, $V_3$.  If $T$ is
a linear mapping from $V_1$ to $V_2$ and $S$ is a linear mapping from
$V_2$ to $V_3$, so that the composition $S \circ T$ defines a linear
mapping from $V_1$ to $V_3$, then 
\begin{equation}
	\|S \circ T \|_{op, 13} \le \|S\|_{op, 23} \, \|T\|_{op, 12}.
\end{equation}
This is easy to check.

\section{Vector spaces of linear mappings}
\label{Vector spaces of linear mappings}
\setcounter{equation}{0}

	Let $k$ be a field, and let $V$, $W$ be vector spaces over $k$.
Recall that $\mathcal{L}(V,W)$ denotes the vector space of linear
mappings from $V$ to $W$.  For notational simplicity, we shall write
$H$ for $\mathcal{L}(V,W)$ in this section.

	Suppose that $\mathcal{A}$ is an algebra of operators on $V$.
We can associate to $\mathcal{A}$ an algebra of operators
$\mathcal{A}_{H,1}$ on $H$, where a linear operator on $H$ lies in
$\mathcal{A}_{H,1}$ if it is of the form $R \mapsto R \circ T$, $R \in
H$, where $T$ lies in $\mathcal{A}$.  Similarly, if $\mathcal{B}$ is
an algebra of linear operators on $W$, then we can associate to it an
algebra $\mathcal{B}_{H,2}$ consisting of operators on $H$ of the form
$R \mapsto S \circ R$, $R \in H$, where $S$ lies in $\mathcal{B}$.
In other words, we are using composition of operators in $H$ with
operators on $V$ or on $W$ to define linear operators on $H$.

\beginremark
\label{mathcal{A} = mathcal{L}(V), mathcal{B} = mathcal{L}(W)}
{\rm As a basic case, let $\mathcal{A}$ be $\mathcal{L}(V)$ and let
$\mathcal{B}$ be all of $\mathcal{L}(W)$.  It is easy to check that
every element of $\mathcal{A}_{H,1}$ commutes with every element of
$\mathcal{B}_{H,2}$.  In fact,
\begin{equation}
	(\mathcal{A}_{H,1})' = \mathcal{B}_{H,2}, \quad
		(\mathcal{B}_{H,2})' = \mathcal{A}_{H,1},
\end{equation}
where the prime refers to the commutant of the algebra as a subalgebra
of $\mathcal{L}(H)$.  This is not difficult to show.
}
\end{remark}

\beginremark
\label{relation with ``expanding an algebra''}
{\rm Let $\mathcal{B}$ be any algebra of operators on $W$.  The
algebra of operators $\mathcal{B}_{H,2}$ on $H$ is essentially the
same as ``expanding'' $\mathcal{B}$ as in Subsection \ref{Expanding an
algebra}, with $n$ equal to the dimension of $V$.  In this way,
Proposition \ref{Expanding an algebra} can be reformulated as saying
that $(\mathcal{B}'')_{H,2} = (\mathcal{B}_{H,2})''$.  The description of
$(\mathcal{B}_{H,2})'$ in the proof of the proposition can be rephrased
as saying that $(\mathcal{B}_{H,2})'$ contains $(\mathcal{B}')_{H,2}$ and
$(\mathcal{L}(V))_{H,1}$, and is generated by them.
}
\end{remark}

\beginremark
\label{algebras from group representations}
{\rm Let $G_1$, $G_2$ be finite groups, and let $\sigma$, $\tau$ be
representations of $G_1$, $G_2$ on $V$, $W$, respectively.  We can
define representations $\widetilde{\sigma}$, $\widetilde{\tau}$ of
$G_1$, $G_2$ on $H$ by
\begin{equation}
	\widetilde{\sigma}_x(R) = R \circ (\sigma_x)^{-1}, 
		\quad \widetilde{\tau}_y(R) = \tau_y \circ R
\end{equation}
for $x \in G_1$, $y \in G_2$, and $R \in H$.  If $\mathcal{A}$ is the
algebra of operators on $V$ generated by $\sigma$, then
$\mathcal{A}_{H,1}$ is the same as the algebra of operators on $H$
generated by $\widetilde{\sigma}$.  Similarly, if $\mathcal{B}$ is the
algebra of operators on $W$ generated by $\tau$, then
$\mathcal{B}_{H,2}$ is the same as the algebra of operators on $H$
generated by $\widetilde{\tau}$.
}
\end{remark}

\beginremark
\label{k a symmetric field, inner products}
{\rm Suppose that $k$ is a symmetric field, and that $V$, $W$ are
equipped with inner products $\langle \cdot, \cdot \rangle_V$,
$\langle \cdot, \cdot \rangle_W$, as in Subsection \ref{Inner product
spaces}.  For every linear mapping $R : V \to W$ there is an
``adjoint''\index{adjoint of a linear operator}\index{$T^*$
(adjoint of an operator $T$)} $R^*$, a linear mapping from $W$
to $V$, characterized by the condition
\begin{equation}
	\langle R(v), w \rangle_W = \langle v, R^*(w) \rangle_V
\end{equation}
for all $v \in V$, $w \in W$.  Note that if $T$ is a linear operator
on $V$, and if $S$ is a linear operator on $W$, then $R \circ T$, 
$S \circ R$ are linear mappings from $V$ to $W$, and
\begin{equation}
	(R \circ T)^* = T^* \circ R^*, \quad (S \circ R)^* = R^* \circ S^*.
\end{equation}
Here $T^*$ is the adjoint of $T$ as an operator on $V$, $S^*$ is the
adjoint of $S$ as an operator on $W$, and $R^*$, $(R \circ T)^*$, and
$(S \circ R)^*$ are the adjoints of $R$, $R \circ T$, and $S \circ R$
as operators from $V$ to $W$.  (In general, if one has three inner
product spaces, a linear mapping from the first to the second, and a
linear mapping from the second to the third, then the adjoint of the
composition is equal to the composition of the adjoints, in the
opposite order.)

	We can define an inner product $\langle \cdot, \cdot
\rangle_H$ on $H$ by
\begin{equation}
\label{def of lange R_1, R_2 rangle_H}
	\langle R_1, R_2 \rangle_H = {\tr}_V R_2^* \circ R_1
				= {\tr}_W R_1 \circ R_2^*,
\end{equation}
where $\tr_V T$ denotes the trace of a linear operator $T$ on $V$ and
$\tr_W S$ denotes the trace of a linear operator $S$ on $W$.  Recall
that if $U_1$ is any linear mapping from $V$ to $W$, and $U_2$ is any
linear mapping from $W$ to $V$, so that $U_2 \circ U_1$ is a linear
operator on $V$ and $U_1 \circ U_2$ is a linear operator on $W$, then
$\tr_V U_2 \circ U_1 = \tr_W U_1 \circ U_2$.  This is a standard
property of the trace, and it implies the second equality in (\ref{def
of lange R_1, R_2 rangle_H}).  To see that $\langle R_1, R_2
\rangle_H$ satisfies the positivity condition required of an inner
product, one can compute ${\tr}_V R^* \circ R$ using an orthogonal
basis of $V$, and reduce to the positivity condition for the inner
product on $W$, since 
\begin{equation}
	\langle (R^* \circ R)(v), v \rangle_V = \langle R(v), R(v) \rangle_W
\end{equation}
for all $v$ in $V$.  Of course one could instead use an orthogonal
basis for $W$ and reduce to the positivity condition for the inner
product on $W$.

	Another way to describe the inner product on $H$ is to observe
that if $R_1$, $R_2$ are rank-$1$ operators given by
\begin{equation}
	R_1(v) = \langle v, v_1 \rangle_V \, w_1,
		\quad R_2(v) = \langle v, v_2 \rangle_V \, w_2,
\end{equation}
where $v_1$, $v_2$ are elements of $V$ and $w_1$, $w_2$ are elements
of $W$, then $R_2^*(w) = \langle w, w_2 \rangle_W \, v_2$ and
\begin{equation}
	\langle R_1, R_2 \rangle_H = 
		\langle v_2, v_1 \rangle_V \, \langle w_1, w_2 \rangle_W.
\end{equation}
If $v_1, \ldots, v_n$ is an orthogonal basis for $V$, and $w_1,
\ldots, w_m$ is an orthogonal basis for $W$, then
\begin{equation}
	\langle \cdot, v_j \rangle_V \, w_l, 
		\quad 1 \le j \le n, \ 1 \le l \le m
\end{equation}
is an orthogonal basis for $H$.

	If $T$ is a linear operator on $V$, then $R \mapsto R \circ T$
defines a linear operator on $H$.  One can verify that the adjoint of
this linear operator, with respect to the inner product just defined
on $H$, is given by $R \mapsto R \circ T^*$.  Similarly, if $S$ is
a linear operator on $W$, then $R \mapsto S \circ R$ defines a linear
operator on $H$, and the adjoint of this operator is given by 
$R \mapsto S^* \circ R$.  

	Thus, if $\mathcal{A}$ is a $*$-algebra of operators on $V$,
then $\mathcal{A}_{H,1}$ is a $*$-algebra of operators on $H$, and if
$\mathcal{B}$ is a $*$-algebra of operators on $W$, then
$\mathcal{B}_{H,2}$ is a $*$-algebra of operators on $H$.  As in
Remark \ref{relation with ``expanding an algebra''}, this also came up
in Subsection \ref{Expanding an algebra}.

}
\end{remark}

	Let $\mathcal{A}$, $\mathcal{B}$ be algebras of operators on
$V$, $W$, respectively.  The \emph{combined algebra of
operators}\index{combined algebra of operators} $\mathcal{C}$ on $H$
is defined to be the algebra generated by $\mathcal{A}_{H,1}$ and
$\mathcal{B}_{H,2}$.  Thus the elements of $\mathcal{C}$ are the
operators on $H$ which can be written as
\begin{equation}
\label{A_1 B_1 + A_2 B_2 + cdots + A_r B_r}
	A_1 \, B_1 + A_2 \, B_2 + \cdots + A_r \, B_r,
\end{equation}
where each $A_j$ lies in $\mathcal{A}_{H,1}$ and each $B_j$ lies in
$\mathcal{B}_{H,2}$.  Because the elements of $\mathcal{A}_{H,1}$ and
$\mathcal{B}_{H,2}$ commute with each other, one does not need more
complicated products in (\ref{A_1 B_1 + A_2 B_2 + cdots + A_r B_r}).

\beginremark
\label{mathcal{A} = mathcal{L}(V), mathcal{B} = mathcal{L}(W), 2}
{\rm If $\mathcal{A} = \mathcal{L}(V)$ and $\mathcal{B} =
\mathcal{L}(W)$, then the combined algebra $\mathcal{C}$ is all of
$\mathcal{L}(H)$.  This is not hard to check.  
}
\end{remark}

\beginremark
\label{combined algebra, inner products, *-algebras}
{\rm Suppose that $k$ is a symmetric field and that $V$ and $W$ are
equipped with inner products, as in Remark \ref{k a symmetric field,
inner products}.  If $\mathcal{A}$ and $\mathcal{B}$ are $*$-algebras
of operators on $V$ and $W$, then the combined algebra $\mathcal{C}$
is a $*$-algebra on $H$ (using the inner product on $H$ defined in Remark
\ref{k a symmetric field, inner products}).  
}
\end{remark}

\beginremark
\label{combined algebra and group representations}
{\rm Suppose that $G_1$, $G_2$ are finite groups, and that $\sigma$,
$\tau$ are representations of $G_1$, $G_2$ on $V$, $W$, respectively.
Consider the product group $G_1 \times G_2$, in which the group
operation is defined componentwise, using the group operations on
$G_1$, $G_2$.  Let $\rho$ be the representation on $H$ obtained from
$\sigma$, $\tau$ by setting
\begin{equation}
	\rho_{(x,y)}(R) = \tau_y \circ R \circ (\sigma_x)^{-1}
\end{equation}
for all $(x,y)$ in $G_1 \times G_2$ and $R$ in $H$.  If $\mathcal{A}$
is the algebra of operators on $V$ generated by $\sigma$, and
$\mathcal{B}$ is the algebra of operators on $W$ generated by $\tau$,
then the combined algebra $\mathcal{C}$ is the same as the algebra of
operators on $H$ generated by the representation $\rho$ of $G_1 \times
G_2$.
}
\end{remark}

\beginlemma
\label{dimension of combined algebra}
If $\mathcal{A}$, $\mathcal{B}$ are algebras of operators on
$V$, $W$ with dimensions $r$, $s$, respectively, as vector spaces over
$k$, then the combined algebra $\mathcal{C}$ has dimension $r s$.
\end{lemma}

	To see this, suppose that $A_1, \ldots, A_r$ is a basis for
$\mathcal{A}$, and that $B_1, \ldots, B_s$ is a basis for
$\mathcal{B}$.  Consider the operators 
\begin{equation}
	R \mapsto B_j \circ R \circ A_l, \quad 1 \le j \le s, \ 1 \le l \le r,
\end{equation}
on $H$.  It is easy to check that the combined algebra $\mathcal{C}$
is spanned by these $r s$ operators, and we would like to show that
they are linearly independent.

	Let $c_{j,l}$, $1 \le j \le s$, $1 \le l \le r$, be an family
of elements of $k$ such that the operator
\begin{equation}
	R \mapsto \sum_{j=1}^s \sum_{l=1}^r c_{j,l} \, B_j \circ R \circ A_l
\end{equation}
on $H$ is equal to $0$, i.e.,
\begin{equation}
\label{sum_{j=1}^s sum_{l=1}^r c_{j,l} B_j circ R circ A_l = 0}
	\sum_{j=1}^s \sum_{l=1}^r c_{j,l} \, B_j \circ R \circ A_l = 0
\end{equation}
for all $R$ in $H$.  If $w$ is any element of $W$ and $f$ is any
linear mapping from $V$ to $k$, then $R(v) = f(v) \, w$ is a linear
mapping from $V$ to $W$, and (\ref{sum_{j=1}^s sum_{l=1}^r c_{j,l} B_j
circ R circ A_l = 0}) can be rewritten for this choice of $R$ as
\begin{equation}
	\sum_{j=1}^s \sum_{l=1}^r c_{j,l} \, B_j(w) \, f(A_l(v)) = 0
\end{equation}
for all $v$ in $V$.  Let us rewrite this again as
\begin{equation}
	\sum_{j=1}^s \Bigl(\sum_{l=1}^r c_{j,l} \, f(A_l(v)) \Bigr) B_j(w) = 0.
\end{equation}
This holds for all $w$ in $W$, $v$ in $V$, and linear mappings $f$
from $V$ to $k$.  For any fixed $v$ and $f$, the linear independence of
the $B_j$'s as operators on $W$ leads to
\begin{equation}
	\sum_{l=1}^r c_{j,l} \, f(A_l(v)) = 0
\end{equation}
for each $j$.  The linear independence of the $A_l$'s now implies that
$c_{j,l} = 0$ for all $j$ and $l$, which is what we wanted.

\section{Division algebras of operators, 2}
\label{Division algebras of operators, 2}
\setcounter{equation}{0}

	Let $k$ be a field, let $V$ be a vector space, and let
$\mathcal{A}$ be an algebra of operators on $V$ which is a
division algebra.  Put
\begin{equation}
	\Gamma = \mathcal{A} \cap \mathcal{A}',
\end{equation}
so that $\Gamma$ is the \emph{center} of $\mathcal{A}$, i.e., the
collection of elements of $\mathcal{A}$ which commute with all other
elements of $\mathcal{A}$.  Thus $\Gamma$ is an algebra of operators
on $V$ which is commutative and a division algebra of operators.

	It will be convenient to write $Y$ for $\mathcal{A}$ viewed as
a vector space over $k$, for the purpose of considering linear
operators on $Y$.  Define $\mathcal{A}_1$ to be the algebra of linear
operators on $Y$ of the form
\begin{equation}
	R \mapsto R \circ T, \quad T \in \mathcal{A},
\end{equation}
and define $\mathcal{A}_2$ to be the algebra of linear operators
on $Y$ of the form
\begin{equation}
	R \mapsto T \circ R, \quad T \in \mathcal{A}.
\end{equation}
When $T$ lies in $\Gamma$, $R \circ T = T \circ R$ for all $R$ in $Y$,
and we write $\Gamma_0$ for the algebra of linear operators on $Y$
of the form
\begin{equation}
	R \mapsto R \circ T = T \circ R, \quad T \in \Gamma.
\end{equation}

\beginlemma
\label{properties of mathcal{A}_1, mathcal{A}_2}
{\rm (a)} $\mathcal{A}_1$, $\mathcal{A}_2$ are division algebras
of operators.

{\rm (b)} Every element of $\mathcal{A}_1$ commutes with every element
of $\mathcal{A}_2$.

{\rm (c)} $\mathcal{A}_1$, $\mathcal{A}_2$ are irreducible algebras
of operators on $Y$.
\end{lemma}

	These properties are easy to verify, just from the definitions
and the assumption that $\mathcal{A}$ is a division algebra of
operators.  Concerning (c), it can be helpful to rephrase the question
as follows: if $R$, $\widetilde{R}$ are elements of $Y$ and $R \ne 0$,
then there are elements of $\mathcal{A}_1$, $\mathcal{A}_2$ which take
$R$ to $\widetilde{R}$.

\beginlemma
\label{Gamma_0 = mathcal{A}_1 cap mathcal{A}_2}
$\Gamma_0 = \mathcal{A}_1 \cap \mathcal{A}_2$.
\end{lemma}

	Indeed, suppose that $S$, $T$ are elements of $\mathcal{A}$
such that the operators
\begin{equation}
	R \mapsto R \circ S, \quad R \mapsto T \circ R
\end{equation}
on $Y$ are the same.  We can take $R$ to be the identity operator
on $V$ to obtain that $S = T$, and then our hypothesis becomes
$R \circ T = T \circ R$ for all $R$ in $\mathcal{A}$.  This says
exactly that $T$ lies in $\Gamma$, as desired.

\beginlemma
\label{(mathcal{A}_1)' = mathcal{A}_2, (mathcal{A}_2)' = mathcal{A}_1}
$(\mathcal{A}_1)' = \mathcal{A}_2$ and $(\mathcal{A}_2)' =
\mathcal{A}_1$, where the primes indicate that we take the commutant
of the given algebra as an algebra of operators on $Y$.
\end{lemma}

	The fact that $\mathcal{A}_2 \subseteq (\mathcal{A}_1)'$
and $\mathcal{A}_1 \subseteq (\mathcal{A}_2)'$ is the same as
part (b) of Lemma \ref{properties of mathcal{A}_1, mathcal{A}_2}.
Conversely, suppose that $\Phi$ is a linear transformation on $Y$
which commutes with every element of $\mathcal{A}_1$.  This is
the same as saying that
\begin{equation}
	\Phi(R \circ S) = \Phi(R) \circ S
\end{equation}
for all $R$, $S$ in $Y$, i.e., in $\mathcal{A}$.  Applying this to
$R$ equal to the identity operator $I$ on $V$ we get that
\begin{equation}
	\Phi(S) = \Phi(I) \circ S
\end{equation}
for all $S$ in $Y$.  This says exactly that $\Phi$ lies in
$\mathcal{A}_2$, so that $(\mathcal{A}_1)' \subseteq \mathcal{A}_2$.
Similarly, $(\mathcal{A}_2)' \subseteq \mathcal{A}_1$, and the lemma
follows.

	Let $\mathcal{A}_{12}$ be the algebra of operators on $Y$
which is generated by $\mathcal{A}_1$, $\mathcal{A}_2$.

\beginlemma
\label{Gamma_0 is the center of mathcal{A}_{12}}
$\Gamma_0$ is the center of $\mathcal{A}_{12}$.
\end{lemma}

	This is easy to check.

\beginproposition
\label{mathcal{A}_{12} = Gamma_0'}
$\mathcal{A}_{12} = \Gamma_0'$.
\end{proposition}

	To prove this we use Lemma \ref{criterion for mathcal{A} =
mathcal{B}', 2}, where now $Y$ has the role that $V$ had before,
$\mathcal{A}_{12}$ has the role that $\mathcal{A}$ had before, and
$\Gamma_0$ has the role that $\mathcal{B}$ had before.  Of course
$\mathcal{A}_{12} \subseteq \Gamma_0'$, which was a standing
assumption for Lemma \ref{criterion for mathcal{A} = mathcal{B}', 2},
mentioned just before the statement of Lemma \ref{criterion for
mathcal{A} = mathcal{B}', 0}.  According to \ref{criterion for
mathcal{A} = mathcal{B}', 2}, it suffices to show that
$\mathcal{A}_{12}$ is irreducible, and that for any elements $R_1$,
$R_2$, $U_1$, $U_2$ of $Y$ such that $R_1$, $R_2$ are
$\Gamma_0$-independent there is a $\Phi$ in $\mathcal{A}_{12}$ such
that $\Phi(R_1) = U_1$ and $\Phi(R_2) = U_2$.  (The notion of
independence employed here was defined near the beginning of Section
\ref{Division algebras of operators}.)  The irreducibility of
$\mathcal{A}_{12}$ is a consequence of the irreducibility of either
$\mathcal{A}_1$ or $\mathcal{A}_2$, and so it remains to verify the
second condition.

	For the second condition, it is enough to show that if $R$,
$U_1$, $U_2$ are elements of $Y$ such that $R$ is not an element of
$\Gamma \subseteq \mathcal{A}$, then there is a $\Psi$ in
$\mathcal{A}_{12}$ such that $\Psi(I) = U_1$, $\Psi(R) = U_2$.  Here
$I$ denotes the identity operator on $V$, as an element of $Y =
\mathcal{A}$, just as $\Gamma$ can be viewed as a subspace of $Y$.  In
other words, if $R_1$, $R_2$ are elements of $Y$ which are
$\Gamma_0$-independent, then they can be mapped to $I$, $R_1^{-1} \,
R_2$ by an element of $\mathcal{A}_2 \subseteq \mathcal{A}_{12}$,
and the $\Gamma_0$-independence of $R_1$, $R_2$ says exactly that
$R_1^{-1} \, R_2$ does not lie in $\Gamma$.

	In fact it is enough to show that for every $R$, $U$ in $Y$
such that $R$ is not in $\Gamma$ there is a $\Psi$ in
$\mathcal{A}_{12}$ which satisfies $\Psi(I) = 0$ and $\Psi(R) = U$.
The reason for this is that for any $U_1$ in $Y$ there is an element
of $\mathcal{A}_1$ (or $\mathcal{A}_2$) which sends $I$ to $U_1$
(corresponding to composition with $U_1$).  To prescribe a value also
at $R$, one can use a $\Psi$ in $\mathcal{A}_{12}$ as in the new
version of the condition.

	Thus we let $R$, $U$ in $Y$ be given, with $R$ not in $\Gamma$.
Because $R$ is not in $\Gamma$, there is an $S$ in $Y = \mathcal{A}$
such that $R \circ S - S \circ R \ne 0$.  Thus the inverse of
$R \circ S - S \circ R$ exists and lies in $\mathcal{A}$.
Let $\Psi$ be the linear
operator on $Y$ defined by
\begin{equation}
	\Psi(A) =
    U \circ (R \circ S - S \circ R)^{-1} \circ (A \circ S - S \circ A),
\end{equation}
$A \in Y = \mathcal{A}$.  It is easy to see that $\Psi$ lies in
$\mathcal{A}_{12}$, i.e., it can be written as a sum of compositions
of operators on $Y$ in $\mathcal{A}_1$, $\mathcal{A}_2$.  Clearly
$\Psi(I) = 0$ and $\Psi(R) = U$, as desired.  This proves Proposition
\ref{mathcal{A}_{12} = Gamma_0'}.

\begincorollary
\label{dimension of mathcal{A}_{12}}
The dimension of $\mathcal{A}_{12}$ is equal to the square of the
dimension of $\mathcal{A} = Y$ divided by the dimension of $\Gamma$,
as vector spaces over $k$.
\end{corollary}

	Compare with Lemma \ref{dimension of mathcal{B}' as a vector
space}, and note that the dimension of $\Gamma$ is equal to the
dimension of $\Gamma_0$.

	Let us now relate this to the set-up in Section \ref{Vector
spaces of linear mappings}.  What were the vector spaces $V$, $W$
before are now both taken to be $V$.  Thus we write $H$ for the vector
space of linear mappings from $V$ to itself.  The algebra of operators
$\mathcal{A}$ on $V$ is now used for both $\mathcal{A}$ and
$\mathcal{B}$ in Section \ref{Vector spaces of linear mappings}.  Let
$\mathcal{A}_{H,1}$, $\mathcal{A}_{H,2}$, and the combined algebra
$\mathcal{C}$ be the algebras of operators defined on $H$ in Section
\ref{Vector spaces of linear mappings} (with $\mathcal{B} =
\mathcal{A}$).

	We can view $Y = \mathcal{A}$ as a vector subspace of $H$.  As
such, it is an invariant subspace for $\mathcal{A}_{H,1}$,
$\mathcal{A}_{H,2}$, and the combined algebra $\mathcal{C}$.  As in
Notation \ref{notation about mathcal{A}(U)}, $\mathcal{A}_{H,1}(Y)$,
$\mathcal{A}_{H,2}(Y)$, and $\mathcal{C}(Y)$ denote the algebras of
operators on $Y$ which occur as restrictions of operators in
$\mathcal{A}_{H,1}$, $\mathcal{A}_{H,2}$, and $\mathcal{C}$ to the
invariant subspace $Y$.  It is easy to check that
\begin{equation}
	\mathcal{A}_1 = \mathcal{A}_{H,1}(Y), 
		\ \mathcal{A}_2 = \mathcal{A}_{H,2}(Y),
		\ \mathcal{A}_{12} = \mathcal{C}(Y).
\end{equation}

	Just as for $\mathcal{A}$, we can associate to $\Gamma$ the
algebras of linear operators $\Gamma_{H,1}$, $\Gamma_{H,2}$ on $H$ of
the form
\begin{equation}
\label{R mapsto R circ T, R mapsto T circ R (T in Gamma)}
	R \mapsto R \circ T, \quad R \mapsto T \circ R,
\end{equation}
respectively, where $T$ lies in $\Gamma$.  Thus 
\begin{equation}
	\Gamma_{H,1} \subseteq \mathcal{A}_{H,1}, \quad
		\Gamma_{H,2} \subseteq \mathcal{A}_{H,2},
\end{equation}
and $Y$ is an invariant subspace for $\Gamma_{H,1}$, $\Gamma_{H,2}$.
For $T$ in $\Gamma$, the linear operators in (\ref{R mapsto R circ T,
R mapsto T circ R (T in Gamma)}) are different as operators on $H$
unless $T$ is a scalar multiple of the identity, but the restrictions
of these operators to $Y$ are the same, since $\Gamma$ is the center
of $\mathcal{A}$.  In particular,
\begin{equation}
	\Gamma_{H,1}(Y) = \Gamma_{H,2}(Y) = \Gamma_0.
\end{equation}

	Let $T_1, \ldots, T_\ell$ be a collection of operators in
$\mathcal{A}$ which are $\Gamma_0$-independent and whose
$\Gamma_0$-span is all of $\mathcal{A}$, as discussed in Section
\ref{Division algebras of operators}, where we think of $\Gamma_0
\simeq \Gamma$ as a division algebra of operators on $\mathcal{A}$.
In fact, we can view $\mathcal{A}$ as a vector space over the field
$\Gamma$, since $\Gamma$ is commutative and commutes with all elements
of $\mathcal{A}$.  In these terms $T_1, \ldots, T_\ell$ is a basis for
$\mathcal{A}$ as a the vector space over $\Gamma$.  If $\gamma_1,
\ldots, \gamma_m$ is a basis for $\Gamma$ as a vector space over $k$,
then the family of products $\gamma_i \circ T_j$, $1 \le i \le m$, $1
\le j \le \ell$, is a basis for $\mathcal{A}$ as a vector space over
$k$.

	As in the proof of Lemma \ref{dimension of combined algebra},
the operators on $H$ given by
\begin{equation}
 R \mapsto \gamma_{i_2} \circ T_{j_2} \circ R \circ \gamma_{i_1} \circ T_{j_1},
		\quad 1 \le i_1, i_2 \le m, \ 1 \le j_1, j_2 \le \ell,
\end{equation}
form a basis for the combined algebra $\mathcal{C}$.  When one
restricts to $Y$, as with $\mathcal{C}(Y)$, then it is enough to
consider the operators of the form
\begin{equation}
 R \mapsto \gamma_i \circ T_{j_2} \circ R \circ T_{j_1},
		\quad 1 \le i \le m, \ 1 \le j_1, j_2 \le \ell,
\end{equation}
since the $\gamma_i$'s commute with $R \in Y$ and the $T_j$'s.  In
other words, on $Y$, the span of these operators is the same as the
span of the previous ones.  This fits with Corollary \ref{dimension of
mathcal{A}_{12}}, which implies that $m \, \ell^2$ is equal to the
dimension of $\mathcal{C}(Y) = \mathcal{A}_{12}$.

\section{Some basic situations}
\label{Some basic situations}
\setcounter{equation}{0}

\subsection{Some algebras of operators}
\label{Some algebras of operators}

	Let $k$ be a field, and let $V$ be a vector space over $k$.
Suppose that $V_1, \ldots, V_\ell$ are vector subspaces of $V$ such
that
\begin{equation}
\label{{0} subseteq V_1 subseteq cdots subseteq V_ell subseteq V}
	\{0\} \subseteq V_1 \subseteq \cdots \subseteq V_\ell \subseteq V,
\end{equation}
where each inclusion is strict.  Consider the algebra of operators
$\mathcal{A}$ on $V$ consisting of the linear mappings $T : V \to V$
such that
\begin{equation}
\label{T(V_j) subseteq V_j, 1 le j le ell}
	T(V_j) \subseteq V_j, \quad 1 \le j \le \ell.
\end{equation}

\beginlemma
\label{lots of subspaces as images and kernels in this situation}
Under the conditions above, if $Z$ and $W$ are arbitrary vector
subspaces of $V$, then there are linear operators $T_1$, $T_2$ in
$\mathcal{A}$ such that $T_1(V) = W$ and the kernel of $T_2$ is $Z$.
\end{lemma}

	To be more precise, one can choose $T_1$ so that it
satisfies
\begin{equation}
	T_1(V_j) = W \cap V_j, \quad 1 \le j \le \ell,
\end{equation}
in addition to $T_1(V) = W$.  This is not hard to arrange, using a
suitable basis for $V$.  Notice that this is not the only possibility,
however; for instance, if the dimension of $W$ is less than or equal
to the codimension of $V_\ell$ in $V$, then one can choose $T_1$ so
that $T_1$ is equal to $0$ on $V_\ell$, but maps a subspace of $V$
whose intersection with $V_\ell$ is trivial onto $W$.  As for $T_2$,
we do need to have
\begin{equation}
	\{v \in V_j : T_2(v) = 0\} = Z \cap V_j, \quad 1 \le j \le \ell,
\end{equation}
in order for the kernel of $T_2$ to be equal to $Z$, and again
this is not hard to arrange, through suitable bases.

\beginproposition
\label{mathcal{A}' = {lambda I : lambda in k} in this situation}
Under the conditions above, the commutant $\mathcal{A}'$ of
$\mathcal{A}$ consists of scalar multiples of the identity operator on
$V$.
\end{proposition}

	Of course $\mathcal{A}'$ always contains the scalar multiples
of the identity.  Conversely, suppose that $S$ lies in $\mathcal{A}'$.
Because of Lemmas \ref{lots of subspaces as images and kernels in this
situation} and \ref{invariant subspaces from kernels and images}, if
$U$ is any vector subspace of $V$, then $S(U) \subseteq U$.  (For this
one might as well even consider only $1$-dimensional subspaces $U$ of
$V$.)  From this it is not hard to see that $S$ must be a scalar multiple
of the identity operator.

\subsection{A construction}
\label{A construction}

	Let $k$ be a field, let $V$ be a vector space, and let $V^*$
denote the \emph{dual space} of $V$, i.e., the vector space of linear
mappings from $V$ into $k$.  If $T$ is a linear mapping from $V$ to
itself, then there is an associated dual linear operator
$\widetilde{T}$ on $V^*$, which is defined by saying that if $\phi$ is
a linear functional on $V$, then $\widetilde{T}(\phi)$ is the linear
functional given by $\widetilde{T}(\phi)(v) = \phi(T(v))$.  The
transformation from a linear operator $T$ on $V$ to the dual linear
operator $\widetilde{T}$ is linear in $T$, and satisfies
$\widetilde{(T_1 \circ T_2)} = \widetilde{T_2} \circ \widetilde{T_1}$
for all linear operators $T_1$, $T_2$ on $V$.  The dual of the
identity operator on $V$ is equal to the identity operator on $V^*$,
and the dual of an invertible operator on $V$ is an invertible
operator on $V^*$, with the inverse of the dual being the dual of the
inverse.  There is a natural identification of $V$ with the second
dual space $V^{**}$, and the dual of the dual of an operator $T$ on
$V$ coincides with $T$ under this identification.

	Assume now that $\mathcal{A}$ is an algebra of linear operators
on $V$.  We can define $\widetilde{\mathcal{A}}$ to be the algebra of
linear operators on the dual space $V^*$ consisting of the duals of the
linear operators in $\mathcal{A}$.  

	Let $W$ be the vector space which is the direct sum of $V$ and
$V^*$.  Thus $W$ consists of ordered pairs $(u, \phi)$ with $u$ in $V$
and $\phi$ in $V^*$, with addition and scalar multiplication of
vectors defined coordinatewise.  Define $\widehat{\mathcal{A}}$ to be
the collection of operators on $W$ of the form
\begin{equation}
\label{(u, phi) mapsto (S(u), widetilde{T}(phi))}
	(u, \phi) \mapsto (S(u), \widetilde{T}(\phi)),
\end{equation}
where $S$, $T$ lie in $\mathcal{A}$.  It is not hard to see that
$\widehat{\mathcal{A}}$ is an algebra of operators on $W$.  The
subspaces $V \times \{0\}$ and $\{0\} \times V^*$ of $W$ are
complementary to each other and invariant under
$\widehat{\mathcal{A}}$, and, in effect, the restriction of
$\widehat{\mathcal{A}}$ to these subspaces reduces to $\mathcal{A}$
and $\widetilde{\mathcal{A}}$, respectively.

	Define bilinear forms $B_1$, $B_2$ on $W$ by
\begin{equation}
	B_1\Big((u,\phi), (v,\psi)\Big) = \psi(u) + \phi(v)
\end{equation}
and
\begin{equation}
	B_2\Big((u,\phi), (v,\psi)\Big) = \psi(u) - \phi(v).
\end{equation}
Notice that $B_1$ is a symmetric bilinear form while $B_2$ is
antisymmetric, i.e.,
\begin{eqnarray}
    B_1\Big((u,\phi), (v,\psi)\Big) & = & B_1\Big((v,\psi), (u,\phi)\Big), 
									\\
    B_2\Big((u,\phi), (v,\psi)\Big) & = & - B_2\Big((v,\psi), (u,\phi)\Big).
								\nonumber
\end{eqnarray}
Both are nondegenerate, which is to say that for $i = 1$ or $2$ and
for any nonzero element $(u,\phi)$ of $W$ there is a $(v,\psi)$ in $W$
such that 
\begin{equation}
	B_i\Big((u,\phi), (v,\psi)\Big) \ne 0.
\end{equation}
On the other hand, it is not true that $B_1$ is definite, because there
are plenty of nonzero elements $(u,\phi)$ of $W$ which satisfy
\begin{equation}
	B_1\Big((u,\phi), (u,\phi)\Big) = 0.
\end{equation}
For $B_2$ we have that
\begin{equation}
	B_2\Big((u,\phi), (u,\phi)\Big) = 0
\end{equation}
for all $(u,\phi)$ in $W$, which is in fact equivalent to the antisymmetry
of $B_2$.

	Suppose that $S$, $T$ are elements of $\mathcal{A}$, and that
$R$ is the element of $\widehat{\mathcal{A}}$ given by $R(u,\phi) =
(S(u), \widetilde{T}(\phi))$, as in (\ref{(u, phi) mapsto (S(u),
widetilde{T}(phi))}).  Observe that
\begin{equation}
    \quad
	B_1\Big(R(u,\phi), (v,\psi)\Big) = \psi(S(u)) + \phi(T(v))
		= B_1\Big((u,\phi), R^t(v,\psi)\Big),
\end{equation}
where $R^t$ is the element of $\widehat{\mathcal{A}}$ defined by
\begin{equation}
	R^t(v,\psi) = (T(v),\widetilde{S}(\psi)).
\end{equation}
Similarly,
\begin{equation}
    \quad
	B_2\Big(R(u,\phi), (v,\psi)\Big) = \psi(S(u)) - \phi(T(v))
		= B_2\Big((u,\phi), R^t(v,\psi)\Big).
\end{equation}

	To put it another way, $R^t$ is the ``transpose'' of $R$ as a
linear mapping from $W$ to itself with respect to each of the bilinear
forms $B_1$ and $B_2$.  Thus we see that $\widehat{\mathcal{A}}$
contains the transposes of all of its elements, which is to say that it
is a ``$t$-algebra'' of operators on $W$, in a sense analogous to that of
Definition \ref{def of t-algebras}, considered before in the setting
of definite scalar products.

	Let us specialize to the case where $\mathcal{A}$ is as in
Subsection \ref{Some algebras of operators}, corresponding to the
subspaces $V_1, \ldots, V_\ell$ of $V$.  For any vector subspace $U$
of $V$, we denote by $U^\perp$ the subspace of $V^*$ consisting of
linear functionals $\phi$ on $V$ such that $\phi(u) = 0$ for all $u
\in U$.  A standard observation from linear algebra is that a vector
$v$ in $V$ lies in $U$ if and only if $\phi(v) = 0$ for all $\phi$ in
$U^\perp$.  From (\ref{{0} subseteq V_1 subseteq cdots subseteq V_ell
subseteq V}) we get that
\begin{equation}
	\{0\} \subseteq V_\ell^\perp \subseteq \cdots 
			\subseteq V_1^\perp \subseteq V^*,
\end{equation}
with the inclusions again being strict since they were strict before.
One can check that a linear operator on $V^*$ lies in
$\widetilde{\mathcal{A}}$ in these circumstances if and only if
each subspace $V_j^\perp$, $1 \le j \le \ell$, of $V^*$ is invariant
under the operator.

	One can also verify that the commutant of
$\widehat{\mathcal{A}}$ in $\mathcal{L}(W)$ consists of linear
combinations of the coordinate projections from $W \simeq V \times
V^*$ onto $V \times \{0\}$ and $\{0\} \times V^*$.

\section{Positive characteristic}
\label{Positive characteristic}
\setcounter{equation}{0}

	Throughout this section, we let $p$ be a prime number, and
we let $k$ be a field of characteristic $p$.

\subsection{A few basic facts}
\label{A few basic facts (section on pos. char.)}

	The most basic example of a field with characteristic $p$ is
${\bf Z}/p {\bf Z}$, the field with $p$ elements obtained by taking
the integers ${\bf Z}$ modulo $p$.  Any field with characteristic $p$
contains a copy of this field in a canonical way, as the set of
elements obtained by taking sums of the multiplicative identity
element.  In particular, any such field can be viewed as a vector
space over ${\bf Z}/p {\bf Z}$, and if a field of characteristic $p$
has a finite number of elements, then that number is a power of $p$.
A well-known result is that for each positive integer $m$ there is a
field of characteristic $p$ with $p^m$ elements, and this field is
unique up to isomorphism.

	In a field $k$ of characteristic $p$ we have that
\begin{equation}
\label{(x + y)^p = x^p + y^p}
	(x + y)^p = x^p + y^p
\end{equation}
for all $x$, $y$ in $k$.  Indeed,
\begin{equation}
	(x + y)^p = \sum_{j=0}^p {p \choose j} x^j \, y^{p-j} = x^p + y^p,
\end{equation}
where the first step is an instance of the binomial theorem and the
second step uses the observation that the binomial coefficient ${p
\choose j}$ is a positive integer which is divisible by $p$ when $1
\le j \le p-1$, and hence gives $0$ in the field $k$.

	Notice that
\begin{equation}
\label{(x - y)^p = x^p - y^p}
	(x - y)^p = x^p - y^p
\end{equation}
for all $x$, $y$ in $k$.  This follows from (\ref{(x + y)^p = x^p +
y^p}) and the fact that $(-1)^p = -1$ in $k$.  More precisely, $(-1)^2
= 1 = -1$ when $p = 2$, while $(-1)^p = -1$ when $p$ is odd because
$(-1)^a = 1$ when $a$ is an even integer, such as $p-1$.  

	By iterating these identities one obtains
\begin{equation}
\label{(x + y)^{p^l} = x^{p^l} + y^{p^l}}
	(x + y)^{p^l} = x^{p^l} + y^{p^l}
\end{equation}
and
\begin{equation}
\label{(x - y)^{p^l} = x^{p^l} - y^{p^l}}
	(x - y)^{p^l} = x^{p^l} - y^{p^l}
\end{equation}
for all positive integers $l$ and all $x$, $y$ in $k$.  As a consequence
of the latter we obtain the following.

\beginlemma
\label{a {p^l}th root of unity is = 1}
If $l$ is a positive integer and $x$ is an element of $k$ which
satisfies $x^{p^l} = 1$, then $x = 1$.
\end{lemma}

	By contrast, there can be nontrivial roots of unity of other
orders.  If $k$ is finite, with $p^m$ elements for some positive
integer $m$, then the set of nonzero elements of $k$ is a group under
multiplication of order $p^m - 1$.  It is well known that this group
is in fact \emph{cyclic}.  The identity
\begin{equation}
	x^{p^m - 1} = 1
\end{equation}
for all nonzero elements $x$ of $k$ can be derived more simply, from
the general result that the order of a group is divisible by the order
of any subgroup, and hence by the order of any element of the group.
Of course $p^m - 1$ is not divisible by $p$.  Since $p^m - 1$ is the
number of nonzero elements in $k$ when $k$ has $p^m$ elements, we see
that the polynomial $x^{p^m - 1} - 1$ completely factors in $k$, i.e.,
it is the product of the $p^m - 1$ linear polynomials $x - \alpha$,
where $\alpha$ runs through the set of nonzero elements of $k$.
Alternatively, one can say that the polynomial $x^{p^m} - x = x
(x^{p^m - 1} - 1)$ completely factors, as the product of the 
linear polynomials $x - \alpha$, where now $\alpha$ runs through all
elements of $k$.

\subsection{Representations}
\label{Representations (finite groups, k of char. p)}

	Let $V$ be a vector space over $k$, let $G$ be a finite group,
and let $\rho$ be a representation of $G$ on $V$.  If $p$ divides the
order of $G$, then the argument described in Subsection
\ref{Reducibility, continued} for finding invariant complements of
invariant subspaces of a representation does not work, and indeed
the result is not true.

	If $W$ is a subspace of $V$ which is invariant under $\rho$,
then one can at least consider the quotient vector space $V / W$, and
the representation of $G$ on $V / W$ obtained from $\rho$ in the
obvious manner.  Recall that $V / W$ is defined using the equivalence
relation $\sim$ on $V$ given by
\begin{equation}
	\hbox{$v_1 \sim v_2$ if and only if $v_1 - v_2 \in W$},
\end{equation}
by passing to the corresponding equivalence classes.  The assumption
that $W$ is invariant under $\rho$ implies that this equivalence
relation is preserved by $\rho$, so that $\rho$ leads to a
representation on the quotient space $V / W$.

	Thus if the representation $\rho$ on $V$ is not irreducible,
so that there is an invariant subspace $W$ which is neither the zero
subspace nor $V$, then we get two representations of smaller positive
degree, namely the restriction of $\rho$ to $W$ and the representation
on $V / W$ obtained from $\rho$.  The sum of the degrees of these two
new representations is equal to the degree of the original
representation, i.e., the sum of the dimensions of $W$ and $V / W$ is
equal to the dimension of $V$.  We can repeat the process, of passing
to subspaces and quotients, to get a family of irreducible
representations, where the sum of the degrees of these representations
is equal to the degree of $\rho$ (the dimension of $V$).  However,
because we do not necessarily have direct sum decompositions at each
step, the family of irreducible representations may not contain as
much information as the original representation.

	Let $F(G)$ denote the vector space of $k$-valued functions on
$G$.  As in Subsection \ref{Representations of finite groups}, one has
the left and right regular representations $L$, $R$, of $G$ on $F(G)$,
and these two representations are isomorphic via the simple
correspondence $f(a) \leftrightarrow f(a^{-1})$ for functions on $G$.
Note that the regular representations have the feature that elements
of $G$ different from the identity element correspond to linear
transformations on $F(G)$ which are different from the identity
transformation.  In other words, these representations give rise
to isomorphisms between $G$ and subgroups of the group of all
invertible linear transformations on $F(G)$.

\beginlemma
\label{irred. rep.'s isomorphic to a subrep. of the reg. rep., 2}
Let $Z$ be a vector space over $k$, $\sigma$ a representation of $G$
on $Z$, and $\lambda$ a nonzero linear mapping from $Z$ to $k$.  For
each $v$ in $Z$, define $f_v(y)$ on $G$ by
\begin{equation}
	f_v(y) = \lambda(\sigma_{y^{-1}}(v)),
\end{equation}
and put
\begin{equation}
	U = \{f_v(y) : v \in Z\}.
\end{equation}
The mapping 
\begin{equation}
\label{v mapsto f_v}
	v \mapsto f_v
\end{equation}
is a nonzero linear mapping from $Z$ onto $U$, and $U$ is a nonzero
vector subspace of $F(G)$.  This mapping intertwines the
representations $\sigma$ on $Z$ and the restriction of the left
regular representation to $U$, and $U$ is invariant under the left
regular representation.  If $\sigma$ is an irreducible representation
of $G$, then (\ref{v mapsto f_v}) is one-to-one and yields an
isomorphism between $\sigma$ and the restriction of the left regular
representation to $U$.
\end{lemma}

	This is essentially the same as Lemma \ref{irred. rep.'s
isomorphic to a subrep. of the reg. rep.}, and the proof is the
same as before.

\subsection{Linear transformations}
\label{Linear transformations (pos. char.)}

	Let $V$ be a vector space over $k$.  Suppose that $R$ and $S$
are two linear transformations on $V$ which commute.  Once again we have
that
\begin{equation}
\label{(R + S)^p = R^p + S^p}
	(R + S)^p = R^p + S^p.
\end{equation}
This can be established in the same manner as before, by expanding
$(R + S)^p$ using the binomial theorem, and then observing that the
intermediate terms vanish because $k$ has characteristic $p$.  We also
get
\begin{equation}
\label{(R - S)^p = R^p - S^p}
	(R - S)^p = R^p - S^p,
\end{equation}
and, for any positive integer $l$,
\begin{equation}
\label{(R + S)^{p^l} = R^{p^l} + S^{p^l}}
	(R + S)^{p^l} = R^{p^l} + S^{p^l}
\end{equation}
and
\begin{equation}
\label{(R - S)^{p^l} = R^{p^l} - S^{p^l}}
	(R - S)^{p^l} = R^{p^l} - S^{p^l}.
\end{equation}

	Now assume that $T$ is a linear transformation on $V$ such that
$T^{p^l} = I$ for some positive integer $l$, where $I$ denotes the
identity transformation on $V$.  Since $T$ and $I$ automatically
commute, we get that
\begin{equation}
	(T - I)^{p^l} = T^{p^l} - I = 0.
\end{equation}
Thus $T - I$ is nilpotent of order (at most) $p^l$.

	If $n$ is a positive integer and $S$ is a linear transformation
on $V$ such that 
\begin{equation}
	S^n = I,
\end{equation}
then there is a well-known way to factor $S$ in order to separate the
factors of $p$ in $n$ from the rest.  We begin by writing $n$ as
$p^l \, m$, where $l$ is a nonnegative integer and $m$ is a positive
integer which is not divisible by $p$.  Because $p^l$ and $m$ are relatively
prime, it is a standard result that there are integers $a$, $b$ such that
\begin{equation}
	a \, m + b \, p^l = 1.
\end{equation}
Thus we can write 
\begin{equation}
	S = S_1 \, S_2, \quad\hbox{where}\quad 
			S_1 = S^{a \, m}, \ S_2 = S^{b \, p^l},
\end{equation}
and where it is understood that $S^c$ is interpreted as being the
identity transformation when $c = 0$.  With these choices we have that
\begin{equation}
   S_1^{p^l} = I, \ S_2^m = I, \quad\hbox{and}\quad S_1 \, S_2 = S_2 \, S_1.
\end{equation}


\end{document}